\documentclass{article}

\RequirePackage[OT1]{fontenc}
\RequirePackage{amsthm,amsmath,bm}
\RequirePackage[colorlinks,citecolor=blue,urlcolor=blue]{hyperref}

\usepackage{amsmath,amssymb, a4wide}
\usepackage{amssymb,latexsym}
\usepackage{dsfont, color}
\usepackage{verbatim}
\usepackage{cancel}
\usepackage{authblk}
\usepackage{appendix}
\usepackage{graphicx}

\numberwithin{equation}{section}
\theoremstyle{plain}
\newtheorem{theorem}{Theorem}
\newtheorem{corollary}{Corollary}
\newtheorem{lemma}{Lemma}

\newtheorem{remark}{Remark}[section]
\newtheorem{example}{Example}

\newcommand{\R}{\mathbb{R}}

\newcommand{\dd}{\mathrm{d}}
\newcommand{\vc}{\mathrm{vec}}
\newcommand{\vch}{\mathrm{vech}}
\newcommand{\tr}{\mathrm{tr}}

\newcommand{\bbeta}{\bm\beta}
\newcommand{\beps}{\bm\epsilon}
\newcommand{\btheta}{\bm\theta}
\newcommand{\bxi}{\bm\xi}
\newcommand{\balpha}{\bm\alpha}

\newcommand{\bgamma}{\bm\gamma}
\newcommand{\bmu}{\bm\mu}
\newcommand{\bSigma}{\mathbf{\Sigma}}

\newcommand{\by}{\mathbf{y}}
\newcommand{\bx}{\mathbf{x}}
\newcommand{\bs}{\mathbf{s}}
\newcommand{\bu}{\mathbf{u}}
\newcommand{\bv}{\mathbf{v}}
\newcommand{\bt}{\mathbf{t}}
\newcommand{\bz}{\mathbf{z}}
\newcommand{\bq}{\mathbf{q}}

\newcommand{\bX}{\mathbf{X}}
\newcommand{\bV}{\mathbf{V}}
\newcommand{\bG}{\mathbf{G}}
\newcommand{\bZ}{\mathbf{Z}}
\newcommand{\bI}{\mathbf{I}}
\newcommand{\bA}{\mathbf{A}}
\newcommand{\bM}{\mathbf{M}}
\newcommand{\bC}{\mathbf{C}}
\newcommand{\bB}{\mathbf{B}}
\newcommand{\bD}{\mathbf{D}}
\newcommand{\bL}{\mathbf{L}}
\newcommand{\bK}{\mathbf{K}}

\newcommand{\bR}{\mathbf{R}}

\newcommand{\bE}{\mathbf{E}}

\newcommand{\bH}{\mathbf{H}}

\newcommand{\E}{\mathbb{E}}

\title{S-estimation in Linear Models with Structured Covariance Matrices~\footnote{This work has been partly supported by the French \textit{Agence Nationale de la Recherche} through the Investments for the Future (Investissements d'Avenir) program, grant ANR-17-EURE-0010.}}
\author[1]{Hendrik Paul Lopuha\"a}
\author[2]{Valerie Gares}
\author[3]{Anne Ruiz-Gazen}
\affil[1]{\emph{Delft University of Technology}}
\affil[2]{\emph{Institut National des Sciences Appliqu\'es de Rennes}}
\affil[3]{\emph{Toulouse School of Economics}}
\date{\today}

\begin{document}
\maketitle
\begin{abstract}
We provide a unified approach to S-estimation in balanced linear models with structured covariance matrices.
Of main interest are S-estimators for linear mixed effects models, but our approach also includes
S-estimators in several other standard multivariate models, such as multiple regression,
multivariate regression, and multivariate location and scatter.
We provide sufficient conditions for the existence of S-functionals and S-estimators,
establish asymptotic properties such as consistency and asymptotic normality,
and derive their robustness properties in terms of breakdown point and influence function.
All the results are obtained for general identifiable covariance structures and are established under mild conditions
on the distribution of the observations, which goes far beyond models with elliptically contoured densities.
Some of our results are new and others are more general than existing ones in the literature.
In this way this manuscript completes and improves results on S-estimation in a wide variety of multivariate models.
We illustrate our results by means of a simulation study and an application to
data from a trial on the treatment of lead-exposed children.
\end{abstract}

\section{Introduction}
\label{sec:introduction}
Linear models are widely used and provide a versatile approach for analyzing correlated responses,
such as longitudinal data, growth data or repeated measurements.
In such models, each subject~$i$, $i=1,\ldots,n$, is observed at $k_i$ occasions,
and the vector of responses $\by_i$ is assumed to arise from the model
\[
\by_i=\bX_i\bbeta+\bu_i,
\]
where $\bX_i$ is the design matrix for the $i$th subject and $\bu_i$ is a vector
whose covariance matrix can be used to model the correlation between the responses.
One possibility is the linear mixed effects model, in which the random effects together with the measurement error
yields a specific covariance structure depending on a vector $\btheta$ consisting of some unknown covariance parameters.
Other covariance structures may arise, for example if the $\bu_i$ are the outcome of a time series,
see e.g.,~\cite{jennrich&schluchter1986} or~\cite{fitzmaurice-laird-ware2011}, for different possible covariance structures.

Maximum likelihood estimation of $\bbeta$ and $\btheta$ has been studied, e.g., in~\cite{hartley&rao1967,rao1972,laird&ware1982},
see also~\cite{fitzmaurice-laird-ware2011,demidenko2013}.
To be resistant against outliers, robust methods have been investigated for linear mixed effects models, e.g.,
in~\cite{pinheiro-liu-wu2001,copt2006high,copt&heritier2007,heritier-cantoni-copt-victoriafeser2009,agostinelli2016composite,chervoneva2014}.
This mostly concerns S-estimators, originally introduced in the multiple regression context by Rousseeuw and Yohai~\cite{rousseeuw-yohai1984}
and extended to multivariate location and scatter in~\cite{davies1987,lopuhaa1989}, to multivariate regression in~\cite{vanaelst&willems2004},
and to linear mixed effects models
in~\cite{copt2006high,heritier-cantoni-copt-victoriafeser2009,chervoneva2014}.
S-estimators are well known smooth versions of the minimum volume ellipsoid estimator~\cite{rousseeuw1985} that
are highly resistant against outliers.
As such, S-estimators have gained popularity as robust estimators,
but they may also serve as initial estimators to further improve the efficiency.
However, the theory about these estimators is far from complete,
even in balanced models where the number of observed responses is the same for all subjects.

In view of this, we provide a unified approach to S-estimation in balanced linear models with structured covariance matrices,
and postpone a unified approach for unbalanced models to a future paper.
The balanced setup is already quite flexible and includes several specific multivariate statistical models.
Of main interest are S-estimators for linear mixed effects models, but our approach also includes
S-estimators in several other standard multivariate models, such as multiple regression,
multivariate regression, and multivariate location and scatter.
We provide sufficient conditions for the existence of S-functionals and S-estimators,
establish their asymptotic properties, such as consistency and asymptotic normality,
and derive their robustness properties in terms of breakdown point and influence function.
All results are obtained for a large class of identifiable covariance structures, and are established under very mild conditions
on the distribution of the observations, which goes far beyond models with elliptically contoured densities.
In this way, some of our results are new and others are more general than existing ones in the literature.

Existence of S-estimators and S-functionals is established under mild conditions.
Although existence of the estimators seems a basic requirement,
such results are missing for instance for multivariate regression in~\cite{vanaelst&willems2005}
and for linear mixed effects models in~\cite{copt2006high,chervoneva2014}.
We obtain robustness properties for S-estimators, such as breakdown point and influence function,
under mild conditions on collections of observations and under mild conditions on the distribution
of the observations.
High breakdown and a bounded influence function seem basic requirements for a robust method,
but both properties are not available for linear mixed effects models~\cite{copt2006high,chervoneva2014}.
For multivariate regression~\cite{vanaelst&willems2005}, the influence function is only determined at distributions with an
elliptical contoured density.
Finally, we establish consistency and asymptotic normality for S-estimators under mild conditions on the distribution of the observations.
A rigorous derivation is missing for multivariate regression~\cite{vanaelst&willems2005}, or is only available
for observations from a normal distribution~\cite{rousseeuw-yohai1984,chervoneva2014}.

We apply our asymptotic results, such as influence function and asymptotic normality, to the special case for which the
distribution of the observations corresponds to an elliptically contoured density.
In this way we retrieve earlier results found in~\cite{rousseeuw-yohai1984,lopuhaa1989,vanaelst&willems2005}.
Somewhat surprisingly, the asymptotic variances of our S-estimators for linear mixed effects models in which the response
has an elliptically contoured density, differ from the ones found in~\cite{copt2006high}.
We investigate this difference by means of a simulation study.

The paper is organized as follows.
In Section~\ref{sec:structured covariance model}, we explain the model in detail and
provide some examples of standard multivariate models that are included in our setup.
In Section~\ref{sec:definitions} we define the S-estimator and S-functional
and in Section~\ref{sec:existence} we give conditions under which they exist.
In Section~\ref{sec:continuity} we establish continuity of the S-functional, which is then used to
obtain consistency of the S-estimator.
Section~\ref{sec:bdp} deals with the breakdown point.
Section~\ref{sec:M-estimator equations} provides the preparation for Sections~\ref{sec:IF}
and~\ref{sec:asymptotic normality} in which we obtain the influence function and establish asymptotic normality.
Finally, in Section~\ref{sec:simulation}, we illustrate our results by means of a simulation
and investigate the performance of our estimators
by means of an application to data from a trial on the treatment of lead-exposed children.
All proofs and some technical lemmas are put in an Appendix at the end of the paper.
Other long and technical proofs are available as supplemental material~\cite{supplement}.

\section{Balanced models with structured covariances}
\label{sec:structured covariance model}
We consider independent observations $(\by_1,\bX_1),\ldots,(\by_n,\bX_n)$,
for which we assume the following model
\begin{equation}
\label{def:model}
\by_i
=
\bX_i\bbeta+\bu_i,
\quad
i=1,\ldots,n,
\end{equation}
where $\by_i\in\R^{k}$ contains repeated measurements for the $i$-th subject,
$\bbeta\in\R^q$ is an unknown parameter vector,
$\bX_i\in\R^{k\times q}$ is a known design matrix, and
$\mathbf{u}_i\in\R^{k}$ are unobservable independent mean zero random vectors with
covariance matrix $\bV\in\text{PDS}(k)$,
the class of positive definite symmetric $k\times k$ matrices.
The model is balanced in the sense that all~$\mathbf{y}_i$ have the same dimension.
Furthermore, we consider a structured covariance matrix, that is,
the matrix $\bV=\bV(\btheta)$ is a known function of unknown covariance parameters combined in a vector $\btheta\in\R^l$.
We first discuss some examples that are covered by this setup.
\begin{example}
\label{ex:linear mixed effects model}
An important case of interest is the (balanced) linear mixed effects model
\begin{equation}
\label{def:linear mixed effects model}
\by_i=\bX_i\bbeta+\bZ\gamma_i+\beps_i,
\quad
i=1,\ldots,n.
\end{equation}
This model arises from $\mathbf{u}_i=\mathbf{Z}\bgamma_i+\beps_i$,
for $i=1,\ldots,n$,
where $\mathbf{Z}\in\R^{k\times g}$ is known and $\bgamma_i\in\R^g$ and~$\beps_i\in\R^k$ are independent mean zero random variables, with
unknown covariance matrices~$\mathbf{G}$ and $\bR$, respectively.
In this case $\bV(\btheta)=\bZ\bG\bZ^T+\bR$ and
$\btheta=(\vch(\bG)^T,\vch(\bR)^T)^T$,
where
\begin{equation}
\label{def:vech}
\vch(\bA)=(a_{11},\ldots,a_{k1},a_{22},\ldots,a_{kk})
\end{equation}
is the unique $k(k+1)/2$-vector that stacks the columns of the lower triangle elements of a symmetric matrix $\bA$.
In full generality, the model is usually overparametrized and one may run into identifiability problems.
A more feasible example is obtained by taking $\bR=\sigma_0^2\bI_k$,
$\bZ=\left[\bZ_1 \,\cdots\,\bZ_r\right]$ and
$\bgamma_i=(\gamma_{i1}\ldots,\gamma_{ir})^T$,
where the $\bZ_j$'s are known $k\times g_j$ design matrices and
the $\gamma_{ij}\in\R^{g_j}$ are independent mean zero random variables with covariance matrix $\sigma_j^2\bI_{g_j}$,
for $j=1,\ldots,r$.
This leads to
\begin{equation}
\label{def:linear mixed effects model Copt}
\by_i=\bX_i\bbeta+\sum_{j=1}^r \bZ_j\gamma_{ij}+\beps_i,
\quad
i=1,\ldots,n,
\end{equation}
with $\bV(\btheta)=\sum_{j=1}^r\sigma_j^2\bZ_j\bZ_j^T+\sigma_0^2\bI_k$ and $\btheta=(\sigma_0^2,\sigma_1^2,\ldots,\sigma_r^2)$.
\end{example}
\begin{example}
\label{ex:multivariate regression}
An example with an unstructured covariance is
the multivariate linear regression model
\begin{equation}
\label{def:multivariate regression model}
\by_i=\bB^T\bx_i+\bu_i,
\qquad
i=1,\ldots,n,
\end{equation}
where $\bB\in\R^{q\times k}$ is a matrix of unknown parameters, $\bx_i\in\R^q$ is known,
and $\mathbf{u}_i$, for $i=1,\ldots,n$, are independent mean zero random variables with
covariance matrix~$\bV(\btheta)=\bC\in\text{PDS}(k)$.
In this case, the vector of unknown covariance parameters is given  by
\begin{equation}
\label{def:theta for unstructured}
\btheta=\vch(\bC)=(c_{11},\ldots,c_{1k},c_{22},\ldots,c_{kk})^T\in\R^{\frac12k(k+1)}.
\end{equation}
The model can be obtained as a special case of~\eqref{def:model}, by taking
$\bX_i=\bx_i^T\otimes \bI_k$ and $\bbeta=\vc(\bB^T)$, where $\otimes$ denotes the Kronecker product
and $\vc(\cdot)$ is the $k^2$-vector that stacks the columns of a matrix.
Clearly, the linear multiple regression model is a special case with $k=1$.
\end{example}
\begin{example}
\label{ex:autoregressive model}
Model~\eqref{def:model} also includes examples, for which
$\bu_1,\ldots,\bu_n$ are generated from a time series.
One example, is the case where $\bu_i$ has a covariance matrix with elements
\begin{equation}
\label{def:autogressive orde 1 covariance}
v_{st}=\sigma^2\rho^{|s-t|},
\quad
s,t=1,\ldots,k.
\end{equation}
This arises when the $\bu_i$'s are generated by an autoregressive process of order one.
The vector of unknown covariance parameters is $\btheta=(\sigma^2,\rho)\in(0,\infty)\times[-1,1]$.
A general stationary process leads to
\begin{equation}
\label{def:autogressive orde 1 covariance}
v_{st}=\theta_{|s-t|+1},
\quad
s,t=1,\ldots,k,
\end{equation}
in which case $\btheta=(\theta_1,\ldots,\theta_k)^T\in\R^k$,
where $\theta_{|s-t|+1}$ represents the autocovariance over lag~$|s-t|$.
\end{example}
\begin{example}
\label{ex:location-scale}
Also the multivariate location-scale model can be obtained as a special case of~\eqref{def:model},
by taking $\bX_i=\bI_k$, the $k\times k$ identity matrix.
In this case, $\bbeta\in\R^k$ is the unknown location parameter and $\bV(\btheta)$ is the unstructured covariance matrix
as in Example~\ref{ex:multivariate regression}, with $\btheta$ as in~\eqref{def:theta for unstructured}.
\end{example}

Throughout the manuscript we will assume that the parameter $\btheta$ is identifiable in the sense that,
\begin{equation}
\label{def:identifiable}
\bV(\btheta_1)=\bV(\btheta_2)
\quad\Rightarrow\quad
\btheta_1=\btheta_2.
\end{equation}
This is true for all models in
Examples~\ref{ex:multivariate regression},
\ref{ex:autoregressive model} and~\ref{ex:location-scale}.
This may not be true in general for the linear mixed effects model
in Example~\ref{ex:linear mixed effects model} with unknown $\vch(\bG)$ and $\vch(\bR)$.
For linear mixed effects models
in~\eqref{def:linear mixed effects model Copt},
identifiability of $\btheta=(\sigma_0^2,\sigma_1^2,\ldots,\sigma_r^2)$
holds for particular choices of the design matrices
$\bZ_1,\ldots,\bZ_r$.

\section{Definitions}
\label{sec:definitions}
We start by representing our observations as points in $\R^k\times\R^{kq}$  in the following way.
For $r=1,\ldots, k$, let $\bx_r^T$ denote the $r$-th row of the $k\times q$ matrix $\bX$,
so that $\bx_r\in\R^q$.
We represent the pair $\mathbf{s}=(\mathbf{y},\bX)$ as an element in $\R^k\times\R^{kq}$ defined by
$\bs^T=(\by^T,  \bx_{1}^T,\ldots,  \bx_{k}^T)$.
In this way our observations can be represented as $\bs_1,\ldots,\bs_n$, with $\bs_i=(\by_i,\bX_i)\in\R^k\times\R^{kq}$.
\subsection{S-estimator}
\label{sec:S-estimator}
S-estimators are defined by means of a function $\rho:\R\to[0,\infty)$ that satisfies the following properties
\begin{itemize}
\item[(R1)]
$\rho$ is symmetric around zero with $\rho(0)=0$ and $\rho$ is continuous at zero;
\item[(R2)]
There exists a finite constant $c_0>0$, such that $\rho$ is non-decreasing on $[0,c_0]$ and constant on $[c_0,\infty)$;
put $a_0=\sup\rho$.
\end{itemize}
The S-estimator $\bxi_n=(\bbeta_n,\btheta_n)$ is defined as the solution to the following minimization problem
\begin{equation}
\label{def:Smin estimator structured}
\begin{split}
&
\min_{\bbeta,\btheta}
\text{det}(\bV(\btheta))\\
&
\text{subject to}\\
&
\frac1n
\sum_{i=1}^n
\rho
\left(\sqrt{(\by_i-\bX_i\bbeta)^T
\bV(\btheta)^{-1}
(\by_i-\bX_i\bbeta)}\right)
\leq
b_0,
\end{split}
\end{equation}
where the minimum is taken over all $\bbeta\in\R^q$ and $\btheta\in\R^l$,
such that $\bV(\btheta)\in\text{PDS}(k)$,
with $\rho$ satisfying (R1)-(R2).

The S-estimator defined by~\eqref{def:Smin estimator structured} for the setup in~\eqref{def:model} includes
several specific cases that have been considered in the literature.
The original regression S-estimator introduced by Rousseeuw and Yohai~\cite{rousseeuw-yohai1984} is obtained
as a special case by
taking $\bX_i=\bx_i^T$ a $1\times q$ vector and $\bV(\btheta)=\sigma^2>0$.
S-estimators for multivariate location and scale, as considered in Davies~\cite{davies1987} and Lopuha\"a~\cite{lopuhaa1989}
can be obtained by taking $\bX_i$ and $\bV(\btheta)$ as in Example~\ref{ex:location-scale}.
For the multivariate regression model in Example~\ref{ex:multivariate regression},
S-estimators have been considered by Van Aelst and Willems~\cite{vanaelst&willems2005}.
Copt and Victoria-Feser~\cite{copt2006high} and Chervoneva and Vishnyakov~\cite{chervoneva2014} consider S-estimators for the parameters in the linear mixed effects model~\eqref{def:linear mixed effects model Copt}.

The constant $0<b_0<a_0$ in~\eqref{def:Smin estimator structured} can be chosen in agreement with an assumed underlying distribution.
For the multivariate regression model in~\cite{vanaelst&willems2005}, it is assumed that
$\by_i\mid\bX_i$ has an elliptically contoured density of the form
\begin{equation}
\label{eq:elliptical}
f_{\bmu,\bSigma}(\by)
=
\text{det}(\bSigma)^{-1/2}
h\left(
(\by-\bmu)^T
\bSigma^{-1}
(\by-\bmu)
\right),
\end{equation}
with $\bmu=\bX_i\bbeta$ and $\bSigma=\bV(\btheta)$ and $h:[0,\infty)\to[0,\infty)$.
For the linear mixed effects model in~\cite{copt2006high}, it is assumed that
$\by_i\mid\bX_i$ has a multivariate normal distribution, which is a special case of~\eqref{eq:elliptical}
with $h(t)=(2\pi)^{-k/2}\exp(-t/2)$.
When the underlying distribution corresponds to a density of the form~\eqref{eq:elliptical},
then a natural choice is $b_0=\E_{\mathbf{0},\bI}\rho(\|\bz\|)$,
where $\bz$ has density~\eqref{eq:elliptical} with~$(\bmu,\bSigma)=(\mathbf{0},\bI_k)$.
Finally, it should be emphasized that the ratio $b_0/a_0$ determines the breakdown point of the S-estimator
(see Theorem~\ref{th:BDP}),
as well as its limiting variance (see Corollary~\ref{cor:Asymp norm elliptical}).
By choosing the constant $c_0$ in~(R2) one then has to make a trade-off between robustness and efficiency.

Note that at this point we do not assume smoothness of $\rho$ or strict monotonicity on $[0,c_0]$.
This means that (R1)-(R2) allow the function $\rho(d)=1-\mathds{1}_{[-c_0,c_0]}(d)$, which corresponds to the minimum volume ellipsoid estimator in location-scale models (see~\cite{rousseeuw1985})
and to the least median of squares estimator in linear regression models (see~\cite{rousseeuw1984}).
Indeed, with $\rho(d)=1-\mathds{1}_{[-c_0,c_0]}(d)$, the S-estimator $(\bbeta_n,\btheta_n)$
corresponds to the smallest cylinder
\begin{equation}
\label{def:cylinder}
\mathcal{C}(\bbeta,\btheta,c_0)=
\left\{
(\by,\bX)\in\R^k\times\R^{kq}:
(\by-\bX\bbeta)^T
\bV(\btheta)^{-1}
(\by-\bX\bbeta)
\leq
c_0^2
\right\}
\end{equation}
that contains at least $n-nb_0$ points.

\begin{remark}
\label{rem:location-scale}
Clearly, the definition of the S-estimator in~\eqref{def:Smin estimator structured} has great similarities with
the S-estimator for multivariate location and covariance (see~\cite{davies1987} and~\cite{lopuhaa1989}),
defined as the solution~$(\mathbf{t}_n,\mathbf{C}_n)$ to the minimization problem
\begin{equation}
\label{def:Smin location covariance}
\begin{split}
&
\min_{\mathbf{t},\mathbf{C}}
\text{\rm det}(\mathbf{C})\\
&
\text{subject to}\\
&
\frac1n
\sum_{i=1}^n
\rho
\left(\sqrt{(\by_i-\mathbf{t})^T\mathbf{C}^{-1}(\by_i-\mathbf{t})}\right)
\leq
b_0,
\end{split}
\end{equation}
where the minimum is taken over all $\mathbf{t}\in\R^k$ and $\mathbf{C}\in\text{PDS}(k)$.
Even more so, if all $\bX_i$ are assumed to be equal to the same design matrix $\bX$ of full rank, as was done in~\cite{copt2006high,copt&heritier2007}.
However, there is a subtle, but important difference between minimization problems~\eqref{def:Smin location covariance} and~\eqref{def:Smin estimator structured}.
The important difference is that in~\eqref{def:Smin location covariance} we minimize over \emph{all} positive definite symmetric $k\times k$ matrices
$\mathbf{C}$, whereas in~\eqref{def:Smin estimator structured}, we \emph{only} minimize over positive definite symmetric $k\times k$ matrices~$\bV(\btheta)$, which can arise as the image of the mapping $\btheta\mapsto \bV(\btheta)$.
The latter collection is a subset of the other:
\[
\left\{
\bV(\btheta)\in\mathrm{PDS}(k):\btheta\in\R^l
\right\}
\subset
\text{PDS}(k),
\]
and will typically be a strictly smaller subset.
This means that the properties of $\bV(\btheta_n)$
and~$\mathbf{C}_n$ are related, but the properties of $\bV(\btheta_n)$ \emph{cannot} simply be derived from
properties of $\mathbf{C}_n$,
not even in the case where all $\bX_i$ are equal to the same $\bX$.
In fact, this will lead to limiting covariances that differ from the ones found
in~\cite{copt2006high}, see Corollary~\ref{cor:Asymp norm elliptical}.
\end{remark}
\subsection{S-functional}
\label{sec:S-functional}
The concept of S-functional is needed to investigate local robustness properties of the corresponding S-estimator,
such as the influence function (see Section~\ref{sec:IF}).
Let $\bs=(\by,\bX)$ have a probability distribution $P$ on $\R^k\times\R^{kq}$.
The S-functional $\bxi(P)=(\bbeta(P),\btheta(P))$
is defined as the solution to the following
minimization problem:
\begin{equation}
\label{def:Smin structured}
\begin{split}
&
\min_{\bbeta,\btheta}
\text{det}(\bV(\btheta))\\
&
\text{subject to}\\
&
\int \rho
\left(\sqrt{(\by-\bX\bbeta)^T\bV(\btheta)^{-1}(\by-\bX\bbeta)}\right)\,\text{d}P(\by,\bX)
\leq
b_0,
\end{split}
\end{equation}
where the minimum is taken over all $\bbeta\in\R^q$ and $\btheta\in\R^l$,
such that $\bV(\btheta)\in\text{PDS}(k)$,
with $\rho$ satisfying (R1)-(R2).

As a special case, we obtain the S-estimator $\bxi_n=(\bbeta_n,\btheta_n)$ by taking~$P=\mathbb{P}_n$,
the empirical measure corresponding to the observations
$(\by_1,\bX_1),\ldots,(\by_n,\bX_n)$.
In view of this connection, existence and consistency of solutions to~\eqref{def:Smin estimator structured}
will follow from general results on the existence and the continuity of
solutions to~\eqref{def:Smin structured}.

The definition of the S-functionals for the multivariate location-scale model given in Lopuha\"a~\cite{lopuhaa1989}
and for the multivariate regression model given by Van Aelst and Willems~\cite{vanaelst&willems2005}
can be obtained as special cases of~\eqref{def:Smin structured},
by choosing $\bX$, $\bbeta$ and $\bV(\btheta)$ as in Examples~\ref{ex:location-scale} and~\ref{ex:multivariate regression},
respectively.
Copt and Victoria-Feser~\cite{copt2006high} do not pay attention to S-functionals
or the influence function in the linear mixed effects model~\eqref{def:linear mixed effects model Copt}.
However, S-functionals for linear mixed effects models can be also be obtained
as a special case of~\eqref{def:Smin structured},
by choosing $\bX$, $\bbeta$ and $\bV(\btheta)$ as in Example~\ref{ex:linear mixed effects model}.

\section{Existence}
\label{sec:existence}
We will first establish existence of the S-functional $\bxi(P)$ defined by~\eqref{def:Smin structured},
under particular conditions on the probability measure $P$.
As a consequence, this will also yield the existence of the S-estimator, defined by~\eqref{def:Smin estimator structured}.
Recall that $(\by_1,\bX_1),\ldots,(\by_n,\bX_n)$ are represented as points in $\R^k\times\R^{kq}$.
Note however, that for linear models with intercept the first column of each $\bX_i$ consists of 1's.
This means that the points $(\by_i,\bX_i)$ are concentrated in a lower dimensional subset of $\R^k\times\R^{kq}$.
A similar situation occurs when all $\bX_i$ are equal to the same design matrix, such as in~\cite{copt2006high}.
In view of this, define $\mathcal{X}\subset\R^{kq}$ as the subset with the lowest dimension
$p=\text{dim}(\mathcal{X})\leq kq$ satisfying
\begin{equation}
\label{def:Xspace}
P(\bX\in \mathcal{X})=1.
\end{equation}
Hence, $P$ is then concentrated on the subset $\R^k\times \mathcal{X}$ of $\R^k\times\R^{kq}$, which
is of dimension~$k+p$, which may be of smaller than $k+kq$.

The first condition we require, expresses the fact that $P$
cannot have too much mass at infinity, in relation to the ratio $r=b_0/a_0$.
\begin{itemize}
\item[$(\mathrm{C1}_\epsilon)$]
There exists a compact set $K_\epsilon\subset\R^k\times \mathcal{X}$,
such that $P(K_\epsilon)\geq r+\epsilon$.
\end{itemize}
The second condition requires that $P$ cannot have too much mass at arbitrarily thin strips in $\R^k\times \mathcal{X}$.
For $\balpha\in\R^{k+kq}$, such that $\|\balpha\|=1$, $\ell\in\R$ and $\delta\geq0$, we define a strip
$H(\balpha,\ell,\delta)$ as follows:
\begin{equation}
\label{def:strip}
H(\balpha,\ell,\delta)
=
\left\{
\mathbf{s}\in\R^k\times\R^{kq}: \ell-\delta/2\leq \balpha^T\mathbf{s}\leq \ell+\delta/2
\right\}.
\end{equation}
Defined in this way, a strip is the area between two parallel hyperplanes
which are symmetric around the hyperplane $H(\balpha,\ell,0)
=
\left\{
\mathbf{s}\in\R^k\times\R^{kq}: \balpha^T\mathbf{s}=\ell
\right\}$.
Since the distance between two parallel hyperplanes
$\balpha^T\mathbf{s}=\ell_1$ and $\balpha^T\mathbf{s}=\ell_2$ is $|\ell_1-\ell_2|$,
the strip $H(\balpha,\ell,\delta)$ defined as in~\eqref{def:strip} has width~$\delta$.
We require the following condition
\begin{itemize}
\item[$(\mathrm{C2}_\epsilon)$]
The value
\[
\delta_\epsilon
=
\inf
\left\{
\delta:
P\left(H(\balpha,\ell,\delta)
\right)\geq \epsilon,
\balpha\in\R^{k+kq},\|\balpha\|=1,\ell\in\R,\delta\geq 0
\right\}
\]
is strictly positive.
\end{itemize}
According to~\eqref{def:Xspace},  in $(\mathrm{C2}_\epsilon)$ one only needs to consider strips in $\R^k\times \mathcal{X}$.

Both conditions are satisfied for any $0<\epsilon\leq 1-r$ by any probability measure $P$
that is absolutely continuous.
Clearly, condition~$(\mathrm{C1}_\epsilon)$ holds for any $0\leq\epsilon\leq 1-r$ for
the empirical measure~$\mathbb{P}_n$ corresponding to a collection of $n$ points
$\mathcal{S}_n=\{\mathbf{s}_1,\ldots,\mathbf{s}_n\}\subset\R^k\times \mathcal{X}$.
Condition~$(\mathrm{C2}_\epsilon)$ for $\epsilon=(k+p+1)/n$
is also satisfied by the empirical measure $\mathbb{P}_n$, when the collection~$\mathcal{S}_n$ is in \emph{general position}, i.e.,
no subset $J\subset \mathcal{S}_n$ of $k+p+1$ points is contained in the same hyperplane
in~$\R^k\times \mathcal{X}$.
Conditions $(\mathrm{C1}_\epsilon)$ and~$(\mathrm{C2}_\epsilon)$ together,
are similar to condition~$(\mathrm{C}_\epsilon)$
in~\cite{lopuhaa1989}.
The reason that $(\mathrm{C1}_\epsilon)$ slightly deviates from~\cite{lopuhaa1989}, is to handle
the presence of $\bX$ in minimization problem~\eqref{def:Smin structured}.
\begin{remark}
\label{rem:condition C2}
Note that condition $(\mathrm{C2}_\epsilon)$ is equivalent with
\begin{equation}
\label{eq:general position}
\omega_\epsilon
=
\inf_{P(J)\geq \epsilon}\,
\inf_{\|\balpha\|=1}\,
\inf_{\ell\in\R}\,
\sup_{\bs\in J}\,
|\balpha^T\bs-\ell|>0,
\end{equation}
where the infima are taken over all subsets $J\subset\R^k\times \mathcal{X}$ with $P(J)\geq \epsilon$,
all vectors $\balpha\in\R^{k+kq}$, with $\|\balpha\|=1$,
and levels $\ell\in\R$.
Details can be found in~\cite{supplement}.
\end{remark}
To establish existence of the S-functional, we follow the reasoning in~\cite{lopuhaa1989}.
The idea is to argue that one can restrict oneself to a compact set for finding solutions to~\eqref{def:Smin structured}.
When the object function in~\eqref{def:Smin structured} is continuous, this immediately yields existence of a solution
of~\eqref{def:Smin structured}.
To this end, we assume the following condition.
\begin{itemize}
\item[(V1)]
The mapping $\btheta\mapsto\bV(\btheta)$ is continuous.
\end{itemize}
The lemma below is fundamental for the existence of the S-functional.
It requires that the identity is in
$\mathcal{V}=\{\bV(\btheta)\in\text{PDS}(k):\btheta\in\R^l\}$
and that $\mathcal{V}$ is closed under multiplication with a positive scalar.
\begin{itemize}
\item[(V2)]
There exists a $\btheta\in\R^l$, such that $\bV(\btheta)=\bI_k$.
For any $\bV(\btheta)\in \mathcal{V}$ and any $\alpha>0$, it holds that
$\alpha\bV(\btheta)=\bV(\btheta')$, for some $\btheta'\in\R^l$.
\end{itemize}
Conditions (V1)-(V2) are not very restrictive.
For example, all models in Examples~\ref{ex:linear mixed effects model} to~\ref{ex:location-scale} satisfy these conditions.

For any $k\times k$ matrix $\bA$, let $\lambda_k(\bA)\leq\cdots\leq\lambda_1(\bA)$ denote the eigenvalues of $\bA$.
We then have the following key lemma for the existence of S-functionals.
The lemma is similar to Lemma 1 in~\cite{lopuhaa1989} and its proof can be found in~\cite{supplement}.

\begin{lemma}
\label{lem:compact structured}
Let $(\bbeta,\btheta)\in\R^q\times\R^l$, $0<m_0<\infty$, $0<c<\infty$, and $0<\epsilon<1$,
and suppose that the mapping $\btheta\mapsto \bV(\btheta)$ satisfies (V2).
Then the following properties hold.
\begin{enumerate}
\item[(i)]
If $P$ satisfies $(\text{C2}_\epsilon)$ and
$P(\mathcal{C}(\bbeta,\btheta,c)))\geq\epsilon$,
then $\lambda_k(\bV(\btheta))\geq a_1>0$,
where $a_1$ only depends on~$c$ and the width $\delta_\epsilon$ from condition $(\text{C2}_\epsilon)$.
\item[(ii)]
Suppose $\int \rho(\|\mathbf{y}\|/m_0)\,\dd P(\mathbf{s})\leq b_0$.
Then for any solution $(\bbeta,\btheta)$ of~\eqref{def:Smin structured},
which is such that $\lambda_k(\bV(\btheta))\geq a_1>0$,
it holds that $\lambda_1(\bV(\btheta))\leq a_2<\infty$,
where $a_2$ only depends on $a_1$ and~$m_0$.
\item[(iii)]
Let $P$ satisfy $(\text{C2}_\epsilon)$
and suppose that $P(\mathcal{C}(\bbeta,\btheta,c))\geq a>0$.
Suppose there exists a compact set $K\in\R^k\times \mathcal{X}$,
such that $P(K)\geq 1-a+\epsilon$.
If $\lambda_1(\bV(\btheta))\leq a_2<\infty$, then
$\|\bbeta\|\leq M<\infty$, where $M$ only depends on~$c$, $a_2$,
the set $K$,
and a constant $\gamma_\epsilon>0$ that can be deduced from condition $(\text{C2}_\epsilon)$.
\end{enumerate}
\end{lemma}
Lemma~\ref{lem:compact structured} will ensure that there exists a compact set
that contains all pairs $(\bbeta,\bV(\btheta))$ that correspond to possible solutions
$(\bbeta,\btheta)$ of~\eqref{def:Smin structured}.
To establish that possible solutions $(\bbeta,\btheta)$ of~\eqref{def:Smin structured} are in a compact set,
we need that the pre-image
$\{\btheta\in\R^l: \bV(\btheta)\in K\}$ of a compact set $K\subset\R^{k\times k}$ is again compact.
Recall that subsets of $\R^l$ are compact if and only if they are closed and bounded,
and note that the pre-image of a continuous mapping of a closed set is closed.
Hence, in view of condition (V1), it suffices to require the following condition.
\begin{itemize}
\item[(V3)]
The mapping $\btheta\mapsto \bV(\btheta)$ is such that the pre-image of a bounded set is bounded.
\end{itemize}
Condition (V3) is satisfied by all models in Examples~\ref{ex:linear mixed effects model} to~\ref{ex:location-scale},
including the linear mixed effects model of Example~\ref{ex:linear mixed effects model}, as long as the matrix $\bZ$ is of full rank.
We then have the following theorem.
\begin{theorem}
\label{th:existence structured}
Consider minimization problem~\eqref{def:Smin structured} with $\rho$ satisfying (R1)-(R2).
Suppose that $P$ satisfies $(\text{C1}_\epsilon)$ and $(\text{C2}_\epsilon)$, for some $0<\epsilon\leq 1-r$,
where $r=b_0/a_0$,
and suppose that $\bV$ satisfies (V1)-(V3).
Then there exists at least one solution to~\eqref{def:Smin structured}.
\end{theorem}

The theorem has a direct corollary for the existence of the S-estimator, when dealing with a collections of points.
Let $\mathcal{S}_n=\{\bs_1,\ldots,\bs_n\}$, with $\bs_i=(\by_i,\bX_i)$ be a collection of $n$ points
in~$\R^k\times \mathcal{X}$.
Define
\begin{equation}
\label{def:k(S)}
\kappa(\mathcal{S}_n)
=
\text{maximal number of points of $\mathcal{S}_n$ lying on the same hyperplane in~$\R^k\times \mathcal{X}$.}
\end{equation}
For example, if the distribution $P$ is absolutely continuous, then
$\kappa(\mathcal{S}_n)\leq k+p$ with probability one.
We then have the following corollary.
\begin{corollary}
\label{cor:existence estimator structured}
Consider minimization problem~\eqref{def:Smin estimator structured} with
$\rho$ satisfying (R1)-(R2),
for a collection $\mathcal{S}_n=\{\bs_1,\ldots,\bs_n\}\subset\R^k\times \mathcal{X}$,
with $\bs_i=(\by_i,\bX_i)$, for $i=1,\ldots,n$.
Suppose that~$\bV$ satisfies (V1)-(V3).
If $\kappa(\mathcal{S}_n)+1\leq n(1-r)$, where $r=b_0/a_0$,
then there exists at least one solution $\bxi_n=(\bbeta_n,\btheta_n)$
to the minimization problem~\eqref{def:Smin estimator structured}.
\end{corollary}
Copt and Victoria-Feser~\cite{copt2006high} consider S-estimators
for the linear mixed effects model~\eqref{def:linear mixed effects model Copt}.
Despite their Proposition~1 about the asymptotic behavior of solutions to their S-minimization problem~\cite[equation (7)]{copt2006high},
the actual existence of such a solution is not established.
However, this now follows from our Corollary~\ref{cor:existence estimator structured}.
In their case, $\bV(\btheta)$ satisfies conditions~(V1) and~(V2).
It can be seen, that if all matrices $\bZ_j$, for $j=1,\ldots,r$, are of full rank, then~$\bV(\btheta)$ also satisfies (V3).
The translated bi-weight $\rho$-function proposed in~\cite{copt2006high} satisfies~(R1)-(R2).
Finally, under their assumption that $\bX_i=\bX$ is the same and $\by_i\mid\bX\sim N_k(\bX\bbeta,\bV(\btheta))$,
it follows that $\kappa(\mathcal{S}_n)\leq k$.
It then follows from Corollary~\ref{cor:existence estimator structured} that with $b_0\leq a_0(n-k-1)/n$,
at least one solution
to their S-minimization problem exists.

For the multivariate regression model from Example~\ref{ex:multivariate regression}, Van Aelst \& Willems \cite{vanaelst&willems2005} do not explicitly
prove existence of the S-estimator.
Since in their case, $\bV(\btheta)=\bC\in\text{PDS}(k)$ satisfies (V1)-(V3) and the conditions imposed in~\cite{vanaelst&willems2005}
on the $\rho$-function
satisfy (R1)-(R2),
the existence of their S-estimator now also follows from Corollary~\ref{cor:existence estimator structured},
when $b_0$ is chosen suitably.

Existence of S-estimators is obtained from existence of S-functionals at the empirical measure $\mathbb{P}_n$.
The following corollary shows that existence can be established in general, for probability measures that are close to $P$.
It requires the following condition on $P$.
\begin{itemize}
\item[(C3)]
Let $\mathfrak{C}$ be the class of all measurable convex subsets of $\R^k\times \R^{kq}$.
Every $C\in \mathfrak{C}$ is a $P$-continuity set, i.e., $P(\partial C)=0$,
where~$\partial C$ denotes the boundary of $C$.
\end{itemize}
\begin{corollary}
\label{cor:existence weak convergence}
Suppose that $\rho$ satisfies (R1)-(R2) and $\bV$ satisfies (V1)-(V3).
Let $P_t$, $t\geq0$ be a sequence of probability measures on $\R^k\times \R^{kq}$ that converges weakly to~$P$, as $t\to\infty$.
Suppose that~$P$ satisfies (C3), as well as $(\text{C1}_{\epsilon'})$ and $(\text{C2}_\epsilon)$, for some $0<\epsilon<\epsilon'\leq 1-r=b_0/a_0$.
Then, for $t$ sufficiently large, the minimization problem~\eqref{def:Smin structured} with
probability measure $P_t$
has at least one solution~$\bxi(P_t)$.
\end{corollary}
Condition (C3) is needed to apply~\eqref{eq:ranga rao structured}.
Clearly, this condition is satisfied if $P$ is absolutely continuous.

\section{Continuity and consistency}
\label{sec:continuity}
Consider a sequence $P_t$, $t\geq0$, of probability measures on $\R^k\times\R^{kq}$ that converges weakly to $P$, as $t\to\infty$.
By continuity of the S-functional $\bxi(P)$ we mean that $\bm{\xi}(P_t)\to\bm{\xi}(P)$, as $t\to\infty$.
An example of such a sequence is the sequence of empirical measures $\mathbb{P}_n$, $n=1,2,\ldots$, that converges weakly to $P$, almost surely.
Continuity of the S-functional for this sequence would then mean that the S-estimator $\bxi_n$ is consistent,
i.e., $\bxi_n=\bxi(\mathbb{P}_n)\to\bxi(P)$,
almost surely.

We require an additional condition for the function~$\rho$.
\begin{itemize}
\item[(R3)]
$\rho$ is continuous and strictly increasing on $[0,c_0]$.
\end{itemize}
For $\bs=(\by,\bX)$ and $\bxi=(\bbeta,\btheta)$, define the Mahalanobis distances by
\begin{equation}
\label{def:Mahalanobis distance structured}
d^2(\bs,\bxi)
=
d^2(\bs,\bbeta,\btheta)
=
(\by-\bX\bbeta)^T\mathbf{V}(\btheta)^{-1}(\by-\bX\bbeta).
\end{equation}
We then have the following theorem for the S-functional $\bxi(P)=(\bbeta(P),\btheta(P))$.
\begin{theorem}
\label{th:continuity structure}
Let $P_t$, $t\geq0$ be a sequence of probability measures on $\R^k\times\R^{kq}$ that converges weakly to~$P$, as $t\to\infty$,
and let $\bxi(P_t)$ be a solution to minimization problem~\eqref{def:Smin structured} with probability measure $P_t$.
Suppose that $\rho$ satisfies (R1)-(R3) and $\bV$ satisfies  (V1)-(V3).
Suppose that~$P$ satisfies (C3), as well as $(\text{C1}_{\epsilon'})$ and $(\text{C2}_\epsilon)$,
for some $0<\epsilon<\epsilon'\leq 1-r=b_0/a_0$.
If the solution~$\bxi(P)$ of~\eqref{def:Smin structured} is unique,
then for any sequence of solutions $\bxi(P_t)$, $t\geq 0$, it holds
\[
\lim_{t\to\infty}
\bxi(P_t)
=
\bxi(P).
\]
\end{theorem}
Theorem~\ref{th:continuity structure} is an extension of Theorem 3.1 in~\cite{lopuhaa1989}
on the continuity of S-functionals for multivariate location and scale.
Continuity of S-functionals for multiple regression has been investigated in~\cite{fasano-maronna-sued-yohai2012}.

Continuity of the S-functional will be used to derive the influence function of the S-estimator in Section~\ref{sec:IF}.
Another nice consequence of the continuity of the S-functional is, that one can directly obtain consistency of the S-estimator.
Consider the S-estimator $\bxi_n$ defined by minimization problem~\eqref{def:Smin estimator structured}.
Recall that $\bxi_n=\bxi(\mathbb{P}_n)$, so that we can use Theorem~\ref{th:continuity structure} to
establish consistency of the S-estimator.
\begin{corollary}
\label{cor:consistency S-estimator}
Let $\bxi_n$ be a solution to minimization problem~\eqref{def:Smin estimator structured}.
Suppose $\rho$ satisfies (R1)-(R3)
and $\bV$ satisfies (V1)-(V3).
Suppose that~$P$ satisfies (C3) as well as $(\text{C1}_{\epsilon'})$ and~$(\text{C2}_\epsilon)$,
for some $0<\epsilon<\epsilon'\leq 1-r=b_0/a_0$.
If the solution $\bxi(P)$ of~\eqref{def:Smin structured} is unique,
then
\[
\lim_{n\to\infty}
\bxi_n
=
\bxi(P),
\]
with probability one.
\end{corollary}
Theorem~\ref{th:continuity structure} and Corollary~\ref{cor:consistency S-estimator} require that $\bxi(P)$ is
the unique solution to minimization problem~\eqref{def:Smin structured}.
An example of a distribution $P$ for which $\bxi(P)$ is unique,
is when $P$ is such that~$\by\mid\bX$ has an elliptically contoured density~\eqref{eq:elliptical}.
This situation is very similar to that of multivariate location-scale S-estimators,
for which Davies~\cite[Theorem 1]{davies1987} shows that the corresponding S-minimization problem~\eqref{def:Smin structured}
has a unique solution.
The next theorem is a direct consequence of that result.
Its proof can be found in~\cite{supplement}.
\begin{theorem}
\label{th:davies}
Suppose that $\rho:\R\to[0,\infty)$ satisfies (R1)-(R2)
and suppose that the probability distribution $P$ of $(\by,\bX)$ is such that $\by\mid\bX$ has an
elliptically contoured density~$f_{\bmu,\bSigma}$ from~\eqref{eq:elliptical}, with $\bmu=\bX\bbeta_0$ and $\bSigma=\bV(\btheta_0)$.
Suppose that $h$ in~\eqref{eq:elliptical} is non-increasing and such that the functions $-\rho$ and $h$ have at least one common
point of decrease $d_0>0$, i.e.,
\[
\rho(s)<\rho(d_0)<\rho(t)
\quad\text{and}\quad
h(s)>h(d_0)>h(t)
\]
for all $s,t\geq 0$, such that $s<d_0<t$.
If $\bX^T\bX$ is non-singular with probability one, then the minimization problem
\begin{equation}
\label{def:Smin regression}
\begin{split}
&
\min_{\bbeta,\btheta}
\text{det}(\bV(\btheta))\\
&
\text{subject to}\\
&
\int
\rho
\left(\sqrt{(\by-\bX\bbeta)^T\mathbf{V}(\btheta)^{-1}(\by-\bX\bbeta)}\right)
f_{\bmu,\bSigma}(\by)
\dd\by
\leq
b_0,
\end{split}
\end{equation}
where the minimum is taken over all $\bbeta\in\R^q$ and $\btheta\in\R^l$,
such that $\bV(\btheta)\in\text{PDS}(k)$,
has the unique solution $(\bbeta,\btheta)=(\bbeta_0,\btheta_0)$ with probability one.
\end{theorem}
Minimization problem~\eqref{def:Smin regression} seems to be slightly different from the one in~\eqref{def:Smin structured}.
However, note that when~$P$ is such that $\bX$ is equal to a single value with probability one,
both minimization problems are identical.
This situation was considered, e.g., in~\cite{copt2006high}.

An elliptically contoured density for $\by_i\mid\bX_i$ in the context of S-estimators for specific cases of the model~\eqref{def:model}
has been assumed
in~\cite{davies1987} for the multivariate location-scale model
of Example~\ref{ex:location-scale},
in~\cite{vanaelst&willems2005} for the multivariate regression model
of Example~\ref{ex:multivariate regression},
and in~\cite{copt2006high} for
the linear mixed effects model~\eqref{def:linear mixed effects model Copt}.
More precisely, in~\cite{copt2006high} it is assumed that $\bX_i=\bX$ and that $\by_i\mid\bX$ has a multivariate normal distribution.
In that case, the function~$h$ in~\eqref{eq:elliptical} satisfies all the conditions of Theorem~\ref{th:davies}.

\section{Global robustness: the breakdown point}
\label{sec:bdp}
Consider a collection of points $\mathcal{S}_n=\{\bs_i=(\by_i,\bX_i),i=1,\ldots,n\}\subset \R^k\times \mathcal{X}$.
To emphasize the dependence on the collection $\mathcal{S}_n$,
denote by $\bxi_n(\mathcal{S}_n)=(\bbeta_n(\mathcal{S}_n),\btheta_n(\mathcal{S}_n))$,
the S-estimator,
as defined in~\eqref{def:Smin estimator structured}.
To investigate the global robustness of S-estimators,
we compute that finite-sample (replacement) breakdown point.
For a given collection $\mathcal{S}_n$ the finite-sample breakdown point
(see Donoho and Huber~\cite{donoho&huber1983})
of a regression S-estimator $\bbeta_n$ is defined as the smallest proportion of points
from~$\mathcal{S}_n$ that one needs to replace in order to
carry the estimator over all bounds.
More precisely,
\begin{equation}
\label{def:BDP beta}
\epsilon_n^*(\bbeta_n,\mathcal{S}_n)
=
\min_{1\leq m\leq n}
\left\{
\frac{m}{n}:
\sup_{\mathcal{S}_m'}
\left\|
\bbeta_n(\mathcal{S}_n)-\bbeta_n(\mathcal{S}_m')
\right\|
=\infty
\right\},
\end{equation}
where the minimum runs over all possible collections $\mathcal{S}_m'$ that can be obtained from $\mathcal{S}_n$
by replacing~$m$ points of $\mathcal{S}_n$ by arbitrary points in $\R^k\times \mathcal{X}$.

The estimator $\btheta_n$ determines the covariance estimator $\bV_n=\bV(\btheta_n)$.
For this reason it seems natural to let the breakdown point of $\btheta_n$ correspond to the breakdown of a covariance estimator.
We define the finite sample (replacement) breakdown point of the S-estimator $\btheta_n$ at a collection~$\mathcal{S}_n$, as
\begin{equation}
\label{def:BDP theta}
\epsilon_n^*(\btheta_n,\mathcal{S}_n)
=
\min_{1\leq m\leq n}
\left\{
\frac{m}{n}:
\sup_{\mathcal{S}_m'}
\text{dist}(\bV(\btheta_n(\mathcal{S}_n))),\bV(\btheta_n(\mathcal{S}_m'))
=\infty
\right\},
\end{equation}
with $\text{dist}(\cdot,\cdot)$ defined as
$\text{dist}(\bA,\mathbf{B})
=
\max\left\{
\left|\lambda_1(\bA)-\lambda_1(\mathbf{B})\right|,
\left|\lambda_k(\bA)^{-1}-\lambda_k(\mathbf{B})^{-1}\right|
\right\}$,
where the minimum runs over all possible collections $\mathcal{S}_m'$ that can be obtained from $\mathcal{S}_n$
by replacing~$m$ points of $\mathcal{S}_n$ by arbitrary points in $\R^k\times \mathcal{X}$.
So the breakdown point of $\btheta_n$ is the smallest proportion of points from~$\mathcal{S}_n$ that one needs to replace in order to
make the largest eigenvalue of $\bV(\btheta(\mathcal{S}_m'))$ arbitrarily large (explosion), or
to make the smallest eigenvalue of~$\bV(\btheta(\mathcal{S}_m'))$ arbitrarily small (implosion).

Good global robustness is illustrated by a high breakdown point.
The breakdown point of the S-estimators is given the theorem below.
It extends the results for S-estimators of multivariate location and scale,
see~\cite{davies1987} and~\cite{lopuhaa&rousseeuw1991}, and S-estimators for multivariate regression, see~\cite{vanaelst&willems2005}.
For S-estimators in the linear mixed effects model considered in~\cite{copt2006high}, the breakdown point
has not been established.
This will now follow as a special case from the next theorem.
Its proof can be found in~\cite{supplement}.
\begin{theorem}
\label{th:BDP}
Consider the minimization problem~\eqref{def:Smin estimator structured} with $\rho$ satisfying (R1)-(R2).
Suppose that~$\bV$ satisfies (V1)-(V3).
Let $\mathcal{S}_n\subset \R^k\times \mathcal{X}$ be a collection of~$n$ points $\bs_i=(\by_i,\bX_i)$,  $i=1,\ldots,n$.
Let $r=b_0/a_0$ and suppose that $0<r\leq (n-\kappa(\mathcal{S}_n))/(2n)$,
where $\kappa(\mathcal{S}_n)$ is defined by~\eqref{def:k(S)}.
Then for any solution $(\bbeta_n,\btheta_n)$ of minimization problem~\eqref{def:Smin structured},
\[
\begin{split}
\frac{\lfloor{(n+1)/2}\rfloor}n
\geq
\epsilon_n^*(\bbeta_n,\mathcal{S}_n)
&\geq
\frac{\lceil{nr}\rceil}n,\\
\epsilon_n^*(\btheta_n,\mathcal{S}_n)
&=
\frac{\lceil{nr}\rceil}n.
\end{split}
\]
\end{theorem}
The largest possible value of the breakdown point occurs when $r=(n-\kappa(\mathcal{S}_n))/(2n)$, in which case
$\lceil nr\rceil/n=\lceil (n-\kappa(\mathcal{S}_n))/2\rceil/n=\lfloor (n-\kappa(\mathcal{S}_n)+1)/2\rfloor/n$.
When the collection $\mathcal{S}_n$ is in general position, then $\kappa(\mathcal{S}_n)= k+p$.
In that case the breakdown point of both estimators is at least equal
to~$\lfloor (n-k-p+1)/2\rfloor/n$.
When all $\bX_i$ are equal to the same $\bX$,  in~\cite{copt2006high,copt&heritier2007},
one has $p=0$ and $\kappa(\mathcal{S}_n)=k$.
In that case, the breakdown point of $\btheta_n$ is equal to $\lfloor (n-k+1)/2\rfloor/n$.
This coincides with the maximal breakdown point for affine equivariant estimators
for $k\times k$ covariance matrices (see~\cite[Theorem 6]{davies1987}).

\begin{remark}
\label{rem:BDP}
Van Aelst \& Willems~\cite{vanaelst&willems2005} also take into account the case $r>(n-\kappa(\mathcal{S}_n))/(2n)$.
For this case, by replacing $\lceil{n-nr}\rceil-\kappa(\mathcal{S}_n)$ points, a specific solution to the S-minimization problem
is constructed that breaks down.
However, since there may be multiple solutions to the S-minimization problem, this does not necessarily mean that all solutions break down.
In the proof of our Theorem~\ref{th:BDP}, for the case $r\leq (n-\kappa(\mathcal{S}_n))/(2n)$, we show that
all solutions to~\eqref{def:Smin estimator structured} do not break down, when replacing at most $\lceil nr\rceil-1$ points,
and that the covariance part of all solutions do break down, when replacing  $\lceil nr\rceil$ points.
For the case $r>(n-\kappa(\mathcal{S}_n))/(2n)$, we can show that all solutions
to~\eqref{def:Smin estimator structured} do not break down,
when replacing at most $\lceil{n-nr}\rceil-\kappa(\mathcal{S}_n)-1$ points.
\end{remark}

\section{Score equations}
\label{sec:M-estimator equations}
Up to this point, properties of S-functionals and S-estimators have been derived from the minimization
problems~\eqref{def:Smin estimator structured} and~\eqref{def:Smin structured}.
To obtain the influence function and to establish the limiting distribution of S-estimators,
we use the score equations that can be found by differentiation of the Lagrangian corresponding to the
constrained minimization problems.
To this end, we require the following additional condition on the function $\rho$,
\begin{itemize}
\item[(R4)]
$\rho$ is continuously differentiable and $u(s)=\rho'(s)/s$ is continuous,
\end{itemize}
and the following condition on the mapping $\btheta\mapsto\bV(\btheta)$,
\begin{itemize}
\item[(V4)]
$\bV(\btheta)$ is continuously differentiable.
\end{itemize}
Obviously, condition (V4) implies the former condition~(V1).
\subsection{General covariance structures}
\label{subsec:general V}
Let $\bxi_P=(\bbeta_P,\btheta_P)$ be a solution to minimization problem~\eqref{def:Smin structured}.
If we denote the corresponding Lagrange multiplier by $\lambda_P$,
then the pair $(\bxi_P,\lambda_P)$ is a zero of all partial derivatives
$\partial L_P/\partial \bbeta$,
$\partial L_P/\partial \btheta$,
and $\partial L_P/\partial \lambda$,
where $L_P$ is the Lagrangian given by
\[
L_P(\bxi,\lambda)
=
\log\text{det}(\bV(\btheta))
-
\lambda
\left\{
\int \rho
\left(\sqrt{(\by-\bX\bbeta)^T\bV(\btheta)^{-1}(\by-\bX\bbeta)}\right)\,\text{d}P(\by,\bX)
-
b_0
\right\}.
\]
If $\E_P\|\bX\|<\infty$, then under conditions~(R4) and~(V4), one may interchange the order of integration and differentiation
in $\partial L_P/\partial \bbeta$ and $\partial L_P/\partial \btheta$, on a neighborhood of $\xi_P$.
It follows that besides the constraint in~\eqref{def:Smin structured},
the pair $(\bxi_P,\lambda_P)$ satisfies
\begin{equation}
\label{eq:M-equations}
\begin{split}
\int
u(d)\bX^T\bV^{-1}(\by-\bX\bbeta)
\,\dd P(\bs)
&=
\mathbf{0}\\
\mathrm{tr}\left(\bV^{-1}\frac{\partial \bV}{\partial \theta_j}\right)
+
\frac{\lambda}{2}
\int
u(d)
(\by-\bX\bbeta)^T\bV^{-1}
\frac{\partial \bV}{\partial \theta_j}
\bV^{-1}(\by-\bX\bbeta)
\,\dd P(\mathbf{s})
&
=0,
\end{split}
\end{equation}
for $j=1,\ldots,l$,
where $u(s)=\rho'(s)/s$ and $d=d(\bs,\bxi)$ is defined by~\eqref{def:Mahalanobis distance structured},
and where we abbreviate~$\bV(\btheta)$ by $\bV$.
To solve $\lambda_P$ from the second set of equations, we multiply the $j$-th equation by~$\theta_j$ and then
sum over $j=1,\ldots,l$.
This leads to
\[
\mathrm{tr}\left(\bV^{-1}\sum_{j=1}^l\theta_j\frac{\partial \bV}{\partial \theta_j}\right)
+
\frac{\lambda}{2}
\int
u(d)
(\by-\bX\bbeta)^T\bV^{-1}
\left(
\sum_{j=1}^l\theta_j\frac{\partial \bV}{\partial \theta_j}
\right)
\bV^{-1}(\by-\bX\bbeta)
\,\dd P(\mathbf{s})
=0,
\]
which is solved by
\[
\lambda_P
=
\frac{-2\mathrm{tr}\left(\bV^{-1}\sum_{j=1}^l\theta_j(\partial \bV/\partial \theta_j)\right)}{
\int
u(d)
(\by-\bX\bbeta)^T\bV^{-1}
\left(
\sum_{j=1}^l\theta_j(\partial \bV/\partial \theta_j)
\right)
\bV^{-1}(\by-\bX\bbeta)
\,\dd P(\mathbf{s})}.
\]
When we insert this back into the second equation in~\eqref{eq:M-equations}, we find
\[
\begin{split}
&
\mathrm{tr}\left(\bV^{-1}\frac{\partial \bV}{\partial \theta_j}\right)
\int
u(d)
(\by-\bX\bbeta)^T\bV^{-1}
\left(
\sum_{t=1}^l\theta_t\frac{\partial \bV}{\partial \theta_t}
\right)
\bV^{-1}(\by-\bX\bbeta)
\,\dd P(\mathbf{s})\\
&\quad-
\mathrm{tr}\left(\bV^{-1}\sum_{t=1}^l\theta_t\frac{\partial \bV}{\partial \theta_t}\right)
\int
u(d)
(\by-\bX\bbeta)^T\bV^{-1}
\frac{\partial \bV}{\partial \theta_j}
\bV^{-1}(\by-\bX\bbeta)
\,\dd P(\mathbf{s})
=0,
\end{split}
\]
or briefly
\begin{equation}
\label{eq:M-equations linear dependent}
\int
u(d)
(\by-\bX\bbeta)^T\bV^{-1}
\bH_j
\bV^{-1}(\by-\bX\bbeta)
\,\dd P(\mathbf{s})
=0,
\quad
j=1,\ldots,l,
\end{equation}
where
\begin{equation}
\label{def:Hj}
\bH_j
=
\mathrm{tr}\left(\bV^{-1}\frac{\partial \bV}{\partial \theta_j}\right)
\left(
\sum_{t=1}^l\theta_t\frac{\partial \bV}{\partial \theta_t}
\right)
-
\mathrm{tr}\left(\bV^{-1}\sum_{t=1}^l\theta_t\frac{\partial \bV}{\partial \theta_t}\right)
\frac{\partial \bV}{\partial \theta_j}.
\end{equation}
Because $\sum_{j=1}^l \theta_j\bH_j=\mathbf{0}$, the system of equations~\eqref{eq:M-equations linear dependent} is linearly dependent.
Similar to~\cite{lopuhaa1989} we subtract the S-constraint from each equation.
For each $j=1,\ldots,l$, we subtract the term
\[
\mathrm{tr}\left(\bV^{-1}\frac{\partial \bV}{\partial \theta_j}\right)(\rho(d)-b_0)
\]
from the left hand side of equation~\eqref{eq:M-equations linear dependent}.
We then find that any solution $\bxi_P$ of~\eqref{def:Smin structured} satisfies
the following equation
\begin{equation}
\label{eq:Psi=0}
\int
\Psi(\bs,\bxi)
\,\dd P(\bs)
=
\mathbf{0},
\end{equation}
where $\Psi=(\Psi_{\bbeta},\Psi_{\btheta})$, with $\Psi_{\btheta}=(\Psi_{\btheta,1},\ldots,\Psi_{\btheta,l})$, where
\begin{equation}
\label{eq:Psi function}
\begin{split}
\Psi_{\bbeta}(\bs,\bxi)
&=
u(d)\bX^T\bV^{-1}(\by-\bX\bbeta)\\
\Psi_{\btheta,j}(\bs,\bxi)
&=
u(d)
(\by-\bX\bbeta)^T\bV^{-1}
\bH_j
\bV^{-1}(\by-\bX\bbeta)
-
\mathrm{tr}\left(\bV^{-1}\frac{\partial \bV}{\partial \theta_j}\right)(\rho(d)-b_0),
\end{split}
\end{equation}
for $j=1,\ldots,l$, where $\bH_j$ and $d=d(\bs,\bxi)$ are defined in~\eqref{def:Hj} and~\eqref{def:Mahalanobis distance structured},
respectively,
and where we abbreviate~$\bV(\btheta)$ by $\bV$.

The regression score equation for $\Psi_{\bbeta}$ with the empirical measure $\mathbb{P}_n$ for $P$ in~\eqref{eq:Psi=0}
coincides with the one for the regression S-estimator in the linear mixed effects model~\eqref{def:linear mixed effects model Copt}
considered in~\cite{copt2006high} (see their equation~(10)).
The empirical regression score equation also coincides with the one
for the regression S-estimator in the multivariate regression model of Example~\ref{ex:multivariate regression} considered in~\cite{vanaelst&willems2005}
(see equation~(2.2) in~\cite{vanaelst&willems2004}).
Similarly, the empirical score equation for $\Psi_{\bbeta}$ coincides with the one
for the location S-estimator of Example~\ref{ex:location-scale} considered in~\cite{lopuhaa1989}.

For general covariance structures the empirical covariance score equation for $\Psi_{\btheta}$ does not compare
directly to existing equations in the literature.
However, as we will see in the next subsection,
similar comparisons are available for models with a linear covariance structure.

\subsection{Linear covariance structures}
\label{subsec:linear V}
In the previous section, we solved $\lambda$ from~\eqref{eq:M-equations} and subtracted the S-constraint,
leading to score equation~\eqref{eq:Psi=0} with $\Psi$ given in~\eqref{eq:Psi function}.
The fact that this was done in a specific way has the following reason.
In cases where~$\bV(\btheta)$ is linear, say
\begin{equation}
\label{def:V linear}
\bV(\btheta)=\sum_{j=1}^l\theta_j\bL_j,
\end{equation}
the function~$\Psi_{\btheta}$ simplifies a lot and can also be related to
the covariance psi-function in~\cite{lopuhaa1989}.
Typical models of interest that have a covariance matrix of this type are
the mixed linear effects model from Example~\ref{ex:linear mixed effects model}
and the multivariate regression model from Example~\ref{ex:multivariate regression}.
But also the multivariate location-scale model from Example~\ref{ex:location-scale}
and the time series model~\eqref{def:autogressive orde 1 covariance} from Example~\ref{ex:autoregressive model}
have linear covariance structures.

When $\bV$ is of the form~\eqref{def:V linear}, then $\partial\bV/\partial\theta_j=\bL_j$
and $\sum_{j=1}^l\theta_j (\partial\bV/\partial\theta_j)=\bV$.
In this case,~\eqref{def:Hj} simplifies to
$\bH_j=
\mathrm{tr}\left(\bV^{-1}\bL_j\right)\bV
-
k\bL_j$,
and $\Psi_{\btheta,j}$ in~\eqref{eq:Psi function} becomes
\[
\Psi_{\btheta,j}(\bs,\bxi)
=
\mathrm{tr}\left(\bV^{-1}\bL_j\right)
v(d)
-
k
u(d)
(\by-\bX\bbeta)^T\bV^{-1}
\bL_j
\bV^{-1}(\by-\bX\bbeta),
\]
where $u(s)$ is defined in (R4) and
\begin{equation}
\label{def:v(d)}
v(s)=u(s)s^2-\rho(s)+b_0.
\end{equation}
Using that $\tr(\bA^T\bB)=\vc(\bA)^T\vc(\bB)$, this can be written as
\[
\Psi_{\btheta,j}(\bs,\bxi)
=
-\vc\Big(
k
u(d)
(\by-\bX\bbeta)(\by-\bX\bbeta)^T
-
v(d)
\mathbf{V}
\Big)^T
\vc\left(
\bV^{-1}
\bL_j
\bV^{-1}
\right).
\]
On the right hand side we recognize
$k
u(d)
(\by-\bX\bbeta)(\by-\bX\bbeta)^T
-
v(d)
\mathbf{V}$,
being the covariance psi-function that also appears in (2.8) in~\cite{lopuhaa1989}.
For our purposes we define
\begin{equation}
\label{def:PsiV}
\Psi_\bV(\bs,\bxi)
=
k
u(d(\bs,\bxi))
(\by-\bX\bbeta)(\by-\bX\bbeta)^T
-
v(d(\bs,\bxi))
\bV.
\end{equation}
The functions $\Psi_{\btheta,j}$, for $j=1,\ldots,l$, can be combined in one expression for the
vector valued function $\Psi_{\btheta}$ as follows.
First note that
\[
\vc\left(
\bV^{-1}
\bL_j
\bV^{-1}
\right)
=
\left(
\bV^{-1}
\otimes
\bV^{-1}
\right)
\vc\left(
\bL_j
\right)
\]
for $j=1,\ldots,l$.
Define the $k^2\times l$ matrix
\begin{equation}
\label{def:L}
\bL
=
\Big[
  \begin{array}{ccc}
    \vc\left(\bL_1\right)
 & \cdots &     \vc\left(\bL_l\right) \\
  \end{array}
\Big].
\end{equation}
Then, the column vector $\Psi_{\btheta}=(\Psi_{\btheta,1},\ldots,\Psi_{\btheta,l})$ can be written as
\[
\Psi_{\btheta}(\bs,\bxi)
=
-\bL^T
\left(
\bV^{-1}
\otimes
\bV^{-1}
\right)
\vc\left(
\Psi_\bV(\bs,\bxi)
\right),
\]
where $\Psi_\bV$ is defined in~\eqref{def:PsiV} and $\bL$ in~\eqref{def:L}.
Note that the dependence on $\bs=(\by,\bX)$ in $\Psi_{\btheta}$ is only through the function $\Psi_\bV$.
We conclude that in the case of a linear covariance structure,
any solution~$\bxi_P$ of~\eqref{def:Smin structured} satisfies~\eqref{eq:Psi=0},
where $\Psi=(\Psi_{\bbeta},\Psi_{\btheta})$, with
\begin{equation}
\label{eq:Psi linear}
\begin{split}
\Psi_{\bbeta}(\bs,\bxi)
&=
u(d)\bX^T\bV^{-1}(\by-\bX\bbeta)\\
\Psi_{\btheta}(\bs,\bxi)
&=
-\bL^T
\left(
\bV^{-1}
\otimes
\bV^{-1}
\right)
\vc\left(
\Psi_\bV(\bs,\bxi)
\right)\end{split}
\end{equation}
where $d=d(\bs,\bxi)$ is defined in~\eqref{def:Mahalanobis distance structured},
and where we abbreviate~$\bV(\btheta)$ by $\bV$.

For the multivariate regression model in Example~\ref{ex:multivariate regression},
one has $\bV(\btheta)=\bC$, where
$\btheta=\vch(\bC)$.
The matrix $\bL=\partial\vc(\bV)/\partial\btheta^T$ is then equal to the so-called duplication matrix~$\mathcal{D}_k$,
which is the unique $k^2\times k(k+1)/2$ matrix, with the properties
$\mathcal{D}_k\vch(\bC)=\vc(\bC)$ and $(\mathcal{D}_k\mathcal{D}_k)^{-1}\mathcal{D}_k^T\vc(\bC)=\vch(\bC)$
(e.g., see~\cite[Ch.~3, Sec.~8]{magnus&neudecker1988}).
Because~$\bV$ has full rank, it follows that equation~\eqref{eq:Psi=0} holds for $\Psi=(\Psi_{\bbeta},\Psi_{\bV})$.
The resulting score equations for the empirical measure $\mathbb{P}_n$ corresponding to observations
$(\by_i,\bX_i)$, for $i=1,\ldots,n$, are then equivalent with the ones found in~\cite{vanaelst&willems2005}.

For the linear mixed effects model~\eqref{def:linear mixed effects model Copt},
the covariance matrix $\bV(\btheta)$ has a linear structure with the vector
$\btheta=(\sigma_0^2,\ldots,\sigma_{r}^2)$
of unknown covariance parameters.
The matrix~$\bL$ is then a $k^2\times (r+1)$ matrix and will typically be of rank $r+1<k^2$.
As a consequence, in this case one cannot further simplify equation~\eqref{eq:Psi=0},
by removing the factor $\bL^T\left(\bV^{-1}\otimes\bV^{-1}\right)$ from the function~$\Psi_{\btheta}$.
The score equation for $\Psi_{\bbeta}$ resulting from the empirical measure $\mathbb{P}_n$
corresponding to observations $(\by_i,\bX_i)$, for $i=1,\ldots,n$, is the same as the one
obtained in~\cite{copt2006high}.
The corresponding score equation for $\Psi_{\btheta}$ differs slightly from the one in~\cite{copt2006high},
because the authors do not subtract a term with $\rho(d)-b_0$
to remove the linear dependency of the equations~\eqref{eq:M-equations linear dependent}.

\section{Local robustness: the influence function}
\label{sec:IF}
For $0<h<1$ and $\bs=(\by,\bX)\in\R^k\times\R^{kq}$ fixed, define the perturbed probability measure
\[
P_{h,\bs}=(1-h)P+h\delta_{\bs},
\]
where $\delta_{\bs}$ denotes the Dirac measure at $\bs\in\R^k\times\R^{kq}$.
The \emph{influence function} of the functional~$\bxi(\cdot)$ at probability measure $P$,
is defined as
\begin{equation}
\label{def:IF}
\text{IF}(\bs;\bxi,P)
=
\lim_{h\downarrow0}
\frac{\bxi((1-h)P+h\delta_{\bs})-\bxi(P)}{h},
\end{equation}
if this limit exists.
In contrast to the global robustness measured by the breakdown point,
the influence function measures the local robustness.
It describes the effect of an infinitesimal contamination at a single point $\bs$ on the functional
(see Hampel~\cite{hampel1974}).
Good local robustness is therefore illustrated by a bounded influence function.

\subsection{The general case}
\label{subsec:IF general}
The theorem below gives the influence function for the S-functional~$\bxi$.
It extends the result for S-functionals of multivariate location and scale~\cite{lopuhaa1989}.
Under the assumption that the limit in~\eqref{def:IF} exists and $P$ has an elliptical contoured density~\eqref{eq:elliptical},
Van Aelst and Willems~\cite{vanaelst&willems2005}
relate the influence function for S-functionals of multivariate regression
to that of S-functionals of multivariate location and scale.
For the linear mixed effects model considered in~\cite{copt2006high},
the influence function has not been established.
The influence function for these functionals now follows as a special case from the theorem below.

We will show that the limit in~\eqref{def:IF} exists and derive its expression at general $P$.
Since the value of $\btheta$ determines the covariance matrix~$\bV(\btheta)$,
we also include the influence function of the covariance functional.
Consider the S-functional at $P_{h,\bs_0}$.
From the Portmanteau theorem~\cite[Theorem 2.1]{billingsley1968} it can easily be seen that $P_{h,\bs_0}\to P$, weakly, as~$h\downarrow0$.
Therefore, under the conditions of Corollary~\ref{cor:existence weak convergence} and Theorem~\ref{th:continuity structure}, it follows
that there exist solutions $\bxi(P_{h,\bs_0})$ and $\bxi(P)$ to minimization problems~\eqref{def:Smin structured}
at $P_{h,\bs_0}$ and $P$, respectively, and that $\bxi(P_{h,\bs_0})\to\bxi(P)$, as $h\downarrow0$.
\begin{theorem}
\label{th:IF}
Let $\bxi(P_{h,\bs_0})$ and $\bxi(P)$ be solutions to minimization problems~\eqref{def:Smin structured}
at~$P_{h,\bs_0}$ and $P$, respectively, and suppose that $\bxi(P_{h,\bs_0})\to\bxi(P)$, as $h\downarrow0$.
Suppose that $\rho$ satisfies (R4) and $\bV$ satisfies~(V4).
Let $\Psi$ be defined in~\eqref{eq:Psi function} and suppose that
\begin{equation}
\label{def:Lambda}
\Lambda(\bxi)
=
\int \Psi(\bs,\bxi)\,\dd P(\bs),
\end{equation}
is continuously differentiable with a non-singular derivative $\mathbf{D}(P)$ at $\bxi(P)$.
Then for $\bs_0\in\R^k\times\R^{kq}$,
\[
\mathrm{IF}(\bs_0;\bxi,P)
=
-\mathbf{D}(P)^{-1}\Psi(\bs_0,\bxi(P)).
\]
For the covariance functional $\bC(P)=\bV(\btheta(P))$, it holds that
\[
\mathrm{IF}(\bs_0;\vc(\bC),P)
=
\left(\frac{\partial\, \vc(\bV(\btheta(P)))}{\partial\btheta^T}\right)
\mathrm{IF}(\bs_0;\btheta,P).
\]
\end{theorem}
To investigate the local robustness of S-estimators, we derive the following bound
on the influence function for $\bxi(P)$.
\begin{corollary}
\label{cor:IF bounded}
Suppose that $\rho$ satisfies (R2) and (R4), and $\bV$ satisfies (V4).
Then there exist $0<C_1<\infty$ and $0<C_2<\infty$, only depending on $P$, such that
for $\bs=(\by,\bX)$ it holds that $\|\mathrm{IF}(\bs,\bxi(P))\|\leq C_1+C_2\|\bX\|$.
\end{corollary}
Its proof can be found in~\cite{supplement}.

\subsection{Elliptically contoured densities}
\label{subsec:IF elliptical}
When $P$ is such that $\by\mid\bX$ has an elliptically contoured density~\eqref{eq:elliptical}
and $\bV(\btheta)$ is linear,
we can obtain a more detailed expression for the influence function,
This requires the following condition on the function $\rho$,
\begin{itemize}
\item[(R5)]
$\rho$ is twice continuously differentiable,
\end{itemize}
and the following condition on the mapping $\btheta\mapsto\bV(\btheta)$,
\begin{itemize}
\item[(V5)]
$\bV(\btheta)$ is twice continuously differentiable.
\end{itemize}
Conditions (R5) and (V5) are needed to establish that $\Lambda$, as defined in~\eqref{def:Lambda},
is continuously differentiable.
Clearly, condition (V5) implies former conditions (V4) and (V1).

Suppose that $P$ is such that $\by\mid\bX$ has an elliptically contoured density $f_{\bmu,\bSigma}$ from~\eqref{eq:elliptical},
with $\bmu\in\R^k$ and $\bSigma\in\text{PDS}(k)$.
When the S-functional is affine equivariant,
it suffices to determine the influence function for the case $(\bmu,\bSigma)=(\mathbf{0},\bI_k)$.
However, this does not hold in general for the S-functionals in our setting.
The reason is that, for a $k\times k$ non-singular matrix~$\bA$ and $\btheta\in\R^l$,
the matrix
$\bA\bV(\btheta)\bA^T$ may not be of the form $\bV(\btheta')$, for some $\btheta'\in\R^l$.
Examples are the (linear)
covariance structure that corresponds to
the linear mixed effects model~\eqref{def:linear mixed effects model Copt}
considered in~\cite{copt2006high}
or the models discussed in Example~\ref{ex:autoregressive model}.

Nevertheless, note that for the general case with $\bmu\in\R^k$ and $\bSigma\in\text{PDS}(k)$,
we can still use the fact that, conditionally on $\bX$, the distribution of $\by$ is the same as that of
$\bSigma^{1/2}\bz+\bmu$, where~$\bz$ has a spherical density $f_{\mathbf{0},\bI_k}$.
As a consequence, we can still obtain the following result, which
enables one to determine the influence functions of the functionals~$\bbeta(P)$
and~$\btheta(P)$ separately.

If $P$ itself is also absolutely continuous, then it satisfies (C3), as well as $(\text{C1}_{\epsilon'})$ and $(\text{C2}_\epsilon)$,
for any $0<\epsilon'<\epsilon\leq 1-r$.
When~$\rho$ and $\bV$ satisfy (R1)-(R3) and (V1)-(V3), it follows from Theorem~\ref{th:existence structured} and
Corollary~\ref{cor:existence weak convergence}
that~$\bxi(P)$ and $\bxi(P_{h,\bs})$ exist, for $h$ sufficiently small.
If $h$ in~\eqref{eq:elliptical} is non-increasing and not constant on $[0,c_0^2]$, then
$\bxi(P)$ is unique, according to Theorem~\ref{th:davies},
so that $\bxi(P_{h,\bs})\to\bxi(P)$, as $h\downarrow0$.
Hence, in order to apply Theorem~\ref{th:IF}, it remains to show that $\Lambda$ in~\eqref{def:Lambda}
is continuously differentiable with a non-singular derivative at~$\bxi(P)$.
As a first step we obtain that the derivative of $\Lambda$ is a block matrix.
\begin{lemma}
\label{lem:block derivative}
Suppose that $P$ is
such that~$\by\mid\bX$ has an elliptically contoured
density~$f_{\bmu,\bSigma}$ from~\eqref{eq:elliptical}
and $\E\|\bX\|^2<\infty$.
Suppose that $\bxi(P)$ is a solution to the corresponding minimization problem~\eqref{def:Smin structured},
such that $(\bX\bbeta(P),\bV(\btheta(P)))=(\bmu,\bSigma)$.
Suppose that $\rho$ satisfies (R2), (R4)-(R5) and that $\bV$ satisfies (V5) and
has a linear structure~\eqref{def:V linear}.
Let $\Lambda$ be defined in~\eqref{def:Lambda} with $\Psi$ defined in~\eqref{eq:Psi linear}.
Then
\[
\frac{\partial \Lambda(\bxi(P))}{\partial \bxi}
=
\left(
  \begin{array}{cc}
\dfrac{\partial \Lambda_{\bbeta}(\bxi(P))}{\partial \bbeta} & \mathbf{0} \\
    \\
\mathbf{0} & \dfrac{\partial \Lambda_{\btheta}(\bxi(P))}{\partial \btheta} \\
  \end{array}
\right),
\]
where
\begin{equation}
\label{def:derivative Lambda beta}
\frac{\partial\Lambda_{\bbeta}(\bxi(P))}{\partial \bbeta}
=
-\alpha
\mathbb{E}\left[\mathbf{X}^T\bSigma^{-1}\mathbf{X}\right],
\end{equation}
with
\begin{equation}
\label{def:alpha}
\alpha
=
\mathbb{E}_{\mathbf{0},\bI_k}
\left[
\left(1-\frac{1}{k}\right)
\frac{\rho'(\|\bz\|)}{\|\bz\|}
+
\frac1k
\rho''(\|\bz\|)
\right],
\end{equation}
and
\[
\frac{\partial\Lambda_{\btheta}(\bxi(P))}{\partial \btheta}
=
\gamma_1\bL^T\left(\bSigma^{-1}\otimes\bSigma^{-1}\right)\bL
-
\gamma_2\bL^T
\vc(\bSigma^{-1})
\vc(\bSigma^{-1})^T
\bL.
\]
where $\bL=\partial\vc(\bV(\btheta(P)))/\partial\btheta^T$ is the $k^2\times l$ matrix given in~\eqref{def:L}
and
\begin{equation}
\label{def:gamma12}
\begin{split}
\gamma_1
&=
\frac{\mathbb{E}_{0,\mathbf{I}_k}
\left[
\rho''(\|\bz\|)\|\bz\|^2+(k+1)\rho'(\|\bz\|)\|\bz\|
\right]}{k+2}\\
\gamma_2
&=
\frac{\mathbb{E}_{0,\mathbf{I}_k}
\left[
2\rho''(\|\bz\|)\|\bz\|^2
+
k\rho'(\|\bz\|)\|\bz\|
\right]}{2k(k+2)},
\end{split}
\end{equation}
\end{lemma}
The proof is tedious, but straightforward, and can be found in~\cite{supplement}.
\begin{remark}
\label{rem:interchange int and diff}
The proof of Lemma~\ref{lem:block derivative} uses the fact
that
\[
\frac{\partial\Lambda(\bxi)}{\partial \bxi}
=
\int
\frac{\partial\Psi(\bs,\bxi)}{\partial \bxi}\,\dd P(\bs).
\]
for all $\bxi$ in a neighborhood of $\bxi(P)$.
This holds for general $P$ and any covariance structure~$\bV(\btheta)$ that satisfies (V2)-(V3) and (V5),
see Lemma~\ref{lem:diff Lambda} in~\cite{supplement}.
Furthermore, Lemma~\ref{lem:block derivative} is obtained for a linear covariance structure.
However, with some additional technicalities, this result can also be shown to hold for $\Psi$ defined in~\eqref{eq:Psi function}
corresponding to general covariance structures.
For general covariance structures one still obtains~\eqref{def:derivative Lambda beta},
and that
\[
\begin{split}
\frac{\partial\Lambda_{\btheta,j}(\bxi(P))}{\partial\theta_s}
=
&-\alpha_1
\tr\left(
\bSigma^{-1}\frac{\partial\bV(\btheta(P))}{\partial\theta_s}\bSigma^{-1}\bH_j
\right)\\
&+
\alpha_2
\mathrm{tr}
\left(\bSigma^{-1}\frac{\partial \bV(\btheta(P))}{\partial \theta_s}\right)
\mathrm{tr}
\left(\bSigma^{-1}\frac{\partial \bV(\btheta(P))}{\partial \theta_j}\right),
\end{split}
\]
for $j,s=1,\ldots,l$, and where
$\alpha_1=\gamma_1/k$ and $\alpha_2=\gamma_1/k-\gamma_2$,
with $\gamma_1,\gamma_2$ from~\eqref{def:gamma12},
and where $\bH_j$ is defined in~\eqref{def:Hj}.
\end{remark}
The next corollary gives expressions for the influence functions of the functionals~$\bbeta(P)$
and~$\btheta(P)$ separately, at a distribution $P$ that is such that
$\by\mid\bX$ has an elliptically contoured density.
The proof is tedious, but straightforward, and can be found in~\cite{supplement}.
\begin{corollary}
\label{cor:IF elliptical}
Suppose that $P$ is such that $\by\mid\bX$ has an elliptically contoured
density~$f_{\bmu,\bSigma}$ from~\eqref{eq:elliptical},
such that $(\bX\bbeta(P),\bV(\btheta(P)))=(\bmu,\bSigma)$.
Let $\bxi(P_{h,\bs_0})$ and $\bxi(P)$ be a solution to minimization problem~\eqref{def:Smin structured}
at $P_{h,\bs_0}$ and $P$, respectively, and suppose that $\bxi(P_{h,\bs_0})\to\bxi(P)$, as $h\downarrow0$.
Suppose that $\E\|\bX\|^2<\infty$ and suppose that~$\rho$ satisfies (R2)-(R5) and that $\bV$ satisfies (V5), and
has a linear structure~\eqref{def:V linear}.
Let~$\alpha$,~$\gamma_1$, and~$\gamma_2$ be defined in~\eqref{def:alpha} and~\eqref{def:gamma12},
and suppose that $\mathbb{E}_{0,\bI_k}
\left[
\rho''(\|\bz\|)
\right]>0$.
If $\bX$ has full rank with probability one, then
\[
\mathrm{IF}(\bs_0,\bbeta,P)
=
\frac{u(d_0)}{\alpha}
\Big(
\E\left[\bX^T\bSigma^{-1}\bX\right]
\Big)^{-1}
\bX_0^T\bSigma^{-1}(\by_0-\bX_0\bbeta(P))
\]
where $d_0^2=(\by_0-\bX_0\bbeta(P))^T\bSigma^{-1}(\by_0-\bX_0\bbeta(P))$ and $u(s)=\rho'(s)/s$.
If $\gamma_1>0$ and the $k^2\times l$ matrix $\bL$, as defined in~\eqref{def:L}, has full rank,
then $\mathrm{IF}(\bs_0,\btheta,P)$ is given by
\[
\begin{split}
&
\frac{ku(d_0)}{\gamma_1}
\Big(\bL^T(\bSigma^{-1}\otimes\bSigma^{-1})\bL)\Big)^{-1}
\bL^T
\vc\left(\bSigma^{-1}(\by_0-\bX_0\bbeta(P))(\by_0-\bX_0\bbeta(P))^T\bSigma^{-1}\right)\\
&\qquad+
\left(
-\frac{u(d_0)d_0^2}{\gamma_1}
+
\frac{\rho(d_0)-b_0}{\gamma_1-k\gamma_2}
\right)
\btheta(P).
\end{split}
\]
\end{corollary}
Note that since $\bL\btheta(P)=\vc(\bV(\btheta(P)))=\vc(\bSigma)$,
we can immediately obtain the influence function for the covariance functional~$\bC(P)=\bV(\btheta(P))$.
From Theorem~\ref{th:IF} it immediately follows that $\mathrm{IF}(\bs_0,\vc(\bC),P)$ is given by
\[
\begin{split}
&
\frac{ku(d_0)}{\gamma_1}
\bL\Big(\bL^T(\bSigma^{-1}\otimes\bSigma^{-1})\bL)\Big)^{-1}
\bL^T
\vc\left(\bSigma^{-1}(\by_0-\bX_0\bbeta)(\by_0-\bX_0\bbeta)^T\bSigma^{-1}\right)\\
&\qquad+
\left(
-\frac{u(d_0)d_0^2}{\gamma_1}
+
\frac{\rho(d_0)-b_0}{\gamma_1-k\gamma_2}
\right)
\vc(\bSigma).
\end{split}
\]
Since the functions $u(s)s=\rho'(s)$,
$u(s)s^2=\rho'(s)s$, and $\rho(s)$ are bounded, it follows that
$\mathrm{IF}(\bs,\btheta,P)$ and $\mathrm{IF}(\bs,\vc(\bC),P)$ are bounded uniformly in both $\by$ and $\bX$,
whereas $\mathrm{IF}(\bs,\bbeta,P)$
is bounded uniformly in $\by$, but not in $\bX$.
This illustrates the phenomenon in linear regression that
leverage points can have a high effect on the regression S-estimator.

For the S-estimators in the linear mixed effects model~\eqref{def:linear mixed effects model Copt} with normal errors
considered in~\cite{copt2006high}, the influence function is not available.
The expression can now be obtained from Corollary~\ref{cor:IF elliptical}.
The expression for $\mathrm{IF}(\bs,\bbeta,P)$ in Corollary~\ref{cor:IF elliptical} coincides with the one found
for the multivariate regression S-functional in~\cite{vanaelst&willems2005},
where $\alpha>0$ is the same constant as the one in the expression
of the influence function for the location S-functional in~\cite{lopuhaa1989}.
Furthermore, for the multivariate regression model, one has
$\btheta=\vch(\bC)$ and the matrix $\bL$ is equal to the duplication matrix~$\mathcal{D}_k$.
From the properties of $\mathcal{D}_k$, the expressions for the influence functions simplify.
One finds in this case that
\[
\mathrm{IF}(\bs,\btheta,P)
=
\frac{ku(d)}{\gamma_1}
\vch\left((\by-\bX\bbeta(P))(\by-\bX\bbeta(P))^T\right)
+
\left(
-\frac{u(d)d^2}{\gamma_1}
+
\frac{\rho(d)-b_0}{\gamma_1-k\gamma_2}
\right)
\btheta(P)
\]
and the influence function of the covariance functional $\bC(P)=\bV(\btheta(P))$ itself
is given by
\[
\mathrm{IF}(\bs,\bC,P)
=
\frac{ku(d)}{\gamma_1}
(\by-\bX\bbeta(P))(\by-\bX\bbeta(P))^T
+
\left(
-\frac{u(d)d^2}{\gamma_1}
+
\frac{\rho(d)-b_0}{\gamma_1-k\gamma_2}
\right)
\bSigma.
\]
This coincides with the expressions found for the covariance S-functionals
in~\cite{vanaelst&willems2005} and in~\cite{lopuhaa1989}.

\section{Asymptotic normality}
\label{sec:asymptotic normality}
To establish asymptotic normality of the S-estimators, we use the score equations
obtained from differentiation of the Lagrangian corresponding to the minimization problem~\eqref{def:Smin estimator structured}.
In the same way as before, we obtain score equation~\eqref{eq:Psi=0}, with~$P$ equal to the empirical measure $\mathbb{P}_n$
corresponding to observations
$\bs_1,\ldots,\bs_n$, with $\bs_i=(\by_i,\bX_i)\in\R^k\times\R^{kq}$.
From~\eqref{eq:Psi=0}, we see that any solution $\bxi_n=\bxi(\mathbb{P}_n)$ to the
S-minimization problem~\eqref{def:Smin estimator structured} must satisfy
\begin{equation}
\label{eq:M-equation estimator}
\int
\Psi(\mathbf{s},\bxi_n)\,\dd \mathbb{P}_n(\mathbf{s})
=
\mathbf{0},
\end{equation}
where $\Psi=(\Psi_{\bbeta},\Psi_{\btheta})$ is defined in~\eqref{eq:Psi function}.
\subsection{General case}
Writing $\bxi_P=\bxi(P)$, we decompose~\eqref{eq:M-equation estimator} as follows
\begin{equation}
\label{eq:decomposition estimator}
\begin{split}
0=
\int \Psi(\mathbf{s},\bxi_n)\,\dd P(\mathbf{s})
&+
\int \Psi(\mathbf{s},\bxi_P)\,\dd (\mathbb{P}_n-P)(\mathbf{s})\\
&+
\int
\left(
\Psi(\mathbf{s},\bxi_n)-\Psi(\mathbf{s},\bxi_P)
\right)
\,\dd (\mathbb{P}_n-P)(\mathbf{s}).
\end{split}
\end{equation}
The essential step in establishing asymptotic normality of $\bxi_n$,
is to show that the third term on the right hand side of~\eqref{eq:decomposition estimator} is of the order $o_P(n^{-1/2})$.
To this end we will apply results from empirical process theory as developed in Pollard~\cite{pollard1984}.
This leads to the following theorem.
\begin{theorem}
\label{th:asymp normal}
Suppose that $\rho$ satisfies (R1)-(R2) and (R4), such that $u(s)$ is of bounded variation,
and suppose that~$\bV$ satisfies (V4).
Let $\bxi_n$ and~$\bxi(P)$ be solutions to minimization problems~\eqref{def:Smin estimator structured}
and~\eqref{def:Smin structured}, and suppose that $\bxi_n\to\bxi(P)$ in probability.
Suppose that $\Lambda$, as defined in~\eqref{def:Lambda} with~$\Psi$ defined in~\eqref{eq:Psi function},
is continuously differentiable with a non-singular derivative $\bD(P)$ at $\bxi(P)$
and suppose that $\E\|\mathbf{X}\|^2<\infty$.
Then $\sqrt{n}(\bxi_n-\bxi(P))$ is asymptotically normal with mean zero and
covariance matrix
\[
\bD(P)^{-1}
\E
\left[
\Psi(\mathbf{s},\bxi(P))
\Psi(\mathbf{s},\bxi(P))^T
\right]
\bD(P)^{-1}.
\]
\end{theorem}
Theorem~\ref{th:asymp normal} is similar to Theorem~4.1 in~\cite{lopuhaa1989}.
Note that Theorem~\ref{th:asymp normal} confirms the
well know heuristic that relates the limiting covariance  of $\sqrt{n}(\bxi_n-\bxi(P))$
to the influence function of the functional $\bxi(\cdot)$ given in Theorem~\ref{th:IF},
\begin{equation}
\label{eq:IF heuristic}
\bD(P)^{-1}
\E
\left[
\Psi(\mathbf{s},\bxi(P))
\Psi(\mathbf{s},\bxi(P))^T
\right]
\bD(P)^{-1}
=
\E
\left[
\text{IF}(\bs,\bxi,P)\text{IF}(\bs,\bxi,P)^T
\right].
\end{equation}
Van Aelst and Willems~\cite{vanaelst&willems2005} consider the limiting behavior of S-estimators in the multivariate regression model
of Example~\ref{ex:multivariate regression}, but only under $P$ for which $\by\mid\bX$ has an elliptical contoured density.
Copt and Victoria-Feser~\cite{copt2006high} consider asymptotic normality for S-estimators in the
linear mixed effects model~\eqref{def:linear mixed effects model Copt} with a constant design matrix $\bX_i=\bX$
and only consider $P$ for which $\by\mid\bX$ has an multivariate normal distribution.

\subsection{Elliptically contoured densities}
\label{subsec:Asymp Norm elliptical}
Consider the special case that $P$ is such that $\by\mid\bX$ has an elliptically contoured density $f_{\bmu,\bSigma}$ from~\eqref{eq:elliptical},
with $\bmu\in\R^k$ and $\bSigma\in\text{PDS}(k)$.
As before, in determining the limiting normal distribution of the individual S-estimators,
we cannot use affine equivariance and restrict ourselves to the case $(\mathbf{0},\bI_k)$.
Instead, we use some of the results obtained in Section~\ref{subsec:IF elliptical}
to establish the limiting normal distributions of the S-estimators $\bbeta_n=\bbeta(\mathbb{P}_n)$,
$\btheta_n=\btheta(\mathbb{P}_n)$, and $\bC_n=\bV(\btheta(\mathbb{P}_n))$.

\begin{corollary}
\label{cor:Asymp norm elliptical}
Suppose that $P$ is such that $\by\mid\bX$ has an elliptically contoured
density~$f_{\bmu,\bSigma}$ from~\eqref{eq:elliptical},
such that $(\bX\bbeta(P),\bV(\btheta(P)))=(\bmu,\bSigma)$.
Let $\bxi_n$ and~$\bxi(P)$ be solutions to minimization problems~\eqref{def:Smin estimator structured}
and~\eqref{def:Smin structured}, and suppose that $\bxi_n\to\bxi(P)$ in probability.
Suppose that $\E\|\bX\|^2<\infty$ and suppose that~$\rho$ satisfies (R2)-(R5), such that $u(s)$ is of bounded variation.
Suppose that $\bV$ satisfies (V5), and has a linear structure~\eqref{def:V linear}.
Let~$\alpha$,~$\gamma_1$, and~$\gamma_2$ be defined in~\eqref{def:alpha} and~\eqref{def:gamma12},
and suppose that $\mathbb{E}_{0,\bI_k}
\left[
\rho''(\|\bz\|)
\right]>0$.
If $\bX$ has full rank with probability one, then
$\sqrt{n}(\bbeta_n-\bbeta(P))$ is asymptotically normal
with mean zero and covariance matrix
\[
\frac{\E_{\mathbf{0},\bI_k}\left[\rho'(\|\bz\|)^2\right]}{k\alpha^2}
\left(
\mathbb{E}\left[\mathbf{X}^T\bSigma^{-1}\mathbf{X}\right]
\right)^{-1}.
\]
If $\gamma_1>0$ and the $k^2\times l$ matrix $\bL$, as defined in~\eqref{def:L}, has full rank,
then $\sqrt{n}(\btheta_n-\btheta(P))$ is asymptotically normal with mean zero and covariance matrix
\[
2\sigma_1\Big(\bL^T\left(\bSigma^{-1}\otimes\bSigma^{-1}\right)\bL\Big)^{-1}
+
\sigma_2
\btheta(P)\btheta(P)^T,
\]
where
\[
\begin{split}
\sigma_1
&=
\frac{k(k+2)\E_{\mathbf{0},\bI_k}\left[u(\|\bz\|)^2\|\bz\|^4\right]}{
\left(
\mathbb{E}_{0,\mathbf{I}_k}
\left[
\rho''(\|\bz\|)\|\bz\|^2+(k+1)\rho'(\|\bz\|)\|\bz\|
\right]
\right)^2}\\
\sigma_2
&=
-\frac2k\sigma_1
+
\frac{4\E_{\mathbf{0},\bI_k}[\left(\rho(\|\bz\|)-b_0\right)^2]}{
\left(\E_{\mathbf{0},\bI_k}\left[\rho'(\|\bz\|)^2\right]\right)^2}
\end{split}
\]
\end{corollary}

Due to the linearity of $\bV$,
we can immediately establish asymptotic normality of the covariance
estimator $\bC_n=\bV(\btheta_n)$.
From Corollary~\ref{cor:Asymp norm elliptical} it follows that
\[
\sqrt{n}
\left(\vc(\bC_n)-\vc(\bSigma)\right)
=
\sqrt{n}
\left(
\bL\btheta_n-\bL\btheta(P)
\right)
=
\bL
\sqrt{n}
\left(\btheta_n-\btheta(P)\right).
\]
It follows that the limiting covariance of
$\sqrt{n}\left(\vc(\bV(\btheta_n))-\vc(\bSigma)\right)$ is given by
\[
2\sigma_1
\bL\Big(\bL^T\left(\bSigma^{-1}\otimes\bSigma^{-1}\right)\bL\Big)^{-1}\bL^T
+
\sigma_2
\vc(\bSigma)\vc(\bSigma)^T.
\]

Corollary~\ref{cor:Asymp norm elliptical} is a direct consequence of Theorem~\ref{th:asymp normal}.
Its proof, in particular the derivations of the expressions for the limiting covariances, can be found in~\cite{supplement}.
Note that the constants
$\E_{\mathbf{0},\bI_k}\left[\rho'(\|\bz\|)^2\right]/(k\alpha^2)$, $\sigma_1$ and $\sigma_2$,
are the same as the ones found in~\cite{lopuhaa1989} for the location and covariance S-estimators, respectively.
In fact, Corollary~\ref{cor:Asymp norm elliptical} is an extension of Corollary~5.1 in~\cite{lopuhaa1989} for S-estimators
in the multivariate location-scale model of Example~\ref{ex:location-scale}.

Asymptotic normality of S-estimators in the multivariate regression model of Example~\ref{ex:multivariate regression}
follows from Corollary~\ref{cor:Asymp norm elliptical}.
These estimators have been considered in~\cite{vanaelst&willems2005},
but asymptotic normality has not been established.
Under the assumption that the heuristic~\eqref{eq:IF heuristic} holds,
asymptotic relative efficiencies are computed on the basis of this heuristic.
Indeed, now that Corollary~\ref{cor:Asymp norm elliptical} has been established,
one may check that~\eqref{eq:IF heuristic} holds.

Finally, note that the limiting covariances of $\sqrt{n}(\bbeta_n-\bbeta(P))$ and $\sqrt{n}(\btheta_n-\btheta(P))$
in Corollary~\ref{cor:Asymp norm elliptical}
differ from the ones found in~\cite{copt2006high} for the linear mixed effects model~\eqref{def:linear mixed effects model Copt}
with $\bX_i=\bX$, for $i=1,\ldots,n$.
The results in~\cite{copt2006high} are obtained by re-parameterizing $\bX\bbeta=\bmu$ and interpreting
the model as a multivariate location-scale model.
Then building on the results in~\cite{lopuhaa1989} for S-estimators of multivariate location-scale,
the limiting covariances in~\cite{copt2006high} are found by application of the delta method.
However, in view of Remark~\ref{rem:location-scale} this does not seem to be a correct approach.

\begin{remark}
\label{rem:copt examples}
Although our expressions for the limiting covariances in Corollary~\ref{cor:Asymp norm elliptical}
differ from the ones found in Proposition~1 in~\cite{copt2006high}, somewhat surprisingly,
they yield the same matrices for the  example discussed in Section~5.1 in~\cite{copt2006high}.
However, this is a consequence of the specific structure of the design matrices $\bX$ and $\bZ$ in this example.
One can easily find other design matrices for which the limiting covariances in Corollary~\ref{cor:Asymp norm elliptical}
yield different matrices as the ones found in~\cite{copt2006high}.
Moreover, the corresponding confidence regions based on the expressions in Corollary~\ref{cor:Asymp norm elliptical}
can be substantially smaller than the ones based on the expressions found in~\cite{copt2006high}.
See the simulation in Section~\ref{sec:simulation}.
\end{remark}

\section{Simulation and data example}
\label{sec:simulation}
We compare the asymptotic results of the S-estimators with their finite sample behavior by means of a simulation.
Moreover we investigate the differences between the expressions found in Corollary~\ref{cor:Asymp norm elliptical}
and the ones in Copt and Victoria-Feser~\cite{copt2006high}.
To this end we will study the behavior of the estimators for samples generated from a model that is close to the one in~\cite{copt2006high}:
\begin{equation}
\label{eq:model copt}
\by_i=\bX\bbeta+\gamma_i\bZ+\beps_i,
\quad
i=1,\ldots,n,
\end{equation}
a linear mixed effects model with $\by_i$ in dimension $k=4$ and all subjects with the same design matrix $\bX$ for the fixed effects
$\bbeta=(\beta_1,\beta_2)^T$.
Following the setup in~\cite{copt2006high}, the matrix $\bX$ is built as follows.
The first column of $\bX$ is taken to be a vector $\mathbf{1}$ consisting of ones of length four.
The four $x$-values in the second column are generated from a standard normal,
and then $\bX$ is rescaled to a new matrix $\bX=[\mathbf{1}\quad \bx]$, such that $\bX^T\bX=4\bI_2$.
For our simulation we used
\[
\bX=\left(
  \begin{array}{rr}
1 & -0.9504967\\
1 & -0.5428346\\
1 & 1.6650521\\
1 & -0.1717207\\
\end{array}
\right).
\]
The random effects $\gamma_i$ are independent $N(0,\sigma_\gamma^2)$ distributed random variables, which are independent from the
measurement error $\beps_i\sim N(\mathbf{0},\sigma_\epsilon^2\bR)$.
This leads to a structured covariance $\bSigma=\sigma_\gamma^2\bZ\bZ^T+\sigma_\epsilon^2\bR$,
with covariance parameter vector $\btheta=(\theta_1,\theta_2)^T$, where $\theta_1=\sigma_\gamma^2$
and $\theta_2=\sigma_\epsilon^2$.
Following the setup in~\cite{copt2006high}, we set $\beta_1=\beta_2=1$ and $\theta_1=\theta_2=1$.

In~\cite{copt2006high}, the authors took $\bZ=(1,1,1,1)^T$ and $\bR=\bI_4$.
With these choices the expression
\begin{equation}
\label{eq:varbeta copt}
\text{Var}_{\mathrm{CVF}}(\bbeta_n)
=
\frac{\E_{\mathbf{0},\bI_k}\left[\rho'(\|\bz\|)^2\right]}{k\alpha^2}
(\bX^T\bX)^{-1}\bX^T\bSigma\bX(\bX^T\bX)^{-1}
\end{equation}
found in~\cite{copt2006high} for the limiting covariance matrix of $\sqrt{n}(\bbeta_n-\bbeta)$
(see (14) in~\cite{copt2006high}),
is equal to our expression
\begin{equation}
\label{eq:varbeta our}
\text{Var}_{\mathrm{LGRG}}(\bbeta_n)
=
\frac{\E_{\mathbf{0},\bI_k}\left[\rho'(\|\bz\|)^2\right]}{k\alpha^2}
(\bX^T\bSigma^{-1}\bX)^{-1},
\end{equation}
found in Corollary~\ref{cor:Asymp norm elliptical},
and similarly for the limiting covariance matrix of $\sqrt{n}(\btheta_n-\btheta)$.
However, this is just the consequence of the extreme simple choices for $\bX$, $\bZ$ and $\bR$.
Already, if we keep~$\bX$ as it is, and only take a slight variation of either $\bZ$ or $\bR$, one finds
severe differences
between~\eqref{eq:varbeta copt} and~\eqref{eq:varbeta our}, and similarly for
the expression of the limiting covariance matrix of $\sqrt{n}(\btheta_n-\btheta)$.

We considered the following two alternatives
\begin{enumerate}
\item
take $\bZ=(1,2,3,4)^T$ and leave $\bX$ and $\bR=\bI_4$ as they are;
\item
take $\bR=(1,4,9,16)^T$ and leave $\bX$ and $\bZ=(1,1,1,1)^T$ as they are.
\end{enumerate}
We generated 10\,000 samples of size $n=100$ according to model~\eqref{eq:model copt}
and computed the value of S-estimators $\bbeta_n$ and $\btheta_n$
by means of Tukey's bi-weight
\begin{equation}
\label{def:biweight}
\rho_{\mathrm{B}}(s;c)
=
\begin{cases}
s^2/2-s^4/(2c^2)+s^6/(6c^4), & |s|\leq c\\
c^2/6 & |s|>c.
\end{cases}
\end{equation}
and $b_0=\E_{\mathbf{0},\bI_k}[\rho_B(\|\bz\|;c_0)]$, with the cut-off value $c_0$ chosen such that $b_0/a_0=0.5$.
According to Theorem~\ref{th:BDP}, this corresponds to (asymptotic) breakdown point 50\%.

\begin{figure}
  \centering
  \includegraphics[width=0.95\textwidth,clip=]{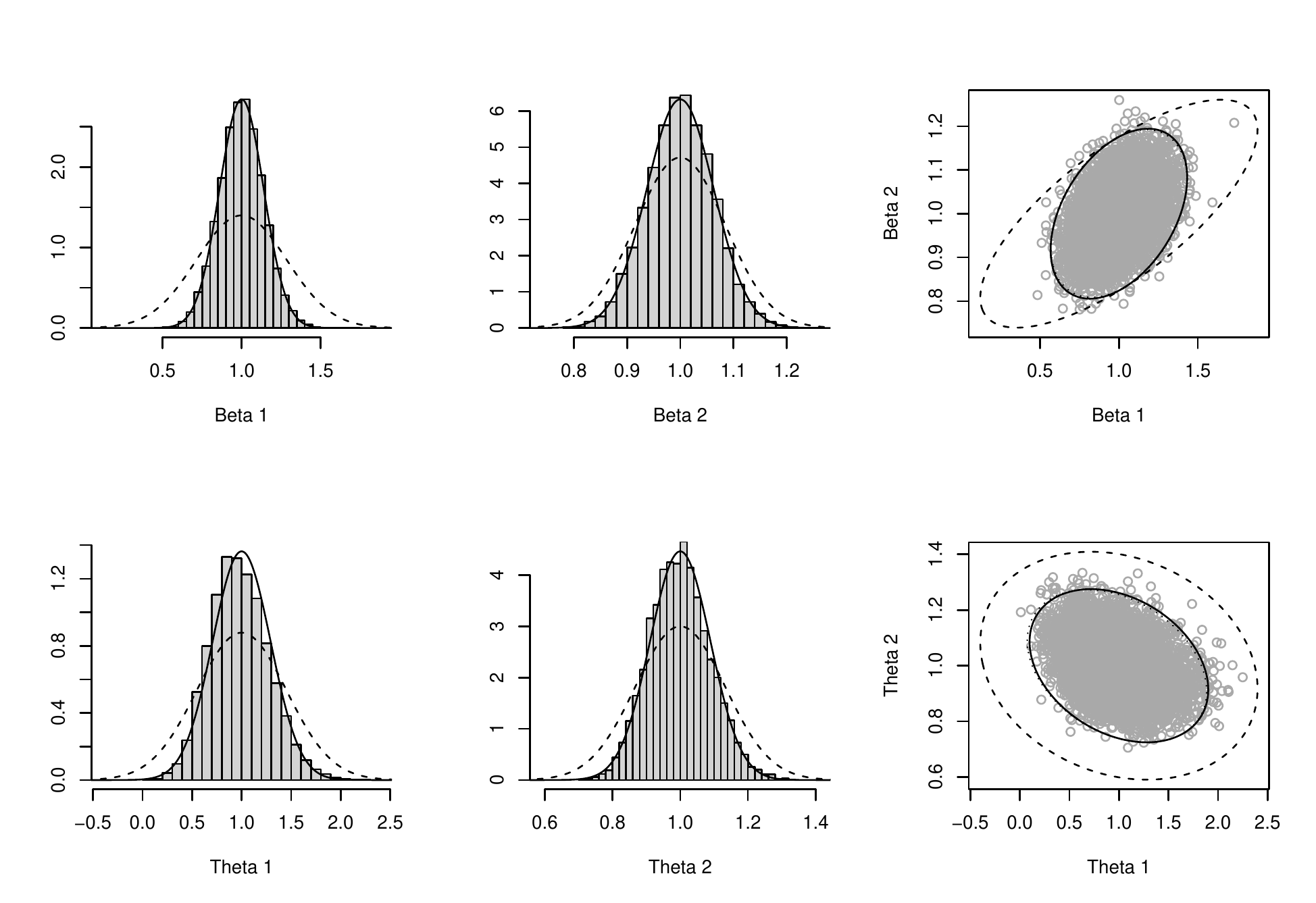}
  \caption{Empirical marginal and joined distributions together with limiting marginal and joined distributions of $\sqrt{n}(\bbeta_n-\bbeta)$ (first row) and
  $\sqrt{n}(\btheta_n-\btheta)$ (second row).}
  \label{fig:variances}
\end{figure}

Figure~\ref{fig:variances} displays the limiting marginal and joined distributions
of $\sqrt{n}(\bbeta_n-\bbeta)$ in the first row, where we generated the samples with alternative 1.
The histograms and scatterplot correspond to the 10\,000 different values of $\sqrt{n}(\bbeta_n-\bbeta)$.
The dashed curves correspond to the densities and 95\% contourlines of the theoretical limiting
marginal and joined normal distributions using the covariance matrix in~\eqref{eq:varbeta copt}.
The solid curves correspond to the marginal and joined normal distributions using the covariance matrix in~\eqref{eq:varbeta our}.
The empirical contourlines based on the sample mean and sample covariance of the 10\,000 estimates are plotted in dotted lines,
but they almost indistinguishable from the solid contourlines.
We find
\[
\text{Var}_{\mathrm{CVF}}(\bbeta_n)
=
\left(
  \begin{array}{cc}
8.13 & 1.78 \\
1.78 & 0.72
\end{array}
\right)
\qquad\text{and}\qquad
\text{Var}_{\mathrm{LGRG}}(\bbeta_n)
=
\left(
  \begin{array}{cc}
1.97 & 0.38\\
0.38 & 0.40
\end{array}
\right).
\]
Clearly, the histograms of the repeated estimates for $\beta_1$ and $\beta_2$ match the graphs of the (marginal)
normal densities with the variances given by $\text{Var}_{\mathrm{LGRG}}(\bbeta_n)$,
and the scatterplot matches with the 95\% contourline corresponding to $\text{Var}_{\mathrm{LGRG}}(\bbeta_n)$.
Note that the differences with $\text{Var}_{\mathrm{CVF}}(\bbeta_n)$ are quite severe.
For example, this yields that the length of the confidence interval for $\beta_1$
based on $\text{Var}_{\mathrm{CVF}}(\bbeta_n)$ will be two times larger than the one
based on $\text{Var}_{\mathrm{LGRG}}(\bbeta_n)$.

The second row in Figure~\ref{fig:variances} displays the limiting distributions
of $\sqrt{n}(\btheta_n-\btheta)$, where we generated the samples with alternative 2.
In~\cite{copt2006high}, the limiting covariance matrix was given by (see (15) in~\cite{copt2006high})
$(\bL^T\bL)^{-1}\bL^T\bV_{\bSigma}\bL(\bL^T\bL)^{-1}$,
where $\bV_{\bSigma}=\sigma_1(\bI_{k^2}+\bK_{k,k})(\bSigma\otimes\bSigma)
+
\sigma_2\vc(\bSigma)\vc(\bSigma)^T$, see Corollary~5.1 in~\cite{lopuhaa1989}.
Because
\[
\begin{split}
(\bL^T\bL)^{-1}\bL^T(\bI_{k^2}+\bK_{k,k})
&=
2(\bL^T\bL)^{-1}\bL^T,\\
(\bL^T\bL)^{-1}\bL^T\vc(\bSigma)
&=
\btheta(P),
\end{split}
\]
the expression given in~\cite{copt2006high} becomes
\[
\text{Var}_{\mathrm{CVF}}(\btheta_n)
=
2\sigma_1
(\bL^T\bL)^{-1}\bL^T
\left(\bSigma\otimes\bSigma\right)
\bL
(\bL^T\bL)^{-1}
+
\sigma_2
\btheta(P)\btheta(P)^T.
\]
This differs from our Corollary~\ref{cor:Asymp norm elliptical}, which gives
\[
\text{Var}_{\mathrm{LGRG}}(\btheta_n)
=
2\sigma_1\Big(\bL^T\left(\bSigma^{-1}\otimes\bSigma^{-1}\right)\bL\Big)^{-1}
+
\sigma_2
\btheta(P)\btheta(P)^T.
\]
For the choices of $\bX$, $\bZ$ and $\bR$ in~\cite{copt2006high}, both expressions are equal.
However, for the alternative choice for $\bR$ made in alternative 2, one finds
\[
\text{Var}_{\mathrm{CVF}}(\btheta_n)
=
\left(
  \begin{array}{rr}
20.63 & -1.22 \\
 -1.22 & 1.77
\end{array}
\right)
\qquad\text{and}\qquad
\text{Var}_{\mathrm{LGRG}}(\btheta_n)
=
\left(
  \begin{array}{rr}
8.57 & -0.82\\
-0.82  & 0.80
\end{array}
\right).
\]
Again the differences are quite large.
For example, as a consequence the length of the confidence interval for $\theta_1$
based on $\text{Var}_{\mathrm{CVF}}(\btheta_n)$ will be 1.5 times larger than the one
based on $\text{Var}_{\mathrm{LGRG}}(\btheta_n)$.

\begin{figure}[t!]
  \centering
  \includegraphics[width=\textwidth]{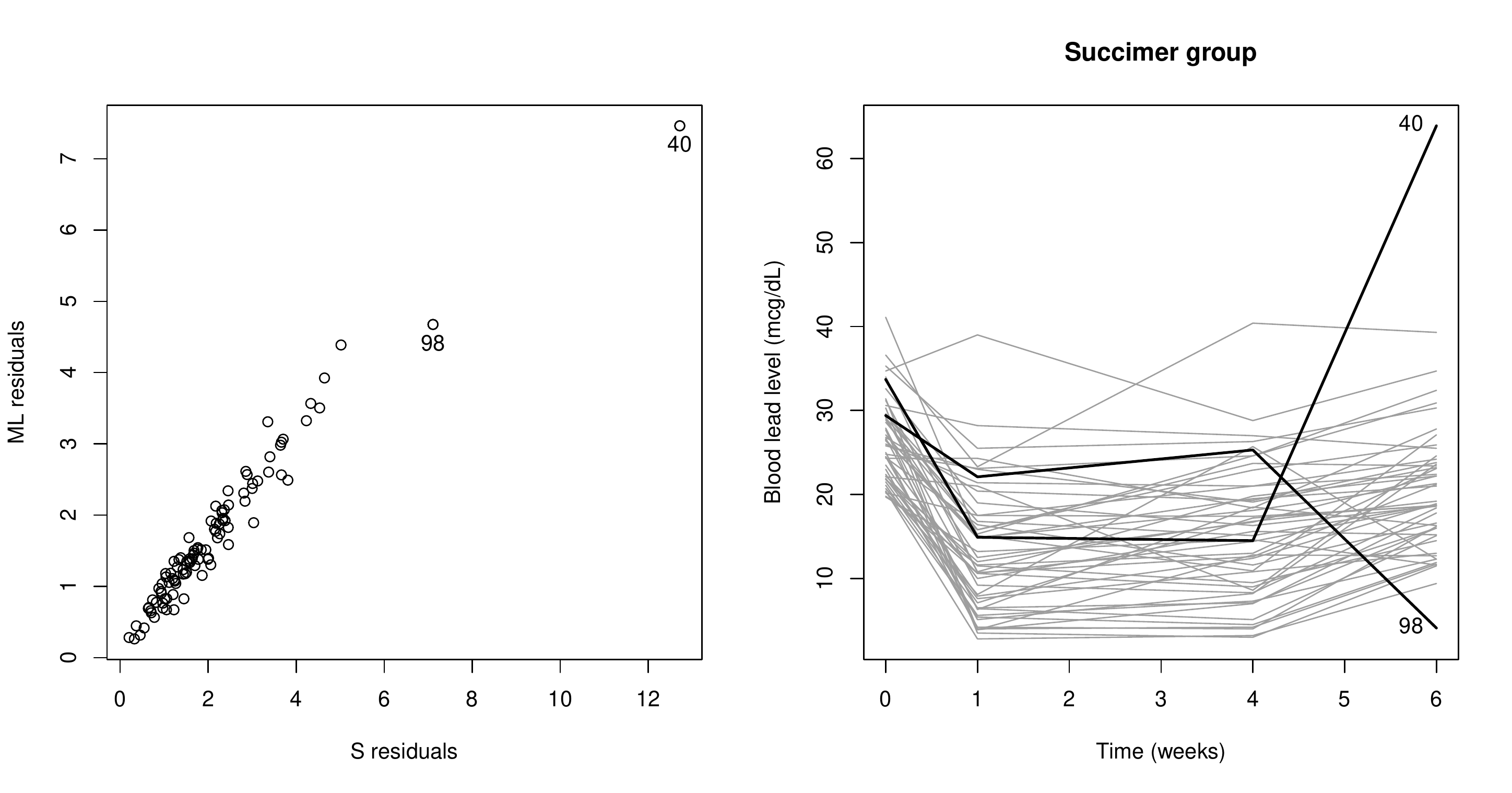}
  \caption{Left picture: standardized residuals for the S-estimates (horizontal axis)
  and the ML estimates (vertical axis).
  Right picture: observations for the subjects in the treatment group.}
  \label{fig:TLC}
\end{figure}
Finally, we illustrate the performance of S-estimators by an application to data from a trial on the treatment of lead-exposed children.
This dataset is discussed in~\cite{fitzmaurice-laird-ware2011} and consists of four repeated measurements of blood lead
levels obtained at baseline (or week~0), week 1, week 4, and week 6 on 100 children who were randomly assigned to chelation treatment
with succimer (a chelation agent) or placebo.
On the basis of a graphical display of the mean response over time, it is suggested in~\cite{fitzmaurice-laird-ware2011} that
a quadratic trend over time seems suitable.
We fitted the following model
\[
y_{ij}=
\beta_0+\beta_1\delta_{i}
+
(\beta_3+\beta_4\delta_i)t_j
+
(\beta_5+\beta_6\delta_i)t_j^2
+
\gamma_{1i}+\gamma_{2i}t_j+\gamma_{3i}t_j^2
+
\epsilon_{ij},
\]
for $i=1,\ldots,100$ and $j=1,\ldots,4$, where
$(t_1,\ldots,t_4)=(0,1,4,6)$ refer to the different weeks,
$y_{ij}$ is the blood lead level (mcg/dL) of subject $i$ obtained at time $t_j$,
and
$\delta_i=0$ if the $i$-th subject is in the placebo group and $\delta_i=1$, otherwise.
The random effects $\bgamma_i=(\gamma_{1i},\gamma_{2i},\gamma_{3i})$, $i=1,\ldots,100$, are assumed to be independent
mean zero normal random vectors with a diagonal covariance matrix consisting of variances $\sigma_{\gamma_1}^2$,
$\sigma_{\gamma_2}^2$ and $\sigma_{\gamma_3}^2$, respectively.
The measurement errors $\beps_i=(\epsilon_{i1},\ldots,\epsilon_{i4})$, $i=1,\ldots,100$, are assumed to be independent
mean zero random vectors with covariance matrix $\sigma_\epsilon^2\bI_4$,
also being independent of the random effects.
In this way we are fitting a balanced linear mixed effects model with unknown parameters $\bbeta=(\beta_1,\ldots,\beta_6)$
and $\btheta=(\sigma_{\gamma_1}^2,\sigma_{\gamma_2}^2,\sigma_{\gamma_3}^2,\sigma_{\epsilon}^2)$, and a linear covariance structure.

We estimated $(\bbeta,\btheta)$ by means of maximum likelihood and by means of the S-estimator corresponding to
Tukey's bi-weight defined in~\eqref{def:biweight}.
The tuning-constant was chosen to be $c=4.097$, which corresponds to asymptotic breakdown point 0.5.
For each estimate $(\widehat\bbeta,\widehat\btheta)$, we determined the estimate $\bV(\widehat\btheta)$ for the structured covariance and the
standardized residuals for each subject
\[
\mathrm{RES}_i=\sqrt{(\by_i-\bX_i\widehat{\bbeta})^T\bV(\widehat\btheta)^{-1}(\by_i-\bX_i\widehat{\bbeta})}
\]
The residuals for both estimation procedures are visible in the left picture of Figure~\ref{fig:TLC},
with the residuals determined from the S-estimate on the horizonal axis and
the ones determined from the ML estimate on the vertical axis.
Both estimates identify subject 40 as an outlier, but only the robust S-estimate also clearly identifies observation 98 as outlier.
The extreme large observation in week 6 seems to be the reason that observation 40 is identified as outlier by both methods.
See the right picture in Figure~\ref{fig:TLC}.
Observation 98 also seems to deviate from the overall quadratic trend, by having a suspicious low observation in week 6.
The corresponding S-residual clearly sticks out from the other S-residuals,
whereas this is much less so for the corresponding ML residual.

\appendix

\numberwithin{lemma}{section}

\section{Proofs and technical lemmas}
\label{sec:appendix}

\paragraph*{Proof of Theorem~\ref{th:existence structured}}
\begin{proof}
Let $(\bbeta,\btheta)\in\R^q\times\R^l$ satisfy the S-constraint in~\eqref{def:Smin structured}.
Then from (R1)-(R2) it follows that
\begin{equation}
\label{eq:enough mass structured}
P(\mathcal{C}(\bbeta,\bV(\btheta),c_0))
\geq
1-\frac{1}{a_0}\int\rho\left(\sqrt{(\by-\bX\bbeta)^T\bV(\btheta)^{-1}
(\by-\bX\bbeta)}\right)\,\dd P(\mathbf{s})
\geq
1-r.
\end{equation}
Since $1-r\geq\epsilon$,
Lemma~\ref{lem:compact structured}(i) then implies that $\lambda_k(\bV(\btheta))\geq a_1>0$.
Because
\[
\lim_{m\to\infty}\int\rho\left(\|\by\|/m\right)\,\dd P(\by,\bX)=0,
\]
we can find $m_0>0$, such that $\int\rho(\|\mathbf{\by}\|/m_0)\,\dd P(\by,\bX)\leq b_0$.
Lemma~\ref{lem:compact structured}(ii) then yields that $\lambda_1(\bV(\btheta))\leq a_2<\infty$.
Application of Lemma~\ref{lem:compact structured}(iii), with $a=1-r$ together with $(\text{C1}_\epsilon)$,
implies that $\|\bbeta\|\leq M<\infty$.
It follows that $\bbeta$ is in a compact subset of $\R^q$ and $\bV(\btheta)$ is in a compact set $K\subset \R^{k\times k}$.

According to~\eqref{def:identifiable}, the mapping $\btheta\mapsto \bV(\btheta)$ is one-to-one,
so that we can restrict~$\btheta$ to the pre-image~$\bV^{-1}(K)$.
Then with conditions (V1) and (V3) it follows that also $\bV^{-1}(K)$ is compact in~$\R^l$.
We conclude that for solving minimization problem~\eqref{def:Smin structured},
we can restrict ourselves to a compact set $K'\subset\R^q\times\R^l$.
As $\text{det}(\bV(\btheta))$ is a continuous function of
$(\bbeta,\btheta)$, due to condition~(V1), it must attain a minimum on~$K'$.
\end{proof}

\paragraph*{Proof of Corollary~\ref{cor:existence estimator structured}}
\begin{proof}
Let $\mathbb{P}_n$ be the empirical measure corresponding to the collection $\mathcal{S}_n$.
Then $\mathbb{P}_n$ satisfies  $(\text{C1}_\epsilon)$ for any $0<\epsilon\leq 1-r$ and
satisfies $(\text{C2}_\epsilon)$, for $\epsilon=(\kappa(\mathcal{S}_n)+1)/n$.
Clearly $0<\epsilon\leq 1-r$, where $r=b_0/a_0$, so according to Theorem~\ref{th:existence structured} there exists at least one solution
to~\eqref{def:Smin structured} with $P=\mathbb{P}_n$.
This means that there exists at least one solution to~\eqref{def:Smin estimator structured}.
\end{proof}

\paragraph*{Proof of Corollary~\ref{cor:existence weak convergence}}
\begin{proof}
First note there exists $0<\eta<\epsilon'-\epsilon$.
According to Ranga Rao~\cite[Theorem 4.2]{rangarao1962} we have
\begin{equation}
\label{eq:ranga rao structured}
\sup_{C\in \mathfrak{C}}
\left|
P_t(C)-P(C)
\right|
\to
0,
\quad
\text{as }t\to\infty.
\end{equation}
Because strips $H(\balpha,\ell,\delta)\in \mathfrak{C}$,
property~\eqref{eq:ranga rao structured} implies that every strip with $P_t(H(\balpha,\ell,\delta))\geq \epsilon+\eta$,
for $t$ sufficiently large, must also satisfy
$P(H(\balpha,\ell,\delta))\geq \epsilon$.
This means that
\[
\inf
\left\{
\delta:P_t(H(\balpha,\ell,\delta))\geq \epsilon+\eta
\right\}
\geq
\inf
\left\{
\delta:P(H(\balpha,\ell,\delta))\geq \epsilon
\right\}>0.
\]
It follows that, for $t$ sufficiently large, $P_t$ satisfies condition $(\text{C2}_{\epsilon+\eta})$.
Next, consider the compact set $K$ from ($\text{C1}_{\epsilon'}$).
Without loss of generality we may assume that it belongs to $\mathfrak{C}$.
Therefore,
as $P(K)\geq r+\epsilon'$, for $t$ sufficiently large $P_t(K)\geq r+\epsilon+\eta$.
It follows that, for~$t$ sufficiently large, $P_t$ satisfies condition $(\text{C1}_{\epsilon+\eta})$.
Since $\epsilon+\eta<1-r$, according to Theorem~\ref{th:existence structured}
at least one solution~$\bxi(P_t)$ exists,
for $t$ sufficiently large.
\end{proof}

\paragraph*{Proof of Theorem~\ref{th:continuity structure}}
\begin{proof}
First note that there exists $0<\eta<\epsilon'-\epsilon$.
Denote $\bxi(P_t)=\bxi_t=(\bbeta_t,\btheta_t)$.
Similar to~\eqref{eq:enough mass structured} we find that $P_t(\mathcal{C}(\bbeta_t,\btheta_t,c_0))\geq 1-r$.
Therefore, as $\mathcal{C}(\bbeta_t,\btheta_t,c_0)\in \mathfrak{C}$
and $1-r>1-r-\eta$,
it follows from~\eqref{eq:ranga rao structured} that
\begin{equation}
\label{eq:bound P(Ct)}
P(\mathcal{C}(\bbeta_t,\btheta_t,c_0))
\geq
P_t(\mathcal{C}(\bbeta_t,\btheta_t,c_0))
-
\sup_{C\in \mathfrak{C}}
\left|
P_t(C)-P(C)
\right|
\geq
1-r-\eta,
\end{equation}
for $t$ sufficiently large.
Since $1-r-\eta>\epsilon$, this means that, according to Lemma~\ref{lem:compact structured}(i),
there exists $a_1>0$ only depending on~$c_0$ and $P$, such that
for $t$ sufficiently large,
\[
\lambda_k(\bV(\btheta_t))\geq a_1>0.
\]
Denote $\bxi(P)=\bxi_0=(\bbeta_0,\btheta_0)$ and let $d_0(\bs)=d(\bs,\bbeta_0,\btheta_0)$.
Then according to Lemma~\ref{lem:billingsley structured} in~\cite{supplement},
for any $\sigma>-1$,
\[
\int
\rho
\left(
\frac{d_0(\bs)}{1+\sigma}
\right)
\,\dd P_t(\bs)
\to
\int
\rho
\left(
\frac{d_0(\bs)}{1+\sigma}
\right)
\,\dd P(\bs),
\]
as $t\to\infty$.
As the limit is strictly decreasing at $\sigma=0$,
and $\bxi_0$ satisfies the constraint in~\eqref{def:Smin structured},
we find that for all $\sigma>0$,
\[
\int
\rho
\left(
\frac{d_0(\bs)}{1+\sigma}
\right)
\,\dd P_t(\bs)
\leq
b_0,
\]
for $t$ sufficiently large.
Hence, similar to the proof of Lemma~\ref{lem:compact structured}(ii) we conclude that
for any possible solution $\bxi_t=(\bbeta_t,\btheta_t)$,
it must hold that
\[
\text{det}(\bV(\btheta_t))
\leq (1+\sigma)^{2k}
\text{det}(\bV(\btheta_0)),
\]
for $t$ sufficiently large.
As $\sigma>0$ can be taken arbitrarily small, we conclude that
\begin{equation}
\label{eq:upper bnd det}
\limsup_{t\to\infty}
\text{det}
(\bV(\btheta_t))
\leq
\text{det}(\bV(\btheta_0)),
\end{equation}
and we find that $\lambda_1(\bV(\btheta_t))\leq \text{det}(\bV(\btheta_0))/a_1^{k-1}<\infty$,
for $t$ large sufficiently large.
Finally, let $K$ be the compact set from $(\text{C1}_{\epsilon'})$, so that
$P(K)\geq r+\epsilon'> r+\epsilon+\eta$.
Then, according to~\eqref{eq:bound P(Ct)}, it
follows from Lemma~\ref{lem:compact structured}(iii) with $a=1-r-\eta$, that there exists
$0<M<\infty$ such that $\|\bbeta_t\|\leq M$, for $t$ sufficiently large.
This means that there exists a compact set $K'$,
such that for $t$ sufficiently large the sequence $\{(\bbeta_t,\bV(\btheta_t))\}\subset K'$.
Then, similar to the second part of the proof of Theorem~\ref{th:existence structured},
the conditions on the mapping $\btheta\mapsto \bV(\btheta)$ yield that there exists a compact set $K''\subset\R^{q+l}$,
such that for $t$ sufficiently large, the sequence $\{\bxi_t\}\subset K''$.

Consider a convergent subsequence $\{\bxi_{t_j}\}$ with $\bxi_{t_j}\to\bxi_L$.
With Lemma~\ref{lem:billingsley structured} in~\cite{supplement} and the fact that~$\bxi_{t_j}$ satisfies the
S-constraint in~\eqref{def:Smin structured} at $P=P_{t_j}$,
we find
\[
\int
\rho(d(\bs,\bxi_L))
\,\dd P(s)
=
\lim_{j\to\infty}
\int
\rho(d(s,\bxi_{t_j}))
\,\dd P_{t_j}(s)
\leq
b_0.
\]
Hence, $\bxi_L$ satisfies the S-constraint in~\eqref{def:Smin structured},
which has solution $\bxi_0$.
This means that  $\text{det}(\bV(\btheta_L))\geq \text{det}(\bV(\btheta_0))$.
But then from~\eqref{eq:upper bnd det}, it follows $\text{det}(\bV(\btheta_L))=\text{det}(\bV(\btheta_0))$.
Uniqueness of $\bxi_0$ together with identifiability~\eqref{def:identifiable} then implies that $\bxi_L=\bxi_0$.
Because~$\{\bxi_t\}$ eventually stays in a compact set, this means that we must have
$\lim_{t\to\infty}\bxi_t=\bxi_0$.
\end{proof}

\paragraph*{Proof of Corollary~\ref{cor:consistency S-estimator}}
\begin{proof}
We apply Theorem~\ref{th:continuity structure} to the sequence $\mathbb{P}_n$, $n=1,2,\ldots$,
of probability measures, where~$\mathbb{P}_n$ is the empirical measure corresponding to
$(\mathbf{y}_1,\mathbf{X}_1),\ldots,(\mathbf{y}_n,\mathbf{X}_n)$.
According to the Portmanteau Theorem (e.g., see Theorem~2.1 in~\cite{billingsley1968}),
$\mathbb{P}_n$ converges weakly to $P$, with probability one.
The corollary then follows from Theorem~\ref{th:continuity structure}.
\end{proof}

\paragraph*{Proof of Theorem~\ref{th:IF}}
\begin{proof}
Denote $\bxi_{h,\bs_0}=\bxi(P_{h,\bs_0})$.
This solution satisfies the score equation~\eqref{eq:Psi=0} for the regression S-functional at $P_{h,\bs_0}$, that is
\[
\int \Psi(\bs,\bxi_{h,\bs_0})\,\dd P_{h,\bs_0}(\bs)=\mathbf{0}.
\]
We decompose as follows
\[
\begin{split}
\mathbf{0}
&=
\int \Psi(\bs,\bxi_{h,\bs_0})\,\dd P_{h,\bs_0}(\bs)\\
&=
(1-h)\int \Psi(\bs,\bxi_{h,\bs_0})\,\dd P(\bs)+h\Psi(\bs_0,\bxi_{h,\bs_0})\\
&=
(1-h)
\Lambda(\bxi_{h,\bs_0})
+
h\Big(\Psi(\bs_0,\bxi_{h,\bs_0})-\Psi(\bs_0,\bxi(P))\Big)
+
h\Psi(\bs_0,\bxi(P)).
\end{split}
\]
We first determine the order of $\bxi_{h,\bs_0}-\bxi(P)$, as $h\downarrow0$.
Because $\bxi\mapsto\Psi(\bs_0,\bxi)$ is continuous, it follows that
\[
\Psi(\bs_0,\bxi_{h,\bs_0})=\Psi(\bs_0,\bxi(P))+o(1),
\qquad
\text{as }
h\downarrow0.
\]
Furthermore, because $\bxi\mapsto\Lambda(\bxi)$ is continuously differentiable at $\bxi(P)$, we have that
\[
\begin{split}
\Lambda(\bxi_{h,\bs_0})
&=
\Lambda(\bxi(P))
+
\bD(P)
(\bxi_{h,\bs_0}-\bxi(P))+
o(\|\bxi_{h,\bs_0}-\bxi(P)\|).
\end{split}
\]
Since $\bxi(P)$ is the S-functional at $P$, it is a zero of the corresponding score equation, i.e.,
$\Lambda(\bxi(P))=0$.
It follows that
\[
\mathbf{0}=
(1-h)
\bD(P)
(\bxi_{h,\bs_0}-\bxi(P))+
o(\|\bxi_{h,\bs_0}-\bxi(P)\|)
+
o(h)
+
h\Psi(\bs_0,\bxi(P)).
\]
Because $\bD(P)$ is non-singular
and $\Psi(\bs_0,\bxi(P))$ is fixed, this implies
$\bxi_{h,\bs_0}-\bxi(P)=O(h)$.
After inserting this in the previous equality, it follows that
\[
\begin{split}
\mathbf{0}
&=
(1-h)
\bD(P)
(\bxi_{h,\bs_0}-\bxi(P))+
h\Psi(\bs_0,\bxi(P))
+o(h)\\
&=
\bD(P)
(\bxi_{h,\bs_0}-\bxi(P))+
h\Psi(\bs_0,\bxi(P))
+o(h).
\end{split}
\]
We conclude
\[
\frac{\bxi_{h,\bs_0}-\bxi(P)}{h}
=
-\bD(P)^{-1}\Psi(\bs_0,\bxi(P))+o(1),
\qquad
\text{as }
h\downarrow0.
\]
This means that the limit of the left hand side exists and
\[
\text{IF}(\bs_0;\bxi,P)
=
\lim_{h\downarrow0}
\frac{\bxi((1-h)P+h\delta_{\bs_0})-\bxi(P)}{h}
=
-\bD(P)^{-1}\Psi(\bs_0,\bxi(P)).
\]
Next, consider the covariance functional $\bC(P)=\bV(\btheta(P))$.
By definition
\[
\mathrm{IF}(\bs_0;\vc(\bC),P)
=
\lim_{h\downarrow0}
\frac{\vc(\bC(P_{h,\bs_0}))-\vc(\bC(P))}{h}.
\]
Due to~(V4), by applying the chain rule
(e.g., see~\cite[Theorem~12, page 108]{magnus&neudecker1988}), we find
\[
\begin{split}
\vc(\bC(P_{h,\bs_0}))-\vc(\bC(P))
&=
\frac{\partial \vc(\bV(\btheta(P)))}{\partial\btheta^T}
\Big(
\btheta(P_{h,\bs_0})-\btheta(P)
\Big)
+
o(\|\btheta(P_{h,\bs_0})-\btheta(P)\|)\\
&=
\left(\frac{\partial \vc(\bV(\btheta(P)))}{\partial\btheta^T}\right)
\Big(
\btheta(P_{h,\bs_0})-\btheta(P)
\Big)
+
o(h).
\end{split}
\]
Dividing by $h$ and letting $h\downarrow0$, finishes the proof.
\end{proof}

\paragraph*{Proof of Theorem~\ref{th:asymp normal}}
\begin{proof}
Write $\bxi_P=\bxi(P)$.
Then from~\eqref{eq:decomposition estimator} and Lemma~\ref{lem:stoch equi} in~\cite{supplement}, it follows that
\begin{equation}
\label{eq:decomp Lambda}
\mathbf{0}=
\Lambda(\bxi_n)
+
\int \Psi(\mathbf{s},\bxi_P)\,\dd (\mathbb{P}_n-P)(\mathbf{s})
+
o_P(n^{-1/2}).
\end{equation}
From Theorem~\ref{th:continuity structure}, we know that $\bxi_n\to\bxi_P$ with probability one.
For the first term on the right hand side we have that
\[
\Lambda(\bxi_n)
=
\bD(P)(\bxi_n-\bxi_P)+o_P(\|\bxi_n-\bxi_P\|)
=
\bD(P)(\bxi_n-\bxi_P)+o_P(\|\bxi_n-\bxi_P\|),
\]
using that $\Lambda(\bxi_P)=\mathbf{0}$, due to the fact that $\bxi_P$ is the solution of~\eqref{def:Smin structured},
see also~\eqref{eq:Psi=0}.
From Lemma~\ref{lem:Psi bounded} in~\cite{supplement} we have that $\|\Psi_{\bbeta}\|\leq C_1\|\bX\|$
and $\|\Psi_{\btheta}\|\leq C_2$, for universal constants $0<C_1,C_2<\infty$.
Hence, from the conditions of the theorem, it follows that
$\E\|\Psi(\mathbf{s},\bxi(P))\|^2<\infty$.
This means that for the second term on the right hand side of~\eqref{eq:decomp Lambda},
\[
\int
\Psi(\mathbf{s},\bxi_P)
\,\dd (\mathbb{P}_n-P)(\mathbf{s})
=
\frac1n\sum_{i=1}^n
\left(
\Psi(\mathbf{s}_i,\bxi_P)-\E\Psi(\mathbf{s},\bxi_P)
\right)
=
O_P(n^{-1/2}),
\]
according to the central limit theorem.
It follows that
\[
\mathbf{0}
=
\bD(P)(\bxi_n-\bxi_P)+o_P(\|\bxi_n-\bxi_P\|)
+
O_P(n^{-1/2}),
\]
so that $\|\bxi_n-\bxi_P\|=O_P(n^{-1/2})$.
If we insert this is in~\eqref{eq:decomp Lambda}, we obtain
\[
\mathbf{0}
=
\bD(P)(\bxi_n-\bxi_P)
+
\frac1n\sum_{i=1}^n
\left(
\Psi(\mathbf{s}_i,\bxi_P)-\E\Psi(\mathbf{s},\bxi_P)
\right)
+
o_P(n^{-1/2}),
\]
from which it follows that
\[
\sqrt{n}(\bxi_n-\bxi_P)
=
-\bD(P)^{-1}
\sqrt{n}
\left(
\frac1n\sum_{i=1}^n
\left(
\Psi(\mathbf{s}_i,\bxi_P)-\E\Psi(\mathbf{s},\bxi_P)
\right)
\right)
+
o_P(1).
\]
After application of the central limit theorem, we conclude that $\sqrt{n}(\bxi_n-\bxi_P)$ is asymptotically normal with mean zero and covariance matrix
\[
\bD(P)^{-1}
\E
\left[
\Psi(\mathbf{s},\bxi_P)\Psi(\mathbf{s},\bxi_P)^T
\right]
\bD(P)^{-1},
\]
where we use that $\E\Psi(\mathbf{s},\bxi_P)=\Lambda(\bxi_P)=\mathbf{0}$.
This proves the theorem.
\end{proof}

\bibliographystyle{abbrv}
\bibliography{SbalancedLME}

\begin{thebibliography}{10}

\bibitem{agostinelli2016composite}
C.~Agostinelli and V.~J. Yohai.
\newblock Composite robust estimators for linear mixed models.
\newblock {\em Journal of the American Statistical Association},
  111(516):1764--1774, 2016.

\bibitem{billingsley1968}
P.~Billingsley.
\newblock {\em Convergence of probability measures}.
\newblock John Wiley \& Sons, Inc., New York-London-Sydney, 1968.

\bibitem{chervoneva2014}
I.~Chervoneva and M.~Vishnyakov.
\newblock Generalized s-estimators for linear mixed effects models.
\newblock {\em Statistica Sinica}, 24(3):1257--1276, 2014.

\bibitem{copt&heritier2007}
S.~Copt and S.~Heritier.
\newblock Robust alternatives to the f-test in mixed linear models based on
  mm-estimates.
\newblock {\em Biometrics}, 63(4):1045--1052, 2007.

\bibitem{copt2006high}
S.~Copt and M.-P. Victoria-Feser.
\newblock High-breakdown inference for mixed linear models.
\newblock {\em Journal of the American Statistical Association},
  101(473):292--300, 2006.

\bibitem{davies1987}
P.~L. Davies.
\newblock Asymptotic behaviour of {$S$}-estimates of multivariate location
  parameters and dispersion matrices.
\newblock {\em Ann. Statist.}, 15(3):1269--1292, 1987.

\bibitem{demidenko2013}
E.~Demidenko.
\newblock {\em Mixed models}.
\newblock Wiley Series in Probability and Statistics. John Wiley \& Sons, Inc.,
  Hoboken, NJ, second edition, 2013.
\newblock Theory and applications with R.

\bibitem{donoho&huber1983}
D.~Donoho and P.~J. Huber.
\newblock The notion of breakdown point.
\newblock In {\em A {F}estschrift for {E}rich {L}. {L}ehmann}, Wadsworth
  Statist./Probab. Ser., pages 157--184. Wadsworth, Belmont, CA, 1983.

\bibitem{fasano-maronna-sued-yohai2012}
M.~V. Fasano, R.~A. Maronna, M.~Sued, and V.~J. Yohai.
\newblock Continuity and differentiability of regression {M} functionals.
\newblock {\em Bernoulli}, 18(4):1284--1309, 2012.

\bibitem{fitzmaurice-laird-ware2011}
G.~M. Fitzmaurice, N.~M. Laird, and J.~H. Ware.
\newblock {\em Applied longitudinal analysis}.
\newblock Wiley Series in Probability and Statistics. John Wiley \& Sons, Inc.,
  Hoboken, NJ, second edition, 2011.

\bibitem{hampel1974}
F.~R. Hampel.
\newblock The influence curve and its role in robust estimation.
\newblock {\em J. Amer. Statist. Assoc.}, 69:383--393, 1974.

\bibitem{hartley&rao1967}
H.~O. Hartley and J.~N.~K. Rao.
\newblock Maximum-likelihood estimation for the mixed analysis of variance
  model.
\newblock {\em Biometrika}, 54:93--108, 1967.

\bibitem{heritier-cantoni-copt-victoriafeser2009}
S.~Heritier, E.~Cantoni, S.~Copt, and M.-P. Victoria-Feser.
\newblock {\em Robust methods in biostatistics}.
\newblock Wiley Series in Probability and Statistics. John Wiley \& Sons, Ltd.,
  Chichester, 2009.

\bibitem{jennrich&schluchter1986}
R.~I. Jennrich and M.~D. Schluchter.
\newblock Unbalanced repeated-measures models with structured covariance
  matrices.
\newblock {\em Biometrics}, 42(4):805--820, 1986.

\bibitem{laird&ware1982}
N.~M. Laird and J.~H. Ware.
\newblock Random-effects models for longitudinal data.
\newblock {\em Biometrics}, 38(4):963--974, 1982.

\bibitem{supplement}
H.~Lopuha\"a, V.~Gares, and A.~Ruiz-Gazen.
\newblock Supplement to ``{S}-estimation in linear models with structured
  covariance matrices''.
\newblock 2022.

\bibitem{lopuhaa1989}
H.~P. Lopuha\"{a}.
\newblock On the relation between {$S$}-estimators and {$M$}-estimators of
  multivariate location and covariance.
\newblock {\em Ann. Statist.}, 17(4):1662--1683, 1989.

\bibitem{lopuhaa1997}
H.~P. Lopuha\"{a}.
\newblock Asymptotic expansion of {$S$}-estimators of location and covariance.
\newblock {\em Statist. Neerlandica}, 51(2):220--237, 1997.

\bibitem{lopuhaa&rousseeuw1991}
H.~P. Lopuha\"{a} and P.~J. Rousseeuw.
\newblock Breakdown points of affine equivariant estimators of multivariate
  location and covariance matrices.
\newblock {\em Ann. Statist.}, 19(1):229--248, 1991.

\bibitem{magnus&neudecker1988}
J.~R. Magnus and H.~Neudecker.
\newblock {\em Matrix differential calculus with applications in statistics and
  econometrics}.
\newblock Wiley Series in Probability and Mathematical Statistics: Applied
  Probability and Statistics. John Wiley \& Sons, Ltd., Chichester, 1988.

\bibitem{nolan&pollard1987}
D.~Nolan and D.~Pollard.
\newblock {$U$}-processes: rates of convergence.
\newblock {\em Ann. Statist.}, 15(2):780--799, 1987.

\bibitem{pinheiro-liu-wu2001}
J.~C. Pinheiro, C.~Liu, and Y.~N. Wu.
\newblock Efficient algorithms for robust estimation in linear mixed-effects
  models using the multivariate {$t$} distribution.
\newblock {\em J. Comput. Graph. Statist.}, 10(2):249--276, 2001.

\bibitem{pollard1984}
D.~Pollard.
\newblock {\em Convergence of stochastic processes}.
\newblock Springer Series in Statistics. Springer-Verlag, New York, 1984.

\bibitem{rangarao1962}
R.~Ranga~Rao.
\newblock Relations between weak and uniform convergence of measures with
  applications.
\newblock {\em Ann. Math. Statist.}, 33:659--680, 1962.

\bibitem{rao1972}
C.~R. Rao.
\newblock Estimation of variance and covariance components in linear models.
\newblock {\em J. Amer. Statist. Assoc.}, 67:112--115, 1972.

\bibitem{rousseeuw1985}
P.~Rousseeuw.
\newblock Multivariate estimation with high breakdown point.
\newblock In {\em Mathematical statistics and applications, {V}ol. {B} ({B}ad
  {T}atzmannsdorf, 1983)}, pages 283--297. Reidel, Dordrecht, 1985.

\bibitem{rousseeuw-yohai1984}
P.~Rousseeuw and V.~Yohai.
\newblock Robust regression by means of {S}-estimators.
\newblock In {\em Robust and nonlinear time series analysis ({H}eidelberg,
  1983)}, volume~26 of {\em Lect. Notes Stat.}, pages 256--272. Springer, New
  York, 1984.

\bibitem{rousseeuw1984}
P.~J. Rousseeuw.
\newblock Least median of squares regression.
\newblock {\em J. Amer. Statist. Assoc.}, 79(388):871--880, 1984.

\bibitem{vanaelst&willems2004}
S.~Van~Aelst and G.~Willems.
\newblock Multivariate regression {$S$}-estimators for robust estimation and
  inference, 2004, preprint received by personal communication.

\bibitem{vanaelst&willems2005}
S.~Van~Aelst and G.~Willems.
\newblock Multivariate regression {$S$}-estimators for robust estimation and
  inference.
\newblock {\em Statist. Sinica}, 15(4):981--1001, 2005.

\end{thebibliography}

\pagebreak

\newpage
\section{Supplemental Material}
\label{sec:supplement}
For later use we first define two important matrix norms and mention some useful properties.
For $m\times n$ real-valued matrices $\bA$, we define the Euclidean norm or Frobenius norm
as
\[
\|\bA\|=\sqrt{\sum_{i=1}^m\sum_{j=1}^n a_{ij}^2}.
\]
and the spectral norm by
\[
\|\bA\|_2=\sup_{\|\bu\|=1}\|\bA\bu\|.
\]
Recall that for real-valued $\bA$, the largest eigenvalue is defined by
\begin{equation}
\label{eq:largest eigenvalue}
\lambda_1(\bA)
=
\sup_{\|\bu\|=1}
\bu^T\bA\bu.
\end{equation}
This means that $\|\bA\|_2^2=\lambda_1(\bA^T\bA)$.
Other useful properties are
\begin{equation}
\label{eq:Frobenius - Spectral}
\|\bA\|_2\leq \|\bA\|\leq \sqrt{\min(m,n)}\|\bA\|_2
\end{equation}
and
\begin{equation}
\label{eq:euclidean - trace}
\|\bA\|^2=\text{tr}(\bA^T\bA)=\text{tr}(\bA\bA^T),
\end{equation}
and
\begin{equation}
\label{eq:euclidean uv}
\|\bu\bv^T\|
=
\|\bu\|
\|\bv\|
\end{equation}
for $\bu\in\R^m$ and $\bv\in\R^n$.
When $\bA$ is symmetric then $\lambda_1(\bA^T\bA)=\lambda_1(\bA^2)=\lambda_1(\bA)^2$.
In that case
\begin{equation}
\label{eq:euclidean - lambda1 symmetric}
\left|\lambda_1(\bA)\right|
=
\|\bA\|_2
\leq
\|\bA\|
\leq
\sqrt{\min(m,n)}
\|\bA\|_2
=
\sqrt{\min(m,n)}
\left|\lambda_1(\bA)\right|.
\end{equation}
Finally, note that both matrix norms are submultiplicative, that is
\begin{equation}
\label{eq:product triangle spectral}
\|\bA\bB\|_2\leq\|\bA\|_2\|\bB\|_2,
\end{equation}
and
\begin{equation}
\label{eq:product triangle}
\|\bA\bB\|\leq \|\bA\|\|\bB\|.
\end{equation}

\subsection{Proofs of Section~\ref{sec:existence}}
\paragraph*{Proof of Lemma~\ref{lem:compact structured}}
\begin{proof}
Note that cylinder $\mathcal{C}(\bbeta,\btheta,c)$ is contained
in some strip $H(\balpha,\ell,2c\sqrt{\lambda_k(\bV(\btheta))})$.  
It then follows from $(\text{C2}_\epsilon)$ that
$\lambda_k(\bV(\btheta))\geq \delta^2_\epsilon/4c^2$.
This proves (i).
Let $\btheta_0$ be such that the pair $(\bbeta,\bV(\btheta_0))=(\mathbf{0},m_0^2\bI_k)$
satisfies the S-constraint in~\eqref{def:Smin structured},
which is possible due to condition (V2).
It then follows that for any solution
$(\bbeta,\bV(\btheta))$ of~\eqref{def:Smin structured}, one must have
$\text{det}(\bV(\btheta))\leq m_0^{2k}$.
The lower bound on the smallest eigenvalue then implies
$\lambda_1(\bV(\btheta))\leq m_0^{2k}/a_1^{k-1}<\infty$,
which proves~(ii).
Note that
\begin{equation}
\label{eq:general position X structured}
\gamma_\epsilon
=
\inf_{P(J)\geq \epsilon}
\inf_{\|\bgamma\|=1}
\sup_{\bs\in J}
\|\bX\bgamma\|
>0,
\end{equation}
where the infima are taken over all subsets $J\subset \R^k\times \mathcal{X}$ with $P(J)\geq\epsilon$
and all unit vectors $\bgamma\in\R^q$.
This can be seen as follows.
Take
$\balpha^T
=
(\mathbf{0}^T,\bgamma^T,\ldots,\bgamma^T)/k$,
where $\mathbf{0}$ is a $k$-vector of zeros, so $\balpha\in\R^k\times \R^{kq}$ and $\|\balpha\|=1$.
Then, with $\ell=0$,
we have
$\balpha^T\bs-\ell
=
(\bgamma^T\mathbf{x}_1+\cdots+\bgamma^T\mathbf{x}_k)/k$.
Note that $\|\bX\bgamma\|^2=\sum_{j=1}^k (\mathbf{x}_j^T\bgamma)^2$.
Therefore, if $\gamma_\epsilon=0$, then $\bgamma^T\mathbf{x}_j=0$, for all $j=1,\ldots,k$,
which means that $\balpha^T\bs-\ell=0$.
This would be in contradiction with~\eqref{eq:general position},
which is equivalent to~$(\text{C2}_\epsilon)$, according to Remark~\ref{rem:condition C2}.
Now, note that
\[
P(\mathcal{C}(\bbeta,\btheta,c)\cap K)
\geq
P(\mathcal{C}(\bbeta,\btheta,c))-P(K^c)
\geq
a-1+1-a+\epsilon=\epsilon.
\]
Hence,
according to~\eqref{eq:general position X structured}, there exists
$\bs_0=(\by_0,\bX_0)\in \mathcal{C}(\bbeta,\btheta,c)\cap K$,
such that $\|\bX_0\bgamma\|\geq\gamma_\epsilon>0$, for all $\bgamma$ with $\|\bgamma\|=1$.
Because $\bs_0\in \mathcal{C}(\bbeta,\btheta,c)$, it holds
\[
\begin{split}
\|\by_0-\bX_0\bbeta\|^2
&\leq
(\by_0-\bX_0\bbeta)^T\bV(\btheta)^{-1}(\by_0-\bX_0\bbeta)
\lambda_1(\bV(\btheta))
\leq
c\lambda_1(\bV(\btheta))
\leq
ca_2,
\end{split}
\]
due to part~(ii).
Because $\bs_0=(\by_0,\bX_0)\in K$, this means that
\[
\|\bX_0\bbeta\|
\leq
\|\by_0-\bX_0\bbeta\|
+
\|\by_0\|
\leq
\sqrt{ca_2}
+
\sup_{(\by,\bX)\in K}\|\by\|.
\]
Because $\bs_0\in \mathcal{C}(\bbeta,\btheta,c)$,
together with~\eqref{eq:general position X structured} we conclude that
\[
\|\bbeta\|
=
\|\bX_0\bbeta\|
\times
\frac{\|\bbeta\|}{\|\bX_0\bbeta\|}
\leq
\frac{1}{\gamma_\epsilon}
\|\bX_0\bbeta\|
\leq
\frac{1}{\gamma_\epsilon}
\left(
\sqrt{ca_2}
+
\sup_{(\by,\bX)\in K}\|\by\|
\right)<\infty.
\]
This proves part (iii).
\end{proof}

\paragraph*{Proof of Remark~\ref{rem:condition C2}}
\begin{proof}
Suppose that $(\mathrm{C2}_\epsilon)$ holds and suppose that $\omega_\epsilon=0$.
Then there exists a sequence $(J_n,\balpha_n,\ell_n)$,
with $J_n\subset\R^k\times \mathcal{X}$, $\balpha_n\in\R^{k+kq}$, $\|\balpha_n\|=1$, and $\ell_n\in\R$,
such that $P(J_n)\geq\epsilon$, for all $n=1,2,\ldots$,
and
\[
\sup_{\bs\in J_n}|\balpha_n^T\bs-\ell_n|\to0.
\]
This means that for $n$ sufficiently large
$\delta_n=2\sup_{\bs\in J_n}|\balpha_n^T\bs-\ell_n|<\delta_\epsilon$.
Then, there exists a strip~$H(\balpha_n,\ell_n,\delta_n)$
containing~$J_n$, such that $\delta_n<\delta_\epsilon$
and
$P(H(\balpha_n,\ell_n,\delta_n)\cap (\R^k\times \mathcal{X}))\geq P(J_n)\geq \epsilon$.
This would be in contradiction with the definition of $\delta_\epsilon$ in~$(\mathrm{C2}_\epsilon)$.
On the other hand, suppose that~\eqref{eq:general position} holds and suppose that $\delta_\epsilon=0$.
That means that we can find a sequence $(\balpha_n,\ell_n,\delta_n)$,
with $\balpha_n\in\R^k\times\R^{kq}$, $\|\balpha_n\|=1$, and $\ell_n\in\R$, such that $\delta_n\downarrow0$
and $P(H(\balpha_n,\ell_n,\delta_n)\cap (\R^k\times \mathcal{X}))\geq \epsilon$, for all $n=1,2,\ldots$.
This means that for $n$ sufficiently large $\delta_n<2\omega_\epsilon$
and then
\[
\sup_{\bs\in H(\balpha_n,\ell_n,\delta_n)\cap \R^k\times \mathcal{X}}
|\balpha_n^T\bs-\ell_n|
=
\delta_n/2
<
\omega_\epsilon,
\]
which would contradict definition~\eqref{eq:general position}.
\end{proof}
\begin{lemma}
\label{lem:billingsley structured}
Let $P_t$, $t\geq0$ be a sequence of probability measures on $\R^k\times \R^{kq}$ that converges weakly to $P$, as $t\to\infty$.
Let $\bxi_t=(\bbeta_t,\btheta_t)$, $t\geq 0$, be a sequence in $\R^q\times\R^l$,
such that $\bxi_t\to\bxi_L$, as $t\to\infty$.
If $g(\bs,\bxi)=\rho(d(\bs,\bxi)/\alpha)$, for some $\alpha>0$ fixed, where
$\rho$ satisfies (R2)-(R3) and $\bV$ satisfies (V1),
then
\[
\lim_{t\to\infty}
\int
g(\bs,\bxi_t)\,\dd P_t(\bs)
=
\int
g(\bs,\bxi_L)\,\dd P(\bs).
\]
\end{lemma}
\begin{proof}
Let $g_t(\bs)=g(\bs,\bxi_t)$ and $g_L(\bs)=g(\bs,\bxi_L)$.
Then for every sequence $\{\bs_t\}$, such that $\bs_t\to \bs$, we have
\[
\lim_{t\to\infty}
g_t(\bs_t)=g_L(\bs).
\]
Now, apply Theorem~5.5 from~\cite{billingsley1968}.
Let $\Gamma:[0,\infty)\to[0,\infty)$ be the function
\[
\Gamma(u)=u1_{[0,a_0]}(u)+a_01_{(a_0,\infty)}(u),
\]
which is bounded and uniformly continuous.
Then as a consequence of $P_t\to P$ weakly, we have
\[
\lim_{t\to\infty}
\int
g(\bs,\bxi_t)\,\dd P_t(\bs)
=
\lim_{t\to\infty}
\int
\Gamma(g_t(\bs))\,\dd P_t(\bs)
=
\int
\Gamma(g_L(\bs))\,\dd P(\bs)
=
\int
g(\bs,\bxi_L)\,\dd P(\bs).
\]
\end{proof}

\subsection{Proofs of Section~\ref{sec:continuity}}
\paragraph*{Proof of Theorem~\ref{th:davies}}
\begin{proof}
Davies~\cite{davies1987} defines location-scale S-estimators by means of a function $\kappa:[0,\infty)\to[0,1]$.
It relates to our $\rho$-function as $\rho(d)=a_0(1-\kappa(d^2))$.
The S-minimization problem considered in~\cite{davies1987} can be formulated in our notation
as follows
\begin{equation}
\label{def:Smin davies}
\begin{split}
&
\min_{\balpha,\bA}
\text{det}(\bA)\\
&
\text{subject to}\\
&
\int \rho
\left(\sqrt{(\mathbf{y}-\balpha)^T\bA^{-1}(\mathbf{y}-\balpha)}\right)
f_{\bmu,\bSigma}(\by)
\,\dd \mathbf{y}
\leq
b_0,
\end{split}
\end{equation}
where the minimum is taken over all $\balpha\in\R^k$ and $\bA\in \mathrm{PDS}(k)$.
The conditions (R1)-(R2) imply the conditions on $\kappa(s)=1-\rho(\sqrt{s})/a_0$ imposed in~\cite{davies1987},
and $\kappa$ and $h$ have a common point of decrease.
It then follows from
Theorem~1 in~\cite{davies1987} that~\eqref{def:Smin davies} has a unique solution
\[
(\balpha^*,\bA^*)=(\bmu,\mathbf{\Sigma})=(\bX\bbeta_0,\bV(\btheta_0)).
\]
Since this solution is unique, candidate solutions to~\eqref{def:Smin davies}
must be of the form $(\balpha,\bA)=(\bX\bbeta,\bV(\btheta))$,
for some $\bbeta\in\R^q$ and $\btheta\in\R^l$.
It follows that minimization problem~\eqref{def:Smin davies} is
equivalent to minimization problem~\eqref{def:Smin regression}.
As a consequence, minimization problem~\eqref{def:Smin regression} has a unique solution~$(\bbeta^*,\btheta^*)$
for which $\bX\bbeta^*=\bX\bbeta_0$ and $\bV(\btheta^*)=\bV(\btheta_0)$.
Since $\bX^T\bX$ is non-singular we can multiply $\bX\bbeta^*=\bX\bbeta_0$ from the left by $(\bX^T\bX)^{-1}\bX^T$.
It then follows that $\bbeta^*=\bbeta_0$.
Finally, from~\eqref{def:identifiable} we find that~$\btheta^*=\btheta_0$.
This proves the theorem.
\end{proof}

\subsection{Proofs for Section~\ref{sec:bdp}}
\paragraph*{Proof of Theorem~\ref{th:BDP}}
\begin{proof}
Without loss of generality we may assume that $c_0=1$ and $\sup\rho=1$, so that $r=b_0$.
The first step is to show that for both estimators
$\epsilon^*_n\geq \lceil nr\rceil/n$.
To this end, consider a collection $\mathcal{S}_m'$ obtained from the original collection $\mathcal{S}_n$ by replacing
at most $m=\lceil nr\rceil-1$ number of points in $\R^k\times \mathcal{X}$.
We must show that at least one solution $(\bbeta_n(\mathcal{S}'_m),\btheta_n(\mathcal{S}'_m))$
to the S-minimization problem~\eqref{def:Smin estimator structured}
exists for the corrupted collection $\mathcal{S}'_m$, and that
all possible solutions do not break down.

Denote a possible solution to the S-minimization problem~\eqref{def:Smin estimator structured} corresponding to the corrupted collection $\mathcal{S}'_m$,
by
\[
\bxi_m'=
(\bbeta'_m,\btheta'_m)
=
(\bbeta_n(\mathcal{S}'_m),\btheta_n(\mathcal{S}'_m)),
\]
and consider the corresponding cylinder
$\mathcal{C}(\bbeta'_m,\btheta'_m,1)$ defined by~\eqref{def:cylinder}.
We apply Lemma~\ref{lem:compact structured} to the empirical measure $\mathbb{P}_m'$ corresponding to the
corrupted collection $\mathcal{S}_m'$ of~$n$ points.
Because $\bxi_m'=(\bbeta_m',\btheta_m')$
must satisfy the S-constraint in~\eqref{def:Smin estimator structured},
one can argue as in~\eqref{eq:enough mass structured},
\[
\begin{split}
\mathbb{P}_m'(\mathcal{C}(\bbeta_m',\btheta'_m,1))
&=
\frac1n
\sum_{\bs_i\in S_m'}
\mathds{1}
\left\{
\bs_i\in
\mathcal{C}(\bbeta_m',\btheta'_m,1)
\right\}\\
&\geq
1-\frac1n
\sum_{\bs_i\in S_m'}
\rho\left(d(\bs_i,\bxi_m')\right)
\geq
1-b_0,
\end{split}
\]
where $d(\bs_i,\bxi_m')$ is defined in~\eqref{def:Mahalanobis distance structured}.
It follows that the cylinder
$\mathcal{C}(\bbeta'_m,\btheta'_m,1)$
must contain at least $\lceil n-nb_0\rceil=\lceil n-nr\rceil$ number of points from the corrupted collection
$\mathcal{S}'_m$.
Furthermore, since $r\leq (n-\kappa(\mathcal{S}_n))/(2n)$,
for any such subset of $\mathcal{S}'_m$ it holds that
it contains
\begin{equation}
\label{eq:points in cylinder}
\lceil n-nr\rceil-m
=
n-\lfloor{nr}\rfloor-\lceil nr\rceil+1
\geq
\kappa(\mathcal{S}_n)+1
\end{equation}
points of the original collection $\mathcal{S}_n$.
It follows that the measure $\mathbb{P}_m'$ satisfies condition $(\text{C2}_\epsilon)$,
for~$\epsilon=(\kappa(\mathcal{S}_n)+1)/n$ and
with the value $\delta_\epsilon>0$ only depending on the original collection $\mathcal{S}_n$.
Moreover, we also have that $\mathbb{P}_m'(\mathcal{C}(\bbeta_m',\btheta'_m,1))\geq \epsilon$,
for $\epsilon=(\kappa(\mathcal{S}_n)+1)/n>0$.
According to Lemma~\ref{lem:compact structured}(i), it then follows that
$\lambda_k(\bV(\btheta'_m))\geq a_1>0$,
where $a_1$ only depends on the original collection $\mathcal{S}_n$.

Because $nb_0-m=nr-\lceil nr\rceil +1>0$ and
\[
\lim_{R\to\infty}
\sum_{(\by_i,\bX_i)\in \mathcal{S}_n}
\rho\left(
\frac{\|\by_i\|}{R}\right)
=
0,
\]
we can find an $R_0>0$, only depending on the original collection $\mathcal{S}_n$,
such that
\[
\sum_{(\by_i,\bX_i)\in \mathcal{S}_n}
\rho\left(
\frac{\|\by_i\|}{R_0}\right)
\leq
nb_0-m.
\]
The collection $\mathcal{S}'_m$ contains $n-m$ points of the original collection $\mathcal{S}_n$.
Consider the smallest~$M>0$, such that
\[
\sum_{(\by_i,\bX_i)\in \mathcal{S}_m'\cap \mathcal{S}_n}
\rho\left(
\frac{\|\by_i\|}{M}\right)
\leq
nb_0-m.
\]
Because $\mathcal{S}_m'\cap \mathcal{S}_n$ has less points than $\mathcal{S}_n$,
it holds that $M\leq R_0$.
It follows that
\[
\begin{split}
\int
\rho\left(\frac{\|\by\|}{R_0}\right)
\,\dd\mathbb{P}_m'(\bs)
&=
\frac1n
\sum_{(\by_i,\bX_i)\in \mathcal{S}_m'}
\rho\left(
\frac{\|\by_i\|}{R_0}\right)\\
&\leq
\frac1n
\sum_{(\by_i,\bX_i)\in S_m'}
\rho\left(
\frac{\|\by_i\|}{M}\right)\\
&\leq
\frac1n
\left(
\sum_{(\by_i,\bX_i)\in S_m'\cap S_n}
\rho\left(
\frac{\|\by_i\|}{M}\right)
+
m
\right)
\leq
b_0.
\end{split}
\]
According to Lemma~\ref{lem:compact structured}(ii), it then follows that
$\lambda_1(\bV(\btheta_m'))\leq a_2<\infty$.
where $a_2$ only depends on~$a_1$ and the collection $\mathcal{S}_n$.

To show that the estimate $\bbeta_m'$ stays bounded, recall that
the cylinder $\mathcal{C}(\bbeta_m',\btheta'_m,1)$ contains a subset $J_0$ of
$\kappa(\mathcal{S}_n)+1$ points from the original collection $\mathcal{S}_n$,
according to~\eqref{eq:points in cylinder}.
By definition, $\kappa(\mathcal{S}_n)+1$ original points cannot be on the same hyperplane, so that
\[
\gamma_n
=
\inf_{J\subset \mathcal{S}_n}
\inf_{\|\bgamma\|=1}\max_{\bs\in J}
\|\bX\bgamma\|>0.
\]
where the first infimum runs over all subsets $J\subset \mathcal{S}_n$ of $\kappa(\mathcal{S}_n)+1$ points.
By definition of~$\gamma_n$,
there exists an original point~$\bs_0\in J_0\subset \mathcal{S}_n\cap \mathcal{C}(\bbeta_m',\btheta'_m,1)$, such that
\[
\|\bbeta'_m\|
=
\|\bX_0\bbeta'_m\|
\times
\frac{\|\bbeta'_m\|}{\|\bX_0\bbeta'_m\|}
\leq
\frac{1}{\gamma_n}
\|\bX_0\bbeta'_m\|.
\]
Because $\bs_0\in \mathcal{C}(\bbeta_m',\btheta'_m,1)$,
similar to the proof of Lemma~\ref{lem:compact structured}(iii), it follows
that $\|\by_0-\bX_0\bbeta'_m\|^2\leq a_2$,
and because $\bs_0\in \mathcal{S}_n$, we have that
\[
\|\bX_0\bbeta'_m\|
\leq
\sqrt{a_2}
+
\max_{(\by_i,\bX_i)\in \mathcal{S}_n}
\|\by_i\|
<\infty.
\]
We conclude that there exists a compact set $K_n$, only depending on the original collection~$\mathcal{S}_n$,
that contains the pair $(\bbeta'_m,\bV(\btheta'_m))$.
Similar to the reasoning in the proof of Theorem~\ref{th:existence structured},
it follows that at least one solution $(\bbeta_n(\mathcal{S}'_m),\btheta_n(\mathcal{S}'_m))$
to the S-minimization problem~\eqref{def:Smin estimator structured}
exists for the collection $\mathcal{S}'_m$, and that
all possible solutions~$\bbeta_n(\mathcal{S}'_m)$ and $\btheta_n(\mathcal{S}'_m)$ do not break down.

We continue by showing
$\epsilon^*_n(\btheta_n)
\leq
\lceil nr\rceil/n$.
Replace $m=\lceil nr\rceil$ points of $\mathcal{S}_n$ to obtain a corrupted collection $\mathcal{S}'_m$
of $n$ points.
Suppose that a solution $\bxi_m'=(\bbeta_m',\btheta_m')=(\bbeta_n(\mathcal{S}'_m),\btheta_n(\mathcal{S}'_m))$
exists to the S-minimization problem~\eqref{def:Smin estimator structured}
corresponding to the corrupted collection $\mathcal{S}'_m$.
We must show that the estimate $\btheta_m'$ breaks down.
Note that the estimates~$\bbeta_m'$ and $\btheta_m'$ satisfy the S-constraint
in~\eqref{def:Smin estimator structured} for the corrupted collection,
\begin{equation}
\label{eq:constraint corrupted}
\sum_{\bs_i\in \mathcal{S}_m'}
\rho\left(d(\bs_i,\bxi_m')\right)
\leq
nr.
\end{equation}
If all $m=\lceil nr\rceil$ replaced points are outside the cylinder $\mathcal{C}(\bbeta_m',\btheta_m',1)$, then
\[
\sum_{\bs_i\in \mathcal{S}_m'}
\rho\left(d(\bs_i,\bxi_m')\right)
=
\sum_{\bs_i\in \mathcal{S}_m'\cap \mathcal{S}_n}
\rho\left(d(\bs_i,\bxi_m')\right)+
\lceil nr\rceil
>
nr,
\]
when $nr\notin \mathds{N}$.
When $nr\in \mathds{N}$, then by assumption $n-m=n-\lceil nr\rceil\geq \kappa(\mathcal{S}_n)+1$.
Hence, there is at least one point $\bs_i\in \mathcal{S}_n$,
for which $d(\bs_i,\bxi_m')>0$.
Because $\rho$ is strictly increasing on $[0,1]$, this implies
\[
\sum_{\bs_i\in \mathcal{S}_m'}
\rho\left(d(\bs_i,\bxi_m')\right)
>
nr.
\]
We conclude that at least one replaced point $\bs_i'=(\by'_i,\bX'_i)\in \mathcal{C}(\bbeta_m',\btheta_m',1)$.
Similarly, if all original points in~$\mathcal{S}'_m$ are outside $\mathcal{C}(\bbeta_m',\btheta_m',1)$, then
\[
\sum_{\bs_i\in \mathcal{S}_m'}
\rho\left(d(\bs_i,\bxi_m')\right)
\geq
n-\lceil nr\rceil
\geq
n-\lceil n-nr\rceil+\kappa(\mathcal{S}_n)
>nr,
\]
which is in contradiction with~\eqref{eq:constraint corrupted}.
Therefore, the cylinder $\mathcal{C}(\bbeta_m',\btheta_m',1)$
must contain a point $\bs_0=(\by_0,\bX_0)$ from the original collection~$\mathcal{S}_n$ as well as one replaced point $\bs_i'=(\by'_i,\bX'_i)$.

Note that for each point $\bs=(\by,\bX)\in \mathcal{C}(\bbeta_m',\btheta_m',1)$ it holds that
\begin{equation}
\label{eq:quadratic form}
\frac{|\bq_1^T(\by-\bX\bbeta_m')|^2}{\lambda_1(\bV(\btheta_m'))}
+\cdots+
\frac{|\bq_k^T(\by-\bX\bbeta_m')|^2}{\lambda_k(\bV(\btheta_m'))}
\leq
1,
\end{equation}
where $\lambda_j(\bV(\btheta_m'))>0$, for $j=1,\ldots,k$, are the eigenvalues of $\bV(\btheta_m')$
and $\bq_1,\ldots,\bq_k$ are the corresponding orthonormal eigenvectors.
Now, replace all $m$ points by $\bs_i'=(\by_i',\bX_i')=(\bz,\mathbf{0})$,
where~$\mathbf{0}$ is a $k\times q$ matrix of zeros and
$\bz=t\sum_{j=1}^k\bq_j$, so that $\bq_j^T\bz=t$, for each $j=1,\ldots,k$.
By sending $t\to\infty$ and the fact that at least one replaced point $(\bz,\mathbf{0})$
satisfies~\eqref{eq:quadratic form}, it follows that
$\lambda_j(\bV(\btheta_m'))\to\infty$, for each $j=1,\ldots,k$.
This means $\btheta_m'$ breaks down.

The upper bound for $\epsilon_n^*(\bbeta_n,\mathcal{S}_n)$ follows from the fact that
$\bbeta_n$ is regression equivariant.
Similar to Theorem~2 in~\cite{lopuhaa&rousseeuw1991} it can be shown that
the maximal breakdown point of regression equivariant estimators is $\lfloor(n+1)/2\rfloor/n$.
This proves the theorem.
\end{proof}

\subsection{Proofs for Section~\ref{sec:IF}}
\begin{lemma}
\label{lem:Psi bounded}
Suppose that $\rho$ satisfies (R2), (R4) and $\bV$ satisfies (V4).
Let $\Psi=(\Psi_{\bbeta},\Psi_{\btheta})$, as defined in~\eqref{eq:Psi function}.
Then there exist $0<C_1<\infty$ and $0<C_2<\infty$, only depending on $P$, such that
$\|\Psi_{\bbeta}(\bs,\bxi(P))\|\leq C_1\|\bX\|$
and~$\|\Psi_{\btheta}(\bs,\bxi(P))\|\leq C_2$.
\end{lemma}
\begin{proof}
Consider the expression of $\Psi_{\bbeta}$ in~\eqref{eq:Psi function}.
Consecutively, we apply~\eqref{eq:euclidean - trace},
\eqref{eq:largest eigenvalue},
\eqref{eq:euclidean - lambda1 symmetric},
\eqref{eq:product triangle spectral},
and
\eqref{eq:Frobenius - Spectral}.
This gives
\begin{equation}
\label{eq:bound XVa}
\begin{split}
\left\|
\bX^T\bV^{-1}(\by-\bX\bbeta)
\right\|^2
&=
(\by-\bX\bbeta)^T
\bV^{-1}
\bX\bX^T\bV^{-1}(\by-\bX\bbeta)\\
&\leq
\lambda_1(\bV^{-1/2}\bX\bX^T\bV^{-1/2})
(\by-\bX\bbeta)^T
\bV^{-1}(\by-\bX\bbeta)\\
&=
d^2\left\|\bV^{-1/2}\bX\bX^T\bV^{-1/2}\right\|_2\\
&\leq
d^2
\left\|\bV^{-1/2}\right\|_2^2
\left\|\bX\bX^T\right\|_2\\
&\leq
d^2\|\bX\|^2\lambda_1(\bV^{-1}),
\end{split}
\end{equation}
where
$d=d(\bs,\bxi)$ as defined by~\eqref{def:Mahalanobis distance structured},
and where we abbreviate~$\bV(\btheta)$ by $\bV$.
This means that
\[
\|\Psi_{\bbeta}(\bs,\bxi(P))\|
\leq
|u(d)d|
\times
\|\bX\|
\times
\sqrt{\lambda_1(\bV(\btheta(P))^{-1})}.
\]
From (R2) and (R4) it follows that $u(s)s=\rho'(s)$ is bounded.
This means that there exists a universal constant $0<C_1<\infty$, such that
\[
\|\Psi_{\bbeta}(\bs,\bxi(P))\|
\leq
C_1\|\bX\|.
\]
For $\Psi_{\btheta}=(\Psi_{\btheta,1},\cdots,\Psi_{\btheta,l})$, we have
\[
\Psi_{\btheta,j}(\bs,\bxi)
=
u(d)
(\by-\bX\bbeta)^T\bV^{-1}
\bH_j
\bV^{-1}(\by-\bX\bbeta)
-
\mathrm{tr}\left(\bV^{-1}\frac{\partial \bV}{\partial \theta_j}\right)(\rho(d)-b_0),
\]
for $j=1,\ldots,l$,
where we write $\bV$ instead of $\bV(\btheta)$.
Recall that
\[
\bH_j
=
\mathrm{tr}\left(\bV^{-1}\frac{\partial \bV}{\partial \theta_j}\right)
\left(
\sum_{t=1}^l\theta_t\frac{\partial \bV}{\partial \theta_t}
\right)
-
\mathrm{tr}\left(\bV^{-1}\sum_{t=1}^l\theta_t\frac{\partial \bV}{\partial \theta_t}\right)
\frac{\partial \bV}{\partial \theta_j}.
\]
To bound $\bH_j$, we first obtain a bound on
\begin{equation}
\label{eq:bound term Dj}
(\by-\bX\bbeta)^T\bV^{-1}
\frac{\partial \bV}{\partial \theta_j}
\bV^{-1}(\by-\bX\bbeta).
\end{equation}
Note that $\partial\bV/\partial\theta_j$ is symmetric, but not necessarily positive definite.
When~\eqref{eq:bound term Dj} is positive, then
application of~\eqref{eq:largest eigenvalue} gives
\[
\begin{split}
0<
(\by-\bX\bbeta)^T\bV^{-1}
\frac{\partial \bV}{\partial \theta_j}
\bV^{-1}(\by-\bX\bbeta)
&\leq
d^2
\lambda_1\left(
\bV^{-1/2}\frac{\partial \bV}{\partial \theta_j}\bV^{-1/2}
\right)\\
&\leq
d^2
\left\|\bV^{-1/2}\right\|^2
\left\|\frac{\partial \bV}{\partial \theta_j}\right\|\\
&\leq
d^2
\left\|\frac{\partial \bV}{\partial \theta_j}\right\|
\lambda_1(\bV^{-1}),
\end{split}
\]
according to~\eqref{eq:product triangle} and~\eqref{eq:euclidean - lambda1 symmetric}.
When~\eqref{eq:bound term Dj} is negative, then similarly
\[
\begin{split}
0<
(\by-\bX\bbeta)^T\bV^{-1}
\left(-\frac{\partial \bV}{\partial \theta_j}\right)
\bV^{-1}(\by-\bX\bbeta)
&\leq
d^2
\lambda_1\left(
\bV^{-1/2}\left(-\frac{\partial \bV}{\partial \theta_j}\right)\bV^{-1/2}
\right)\\
&\leq
d^2
\left\|\frac{\partial \bV}{\partial \theta_j}\right\|
\lambda_1(\bV^{-1}).
\end{split}
\]
It follows that
\begin{equation}
\label{eq:bound aMa}
\left|
(\by-\bX\bbeta)^T\bV^{-1}
\frac{\partial \bV}{\partial \theta_j}
\bV^{-1}(\by-\bX\bbeta)
\right|
\leq
d^2
\left\|\frac{\partial \bV}{\partial \theta_j}\right\|
\lambda_1(\bV^{-1}).
\end{equation}
Furthermore, according to~(V4),
the mapping $\btheta\mapsto\bV(\btheta)$ is continuously differentiable.
This means that
there exists a universal constant $0<M_1<\infty$, such that
\begin{equation}
\label{eq:bound max dV}
\max_{1\leq j\leq l}
\left\|
\frac{\partial \bV(\btheta(P))}{\partial \theta_j}\right\|\leq M_1.
\end{equation}
Because $\partial\bV/\partial\theta_j$ is symmetric, according to~\eqref{eq:euclidean - lambda1 symmetric}
and~\eqref{eq:product triangle spectral}
\[
\left|
\lambda_1\left(\bV^{-1}\frac{\partial \bV}{\partial \theta_j}\right)
\right|
\leq
\left\|\bV^{-1}\right\|_2
\left\|\frac{\partial \bV}{\partial \theta_j}\right\|_2
\leq
\sqrt{k}
\lambda_1\left(\bV^{-1}\right)
\left\|
\frac{\partial \bV}{\partial \theta_j}
\right\|.
\]
Together with~\eqref{eq:bound max dV},
we find
\begin{equation}
\label{eq:bound trace}
\left|
\mathrm{tr}\left(\bV^{-1}\frac{\partial \bV}{\partial \theta_j}\right)
\right|
\leq
k
\left|
\lambda_1\left(\bV^{-1}\frac{\partial \bV}{\partial \theta_j}\right)
\right|
\leq
k^{3/2}
M_1
\lambda_1\left(\bV^{-1}\right),
\end{equation}
where we abbreviate $\bV(\btheta(P))$ by $\bV$.
It then follows there exists a constant $0<M_2<\infty$,
only depending on $P$, such that at $\btheta(P)$,
\begin{equation}
\label{eq:bound Dj}
\max_{1\leq j\leq l}
\left|
(\by-\bX\bbeta)^T\bV^{-1}
\bH_j
\bV^{-1}(\by-\bX\bbeta)
\right|
\leq
d^2M_2.
\end{equation}
From (R2) and (R4) it follows that $u(s)s^2=\rho'(s)s$ and $\rho(s)-b_0$ are bounded.
Together with~\eqref{eq:bound trace},
it follows
that
there exists a universal constant $0<C_2<\infty$, such that
\[
\|\Psi_{\btheta,j}(\bs,\bxi(P))\|
\leq
C_2,
\]
for all $j=1,\ldots,l$.
This finishes the proof.
\end{proof}
\paragraph*{Proof of Corollary~\ref{cor:IF bounded}}
\begin{proof}
Take $\bs=(\by,\bX)$ fixed and consider $\mathrm{IF}(\bs;\bxi,P)$.
Since $\bD(P)$ does not depend on $\bs$, from Theorem~\ref{th:IF} and Lemma~\ref{lem:Psi bounded},
it follows immediately that $\text{IF}(\bs;\bxi,P)$ remains bounded in $\by$,
but not necessarily in $\bX$.
\end{proof}
\begin{lemma}
\label{lem:diff Lambda}
Consider $\Lambda$ as defined by~\eqref{def:Lambda} with~$\Psi$ defined in~\eqref{eq:Psi function}.
Suppose that $\rho$ satisfies (R2) and (R5) and $\bV$ satisfies (V5).
Furthermore, suppose that $\E\|\bX\|^2<\infty$.
Let~$\bxi(P)$ be a solution to~\eqref{def:Smin structured} and let $N$ be an open neighborhood of~$\bxi(P)$.
Then,
$\Lambda$ is continuous differentiable at $\bxi(P)$ and
for all $\bxi\in N$,
\[
\frac{\partial\Lambda(\bxi)}{\partial \bxi}
=
\int
\frac{\partial\Psi(\bs,\bxi)}{\partial \bxi}\,\dd P(\bs).
\]
\end{lemma}
\begin{proof}
Write $\partial \Lambda/\partial\bxi$ as the block matrix
\begin{equation}
\label{def:Lambda derivative}
\frac{\partial \Lambda(\bxi)}{\partial \bxi}
=
\left(
  \begin{array}{cc}
\dfrac{\partial \Lambda_{\bbeta}(\bxi)}{\partial \bbeta} & \dfrac{\partial \Lambda_{\bbeta}(\bxi)}{\partial \btheta} \\
    \\
\dfrac{\partial \Lambda_{\btheta}(\bxi)}{\partial \bbeta} & \dfrac{\partial \Lambda_{\btheta}(\bxi)}{\partial \btheta} \\
  \end{array}
\right),
\end{equation}
where
\[
\begin{split}
\Lambda_{\bbeta}(\bxi)
&=
\int
\Psi_{\bbeta}(\bs,\bxi)\,\dd P(\bs),\\
\Lambda_{\btheta}(\bxi)
&=
\int
\Psi_{\btheta}(\bs,\bxi)\,\dd P(\bs).
\end{split}
\]
We prove the lemma for each block separately.
Consider $\partial \Lambda_{\bbeta}/\partial\bbeta$.
We have
\begin{equation}
\label{eq:derivative psi beta wrt beta}
\frac{\partial\Psi_{\bbeta}(\bs,\bxi)}{\partial\bbeta}
=
-
\frac{u'(d)}{d}
\bX^T\bV^{-1}(\by-\bX\bbeta)
(\by-\bX\bbeta)^T
\bV^{-1}\bX
-
u(d)\bX^T\bV^{-1}\bX,
\end{equation}
where $d=d(\bs,\bxi)$ is defined by~\eqref{def:Mahalanobis distance structured}
and where we abbreviate~$\bV(\btheta)$ by $\bV$.
First note that, according to (R5), $\bxi\mapsto u'(d(\bs,\bxi))/d(\bs,\bxi)$ is continuous at $\bxi(P)$, for each $\bs$ fixed such that $d(\bs,\bxi(P))\ne0$.
Together with (V1), this means that for such $\bs$ fixed, $\bxi\mapsto \partial\Psi_{\bbeta}(\bs,\bxi)/\partial\bbeta$ is continuous at $\bxi(P)$.
For the first term on the right hand side of~\eqref{eq:derivative psi beta wrt beta},
we apply~\eqref{eq:euclidean uv} and~\eqref{eq:bound XVa}.
This gives
\[
\left\|
\bX^T\bV^{-1}(\by-\bX\bbeta)
(\by-\bX\bbeta)^T
\bV^{-1}\bX
\right\|
=
\left\|
\bX^T\bV^{-1}(\by-\bX\bbeta)
\right\|^2
\leq
d^2\|\bX\|^2\lambda_1(\bV^{-1}).
\]
Similarly, for the second term on the right hand side of~\eqref{eq:derivative psi beta wrt beta},
after application of~\eqref{eq:product triangle} and~\eqref{eq:euclidean - lambda1 symmetric}, we get
\[
\left\|
\bX^T\bV^{-1}\bX
\right\|
\leq
\left\|
\bX
\right\|^2
\left\|\bV^{-1}\right\|
\leq
\sqrt{q}\|\bX\|^2\lambda_1(\bV^{-1}).
\]
Since $\lambda_1(\bV^{-1})$ is bounded uniformly on the neighborhood $N$ of $\bxi(P)$ and
because $u(s)$ and $u'(s)s=\rho''(s)-u(s)$ are bounded,
due to (R2), it follows that there exists a constant $0<C_1<\infty$,
only depending on~$P$, such that
\[
\left\|
\frac{\partial\Psi_{\bbeta}(\bs,\bxi)}{\partial\bbeta^T}
\right\|
\leq
C_1
\|\bX\|^2.
\]
Since $\E\|\bX\|^2<\infty$, it follows by dominated convergence that for $\bxi$ in the neighborhood $N$
of~$\bxi(P)$, it holds that
\begin{equation}
\label{eq:derivative Lambda beta}
\frac{\partial\Lambda_{\bbeta}(\bxi)}{\partial\bbeta}
=
\int
\frac{\partial\Psi_{\bbeta}(\bs,\bxi)}{\partial \bbeta}\,\dd P(\bs),
\end{equation}
and that $\partial\Lambda_{\bbeta}/\partial\bbeta$ is continuous at $\bxi(P)$.

Next consider $\partial\Psi_{\bbeta}/\partial\btheta$.
For each $j=1,\ldots,l$ fixed, we have
\begin{equation}
\label{eq:derivative psi beta wrt theta}
\begin{split}
\frac{\partial\Psi_{\bbeta}(\bs,\bxi)}{\partial\theta_j}
&=
\frac{u'(d)}{2d}
(\by-\bX\bbeta)^T
\bV^{-1}\frac{\partial\bV}{\partial\theta_j}
\bV^{-1}(\by-\bX\bbeta)
\cdot \bX^T\bV^{-1}(\by-\bX\bbeta)\\
&\quad+
u(d)\cdot
\bX^T\bV^{-1}\frac{\partial\bV}{\partial\theta_j}\bV^{-1}(\by-\bX\bbeta).
\end{split}
\end{equation}
First note that, similar to~\eqref{eq:derivative psi beta wrt beta},
$\bxi\mapsto \partial\Psi_{\bbeta}(\bs,\bxi)/\partial\theta_j$ is continuous at $\bxi(P)$,
for each $\bs$ fixed such that $d(\bs,\bxi(P))\ne0$.
Consider the first term on the right hand side of~\eqref{eq:derivative psi beta wrt theta}.
From~\eqref{eq:bound XVa}, we have
\[
\|\bX^T\bV^{-1}(\by-\bX\bbeta)\|\leq d\|\bX\|\sqrt{\lambda_1(\bV^{-1})}.
\]
Moreover, similar to the reasoning in~\eqref{eq:bound aMa}, we find
\begin{equation}
\label{eq:bound aVMVa}
\left|
(\by-\bX\bbeta)^T
\bV^{-1}
\frac{\partial\bV}{\partial\theta_j}
\bV^{-1}(\by-\bX\bbeta)
\right|
\leq
d^2
\left\|
\frac{\partial\bV}{\partial\theta_j}
\right\|
\lambda_1(\bV^{-1}).
\end{equation}
For the second term on the right hand side of~\eqref{eq:derivative psi beta wrt theta},
similar to the reasoning in~\eqref{eq:bound XVa}, we have
\[
\begin{split}
\left\|
\bX^T\bV^{-1}\frac{\partial\bV}{\partial\theta_j}\bV^{-1}(\by-\bX\bbeta)
\right\|^2
&=
(\by-\bX\bbeta)\bV^{-1}\frac{\partial\bV}{\partial\theta_j}\bV^{-1}\bX
\bX^T\bV^{-1}\frac{\partial\bV}{\partial\theta_j}\bV^{-1}(\by-\bX\bbeta)\\
&\quad\leq
d^2
\lambda_1
\left(
\bV^{-1/2}\frac{\partial\bV}{\partial\theta_j}\bV^{-1}\bX
\bX^T\bV^{-1}\frac{\partial\bV}{\partial\theta_j}\bV^{-1/2}
\right)\\
&\quad\leq
d^2
\left\|
\bV^{-1/2}\frac{\partial\bV}{\partial\theta_j}\bV^{-1}\bX
\bX^T\bV^{-1}\frac{\partial\bV}{\partial\theta_j}\bV^{-1/2}
\right\|_2\\
&\quad\leq
d^2
\|\bX\|^2
\left\|\frac{\partial\bV}{\partial\theta_j}\right\|^2
\lambda_1(\bV^{-1})^3.
\end{split}
\]
According to (V4), the mapping $\bV(\btheta)$ is continuously differentiable.
This means that $\left\|\partial\bV/\partial\theta_j\right\|$ is bounded on the neighborhood~$N$ of~$\bxi(P)$.
Since $\lambda_1(\bV^{-1})$ is bounded uniformly on $N$ and
because~$u(s)s=\rho'(s)$ and~$u'(s)s^2=\rho''(s)s-\rho'(s)$ are bounded,
it follows that there exists a
constant $0<C_2<\infty$, only depending on $P$, such that
for $j=1,\ldots,l$,
\[
\left\|
\frac{\partial\Psi_{\bbeta}(\bs,\bxi)}{\partial\theta_j}
\right\|
\leq
C_2
\|\bX\|^2.
\]
As before, it follows by dominated convergence that for $\bxi$ in the neighborhood $N$ of~$\bxi(P)$, it holds that
\begin{equation}
\label{eq:derivative Lambda beta}
\frac{\partial\Lambda_{\bbeta}(\bxi)}{\partial\btheta}
=
\int
\frac{\partial\Psi_{\bbeta}(\bs,\bxi)}{\partial \btheta}\,\dd P(\bs),
\end{equation}
and that $\partial\Lambda_{\bbeta}/\partial\btheta$ is continuous at $\bxi(P)$.

Next consider $\partial\Psi_{\btheta,j}/\partial\bbeta$, for $j=1,\ldots,l$.
We have
\begin{equation}
\label{eq:derivative psi theta wrt beta}
\begin{split}
\frac{\partial\Psi_{\btheta,j}(\bs,\bxi)}{\partial\bbeta}
&=
\frac{u'(d)}{d}
\bX^T\bV^{-1}(\by-\bX\bbeta)
\cdot
(\by-\bX\bbeta)^T\bV^{-1}
\bH_j
\bV^{-1}(\by-\bX\bbeta)
\\
&\quad-
u(d)
\cdot
2\bX^T\bV^{-1}\bH_j
\bV^{-1}(\by-\bX\bbeta)\\
&\quad+
\mathrm{tr}\left(\bV^{-1}\frac{\partial \bV}{\partial \theta_j}\right)
u(d)
\bX^T\bV^{-1}(\by-\bX\bbeta),
\end{split}
\end{equation}
where $\bH_j$ is defined in~\eqref{def:Hj}.
As before,
$\bxi\mapsto \partial\Psi_{\theta_j}(\bs,\bxi)/\partial\bbeta$ is continuous at $\bxi(P)$,
for each $\bs$ fixed such that $d(\bs,\bxi(P))\ne0$.
Consider the first term on the right hand side of~\eqref{eq:derivative psi theta wrt beta}.
From~\eqref{eq:bound Dj},
\begin{equation}
\label{eq:bound aMaDj}
\left|
(\by-\bX\bbeta)^T\bV^{-1}
\bH_j
\bV^{-1}(\by-\bX\bbeta)
\right|
\leq
d^2
M_2.
\end{equation}
Because $u'(s)s^2$ is bounded,
together with~\eqref{eq:bound XVa},
the norm of the first term on the right hand side of~\eqref{eq:derivative psi theta wrt beta} is bounded
by a constant times $\|\bX\|\lambda_1(\bV^{-1})^{1/2}$.
Similar to~\eqref{eq:bound XVa}, for the second term on the right hand side of~\eqref{eq:derivative psi theta wrt beta},
\[
\begin{split}
\|
\bX^T\bV^{-1}\bH_j\bV^{-1}(\by-\bX\bbeta)
\|^2
&\leq
d^2
\left\|\bV^{-1/2}\bH_j\bV^{-1}\right\|_2^2
\left\|\bX\bX^T\right\|_2^2\\
&\leq
kd^2
\|\bX\|^2
\|\bH_j\|^2
\lambda_1(\bV^{-1})^3,
\end{split}
\]
and for the third term on the right hand side of~\eqref{eq:derivative psi theta wrt beta},
we can use~\eqref{eq:bound XVa} and~\eqref{eq:bound trace}.
As before, since $u'(s)s=\rho''(s)-u(s)$ and $u(s)s^2=\rho'(s)s$ are bounded, it follows that there exists a
constant $0<C_3<\infty$, only depending on $P$,
such that for $j=1,\ldots,l$,
\[
\left\|
\frac{\partial\Psi_{\btheta,j}(\bs,\bxi)}{\partial\bbeta}
\right\|
\leq
C_3\|\bX\|.
\]
Since $\E\|\bX\|<\infty$, if follows by dominated convergence that for $\bxi$ in the neighborhood~$N$ of~$\bxi(P)$, it holds that
\begin{equation}
\label{eq:derivative Lambda theta beta}
\frac{\partial\Lambda_{\btheta}(\bxi)}{\partial\bbeta}
=
\int
\frac{\partial\Psi_{\btheta}(\bs,\bxi)}{\partial \bbeta}\,\dd P(\bs),
\end{equation}
and that $\partial\Lambda_{\btheta}/\partial\bbeta$ is continuous at $\bxi(P)$.

Finally, consider $\partial\Psi_{\btheta,j}/\partial\theta_t$, for $j,t=1,\ldots,l$.
We find
\begin{equation}
\label{eq:derivative psi theta wrt theta}
\begin{split}
\frac{\partial\Psi_{\btheta,j}}{\partial\theta_t}
=&-
\frac{u'(d)}{2d}
(\by-\bX\bbeta)^T\bV^{-1}
\frac{\partial\bV}{\partial\theta_t}\bV^{-1}(\by-\bX\bbeta)\\
&\qquad\qquad\qquad\qquad\qquad\qquad
\cdot
(\by-\bX\bbeta)^T\bV^{-1}\bH_j\bV^{-1}(\by-\bX\bbeta)
\\
&\quad-
u(d)
(\by-\bX\bbeta)^T\bV^{-1}\frac{\partial\bV}{\partial\theta_t}\bV^{-1}
\cdot
\bH_j\bV^{-1}(\by-\bX\bbeta)\\
&\quad+
u(d)
(\by-\bX\bbeta)^T\bV^{-1}\frac{\partial\bH_j}{\partial\theta_t}
\cdot
\bV^{-1}(\by-\bX\bbeta)\\
&\quad-
u(d)
(\by-\bX\bbeta)^T\bV^{-1}\bH_j
\cdot
\bV^{-1}\frac{\partial\bV}{\partial\theta_t}\bV^{-1}(\by-\bX\bbeta)\\
&\qquad+
\mathrm{tr}\left(\bV^{-1}\frac{\partial\bV}{\partial\theta_t}\bV^{-1}
\cdot
\frac{\partial \bV}{\partial \theta_j}\right)
\left(\rho(d)-b_0\right)\\
&\qquad-
\mathrm{tr}\left(\bV^{-1}
\cdot
\left(\frac{\partial^2\bV}{\partial\theta_j\theta_t}\right)\right)
\left(\rho(d)-b_0\right)\\
&\qquad+
\mathrm{tr}\left(\bV^{-1}\frac{\partial \bV}{\partial \theta_j}\right)
\cdot
\frac{u(d)}{2}
(\by-\bX\bbeta)^T\bV^{-1}\frac{\partial\bV}{\partial\theta_t}\bV^{-1}(\by-\bX\bbeta).
\end{split}
\end{equation}
As before, together with (V5),
$\bxi\mapsto \partial\Psi_{\theta_j}(\bs,\bxi)/\partial\theta_t$ is continuous at $\bxi(P)$,
for each $\bs$ fixed such that $d(\bs,\bxi(P))\ne0$.
From~\eqref{eq:bound aMa} and~\eqref{eq:bound Dj}, it follows that the first term
on the right hand side
of~\eqref{eq:derivative psi theta wrt theta} is bounded by
\[
\frac{|u'(d)d^3|}{2}
\left\|
\frac{\partial\bV}{\partial\theta_t}
\right\|
\lambda_1(\bV^{-1})
M_2.
\]
Because $u'(s)s^3=\rho''(s)s^2-\rho'(s)s$ is bounded,
together with~\eqref{eq:bound max dV}, we conclude this is
bounded on the neighborhood $N$ of $\bxi(P)$.
Similar to~\eqref{eq:bound XVa},
the second term on the right hand side
of~\eqref{eq:derivative psi theta wrt theta} is bounded by
\[
|u(d)|d^2
\left\|
\bV^{-1/2}\frac{\partial\bV}{\partial\theta_t}\bV^{-1}\bH_j\bV^{-1/2}
\right\|_2^2
\leq
|u(d)|d^2
\left\|
\frac{\partial\bV}{\partial\theta_t}
\right\|^2
\left\|
\bH_j
\right\|^2
\lambda_1\left(\bV^{-1}\right)^4.
\]
Since $u(s)s^2=\rho'(s)s$ is bounded,
together with~\eqref{eq:bound trace},
we again find that this is bounded on the neighborhood $N$ of $\bxi(P)$.
The same holds for the fourth term on the right hand side
of~\eqref{eq:derivative psi theta wrt theta}.

We continue with the third term on the right hand side
of~\eqref{eq:derivative psi theta wrt theta}.
Similar to~\eqref{eq:bound XVa},
this is bounded by
\[
|u(d)|d^2
\left\|
\bV^{-1/2}\frac{\partial\bH_j}{\partial\theta_t}\bV^{-1/2}
\right\|_2^2
\leq
|u(d)|d^2
\lambda_1\left(\bV^{-1}\right)^2
\left\|
\frac{\partial\bH_j}{\partial\theta_t}
\right\|.
\]
We have that
\[
\begin{split}
\frac{\partial\bH_j}{\partial\theta_t}
&=
\left\{
\mathrm{tr}\left(-\bV^{-1}\frac{\partial\bV}{\partial\theta_t}\bV^{-1}\frac{\partial \bV}{\partial \theta_j}\right)
+
\mathrm{tr}\left(\bV^{-1}\frac{\partial^2 \bV}{\partial \theta_j\theta_t}\right)
\right\}
\left(
\sum_{s=1}^l\theta_s\frac{\partial \bV}{\partial \theta_s}
\right)\\
&\quad+
\mathrm{tr}\left(\bV^{-1}\frac{\partial\bV}{\partial \theta_j}\right)
\left\{
\left(
\frac{\partial \bV}{\partial \theta_t}
\right)
+
\left(
\sum_{s=1}^l\theta_s\frac{\partial^2 \bV}{\partial \theta_s\theta_t}
\right)
\right\}
\\
&\quad
-
\left\{
\mathrm{tr}\left(
-\bV^{-1}\frac{\partial\bV}{\partial\theta_t}\bV^{-1}
\sum_{s=1}^l\theta_s\frac{\partial \bV}{\partial \theta_s}\right)
+
\mathrm{tr}\left(\bV^{-1}
\frac{\partial \bV}{\partial \theta_t}\right)
+
\mathrm{tr}\left(\bV^{-1}
\sum_{s=1}^l\theta_s\frac{\partial^2 \bV}{\partial \theta_s\theta_t}\right)
\right\}
\frac{\partial \bV}{\partial \theta_j}\\
&\quad-
\mathrm{tr}\left(\bV^{-1}
\sum_{s=1}^l\theta_s\frac{\partial \bV}{\partial \theta_s}\right)
\left(\frac{\partial^2 \bV}{\partial \theta_j\theta_t}\right).
\end{split}
\]
With~\eqref{eq:euclidean - lambda1 symmetric}, \eqref{eq:product triangle spectral}, and~\eqref{eq:bound max dV}, we find
\begin{equation}
\label{eq:bound term1 Dj}
\begin{split}
\left|
\mathrm{tr}\left(\bV^{-1}\frac{\partial\bV}{\partial\theta_t}\bV^{-1}
\frac{\partial \bV}{\partial \theta_s}\right)
\right|
&\leq
k
\left\|
\bV^{-1}\frac{\partial\bV}{\partial\theta_t}\bV^{-1}
\frac{\partial \bV}{\partial \theta_s}
\right\|_2
\\
&\leq
k
\left\|
\bV^{-1}
\right\|_2^2
\left\|\frac{\partial\bV}{\partial\theta_t}\right\|
\left\|\frac{\partial \bV}{\partial \theta_s}\right\|\\
&\leq
k\lambda_1(\bV^{-1})M_1^2,
\end{split}
\end{equation}
which is uniformly bounded on the neighborhood $N$ of $\bxi(P)$,
and similar to~\eqref{eq:bound trace} we find
\begin{equation}
\label{eq:bound term2 Dj}
\left|
\mathrm{tr}\left(\bV^{-1}
\frac{\partial^2\bV}{\partial\theta_j\theta_t}\right)
\right|
\leq
k^{3/2}
\lambda_1(\bV^{-1})
\left\|
\frac{\partial^2\bV}{\partial\theta_j\theta_t}\right\|.
\end{equation}
Because, according to~(V5), the mapping $\btheta\mapsto\bV(\btheta)$ is twice continuously differentiable,
it follows that
$\left\|\partial^2 \bV/\partial \theta_j\theta_t\right\|$
is uniformly bounded on the neighborhood $N$ of $\bxi(P)$.
Together with the fact that with~\eqref{eq:bound max dV},
\[
\left\|
\sum_{s=1}^l\theta_s\frac{\partial \bV}{\partial \theta_s}
\right\|
\leq
\sum_{s=1}^l
\|\theta_s\|
\left\|
\frac{\partial \bV}{\partial \theta_s}
\right\|
\leq
M_1
\sum_{s=1}^l
\|\theta_s\|,
\]
it follows that the first term of $\partial \bH_j/\partial \theta_t$
is bounded on the neighborhood $N$ of $\bxi(P)$.
The traces in the other terms can be handled in the same way, which yields that
$\|\partial\bH_j/\partial\theta_t\|$ is bounded on the neighborhood $N$ of $\bxi(P)$.
Because $u(s)s^2=\rho'(s)s$ is bounded, it follows that the third term
on the right hand side of~\eqref{eq:derivative psi theta wrt theta} is bounded.

Next, consider the fifth term on the right hand side of~\eqref{eq:derivative psi theta wrt theta}.
From~\eqref{eq:bound term1 Dj},
\[
\left|
\mathrm{tr}\left(\bV^{-1}\frac{\partial\bV}{\partial\theta_t}\bV^{-1}
\frac{\partial \bV}{\partial \theta_j}\right)
\right|
\leq
k\lambda_1(\bV^{-1})M_1^2,
\]
which is uniformly bounded on the neighborhood $N$ of $\bxi(P)$.
Because $\rho(s)$ is bounded, it follows that
the fifth term on the right hand side of~\eqref{eq:derivative psi theta wrt theta} is bounded.
For the sixth term on the right hand side of~\eqref{eq:derivative psi theta wrt theta},
from~\eqref{eq:bound term2 Dj} we have
\[
\left|
\mathrm{tr}\left(\bV^{-1}
\frac{\partial^2\bV}{\partial\theta_j\theta_t}\right)
\right|
\leq
k^{3/2}
\lambda_1(\bV^{-1})
\left\|
\frac{\partial^2\bV}{\partial\theta_j\theta_t}\right\|.
\]
Because $\bV(\btheta)$ is twice continuously differentiable and $\rho(s)$ is bounded,
we conclude that
the sixth term on the right hand side of~\eqref{eq:derivative psi theta wrt theta} is bounded.
Finally, from~\eqref{eq:bound trace} and~\eqref{eq:bound aVMVa} together with the fact that
$u(s)s^2=\rho'(s)s$ is bounded,
it also follows that
the last term on the right hand side of~\eqref{eq:derivative psi theta wrt theta} is bounded.
By putting everything together, it follows that there exists a
constant $0<C_4<\infty$, only depending on $P$,
such that
for $j,t=1,\ldots,l$,
\[
\left\|
\frac{\partial\Psi_{\btheta,j}(\bs,\bxi)}{\partial\theta_t}
\right\|
\leq
C_4.
\]
It follows by dominated convergence that for $\bxi$ in the neighborhood $N$ of $\bxi(P)$, it holds
\begin{equation}
\label{eq:derivative Lambda theta theta}
\frac{\partial\Lambda_{\btheta}(\bxi)}{\partial\btheta}
=
\int
\frac{\partial\Psi_{\btheta}(\bs,\bxi)}{\partial \btheta}\,\dd P(\bs),
\end{equation}
and that $\partial\Lambda_{\btheta}/\partial\btheta$ is continuous at $\bxi(P)$.
This finishes the proof.
\end{proof}
For convenience we state the following result from~\cite{lopuhaa1989}
about spherically contoured densities, see Lemma~5.1 in~\cite{lopuhaa1989}.
This lemma uses the commutation matrix $\mathbf{K}_{k,k}$, which
is the $k^2\times k^2$ block matrix with the $(i,j)$-block being equal to the $k\times k$ matrix $\mathbf{\Delta}_{ji}$
consisting of zero's except a 1 at entry $(j,i)$.
A useful property (e.g., see~\cite[Section 3.7]{magnus&neudecker1988}) is that for any $k\times k$ matrix $\bA$,
it holds that
\begin{equation}
\label{eq:prop K}
\bK_{k,k}\vc(\bA)=\vc(\bA^T).
\end{equation}
\begin{lemma}
\label{lem:Lemma 5.1}
Suppose that $\bz$ has a $k$-variate elliptical contoured density
defined in~\eqref{eq:elliptical}, with parameters $\bmu=\mathbf{0}$ and $\bSigma=\bI_k$.
Then $\bu=\bz/\|\bz\|$ is independent of $\|\bz\|$,
has mean zero and covariance matrix $(1/k)\bI_k$.
Furthermore, $\mathbb{E}_{\mathbf{0},\bI_k}\bu\bu^T\bu=\mathbf{0}$
and
\[
\mathbb{E}_{\mathbf{0},\bI_k}\vc(\bu\bu^T)\vc(\bu\bu^T)^T
=
\sigma_1(\bI_{k^2}+\mathbf{K}_{k,k})+\sigma_2\vc(\bI_k)\vc(\bI_k)^T,
\]
where $\sigma_1=\sigma_2=(k(k+2))^{-1}$.
\end{lemma}
\paragraph*{Proof of Lemma~\ref{lem:block derivative}}
\begin{proof}
Write $\partial\Lambda/\partial\bxi$ as in~\eqref{def:Lambda derivative}.
We determine each block separately and apply Lemma~\ref{lem:diff Lambda}.
Let~$\bxi_P=(\bbeta_P,\btheta_P)=(\bbeta(P),\btheta(P))$.
When we also write $\bV_P$ instead of $\bV(\btheta(P))$, then according to Lemma~\ref{lem:diff Lambda}
we have
\[
\begin{split}
\frac{\partial\Lambda_{\bbeta}(\bxi_P)}{\partial \bbeta}
&=
\int
\frac{\partial\Psi_{\bbeta}(\bs,\bxi_P)}{\partial \bbeta}\,\dd P(\bs)\\
&=
-\E
\left[
\frac{u'(d)}{d}
\bX^T\bV_P^{-1}(\by-\bX\bbeta_P)
(\by-\bX\bbeta_P)^T
\bV_P^{-1}\bX
+
u(d)\cdot\bX^T\bV_P^{-1}\bX
\right]\\
&=
-\E
\left[
\E
\left[
\frac{u'(d)}{d}
\bX^T\bSigma^{-1}(\by-\bmu)
(\by-\bmu)^T
\bSigma^{-1}\bX
+
u(d)\cdot\bX^T\bSigma^{-1}\bX
\,\Bigg|\,
\bX
\right]
\right],\\
\end{split}
\]
where $d^2=(\by-\bX\bbeta_P)^T\bV_P^{-1}(\by-\bX\bbeta_P)$.
The inner expectation on the right hand side is the conditional expectation of $\by\mid\bX$,
which has the same distribution as $\bSigma^{1/2}\bz+\bmu$,
where $\bz$ has a spherical density $f_{\mathbf{0},\bI_k}$.
This implies that the inner expectation on the right hand side is equal to
\[
\bX^T
\bSigma^{-1/2}
\mathbb{E}_{\mathbf{0},\bI_k}
\left[
\frac{u'(\|\bz\|)}{\|\bz\|}
\bz\bz^T
+
u(\|\bz\|)
\bI_k
\right]
\bSigma^{-1/2}
\bX.
\]
From Lemma~\ref{lem:Lemma 5.1}, we find
\[
\begin{split}
\mathbb{E}_{\mathbf{0},\bI_k}
\left[
\frac{u'(\|\bz\|)}{\|\bz\|}
\bz\bz^T
+
u(\|\bz\|)
\bI_k
\right]
&=
\mathbb{E}_{\mathbf{0},\bI_k}
\left[
u'(\|\bz\|)\|\bz\|
\frac{\bz\bz^T}{\|\bz\|^2}
+
u(\|\bz\|)
\bI_k
\right]\\
&=
\mathbb{E}_{\mathbf{0},\bI_k}
\Big[
u'(\|\bz\|)\|\bz\|
\Big]
\mathbb{E}_{\mathbf{0},\bI_k}
\left[
\bu\bu^T
\right]
+
\mathbb{E}_{\mathbf{0},\bI_k}
\left[
u(\|\bz\|)
\right]\bI_k
\\
&=
\mathbb{E}_{\mathbf{0},\bI_k}
\left[
u'(\|\bz\|)\|\bz\|
\right]
\frac1k \bI_k
+
\mathbb{E}_{\mathbf{0},\bI_k}
\left[
u(\|\bz\|)
\right]\bI_k
\\
&=
\alpha \bI_k,
\end{split}
\]
where $\bu=\bz/\|\bz\|$ and
\[
\alpha
=
\mathbb{E}_{\mathbf{0},\bI_k}
\left[
\frac1{k}
u'(\|\bz\|)\|\bz\|
+
u(\|\bz\|)
\right]
=
\mathbb{E}_{\mathbf{0},\bI_k}
\left[
\left(1-\frac{1}{k}\right)
\frac{\rho'(\|\bz\|)}{\|\bz\|}
+
\frac1k
\rho''(\|\bz\|)
\right].
\]
It follows that
\[
\frac{\partial\Lambda_{\bbeta}(\bxi_P)}{\partial \bbeta}
=
-\alpha
\mathbb{E}\left[\mathbf{X}^T\bSigma^{-1}\mathbf{X}\right].
\]
Next, for $j=1,\ldots,l$, consider
\begin{equation}
\label{eq:derivative Lambda-beta wrt theta}
\begin{split}
&
\frac{\partial\Lambda_{\bbeta}(\bxi_P)}{\partial\theta_j}
=
\int
\frac{\partial\Psi_{\bbeta}(\bs,\bxi_P)}{\partial\theta_j}
\,\dd P(\bs)\\
&=
\E\left[
\E\left[
\frac{u'(d)}{2d}
(\by-\bX\bbeta_P)^T
\bV_P^{-1}\frac{\partial\bV_P}{\partial\theta_j}
\bV_P^{-1}(\by-\bX\bbeta_P)
\cdot
\bX^T\bV_P^{-1}(\by-\bX\bbeta_P)
\,\Bigg|\,
\bX
\right]
\right]\\
&\qquad\qquad+
\E\left[
\E\left[
u(d)
\cdot
\bX^T\bV_P^{-1}\frac{\partial\bV_P}{\partial\theta_j}
\bV_P^{-1}(\by-\bX\bbeta_P)
\,\Bigg|\,
\bX
\right]
\right].
\end{split}
\end{equation}
According to Lemma~\ref{lem:Lemma 5.1}, the inner conditional expectation of the first term
on the right hand side of~\eqref{eq:derivative Lambda-beta wrt theta} can be written as
\[
\begin{split}
&
\bX^T
\E_{\mathbf{0},\bI_k}
\left[
\frac{u'(\|\bz\|)}{2\|\bz\|}
\bz^T
\bSigma^{-1/2}
\frac{\partial\bV_P}{\partial\theta_j}
\bSigma^{-1/2}
\bz
\bSigma^{-1/2}
\bz
\right]\\
&=
\bX^T
\E_{\mathbf{0},\bI_k}
\left[
\frac{u'(\|\bz\|)\|\bz\|^2}{2}
\right]
\E_{\mathbf{0},\bI_k}\left[
\bu^T
\bSigma^{-1/2}
\frac{\partial\bV_P}{\partial\theta_j}
\bSigma^{-1/2}
\bu
\bSigma^{-1/2}\bu
\right],
\end{split}
\]
Since the second term on the right hand side is the expectation with respect to a spherical density of an odd function of $\bu$,
this expectation is equal
to zero due to Lemma~\ref{lem:Lemma 5.1}.
Similarly, the second term on the right hand side of~\eqref{eq:derivative Lambda-beta wrt theta}
has inner conditional expectation
\[
\bX^T
\E_{\mathbf{0},\bI_k}\left[
u(\|\bz\|)
\bSigma^{-1}
\frac{\partial\bV_P}{\partial\theta_j}
\bSigma^{-1/2}\bz
\right]
=
\bX^T
\E_{\mathbf{0},\bI_k}\left[
u(\|\bz\|)\|\bz\|
\right]
\bSigma^{-1}
\frac{\partial\bV_P}{\partial\theta_j}
\bSigma^{-1/2}
\E_{\mathbf{0},\bI_k}\left[
\bu
\right]
=\mathbf{0},
\]
due to Lemma~\ref{lem:Lemma 5.1}.
It follows that
\[
\frac{\partial\Lambda_{\bbeta}(\bxi_P)}{\partial\btheta}=\mathbf{0}.
\]
Next, using Lemma~\ref{lem:diff Lambda} and~\eqref{eq:derivative psi theta wrt beta},
for all $j=1,\ldots,l$, consider
\begin{equation}
\label{eq:derivative Lambda-theta wrt beta}
\begin{split}
&
\frac{\partial \Lambda_{\btheta,j}}{\partial \bbeta^T}
=
\int
\frac{\partial \Psi_{\btheta,j}(\bs,\bxi_P)}{\partial \bbeta^T}
\,\dd P(\bs)\\
&=
-
\E\left[
\E\left[
\frac{u'(d)}{d}
\bX^T\bV_P^{-1}(\by-\bX\bbeta_P)
\cdot
\,(\by-\bX\bbeta_P)^T
\bV_P^{-1}\bH_j\bV_P^{-1}
(\by-\bX\bbeta_P)
\,\Bigg|\,
\bX
\right]
\right]\\
&\quad-
\E\left[
\E\left[
u(d)
\cdot
2
\bX^T\bV_P^{-1}\bH_j\bV_P^{-1}(\by-\bX\bbeta_P)
\,\Bigg|\,
\bX
\right]
\right]\\
&\quad-
\mathrm{tr}\left(\bV^{-1}\frac{\partial \bV}{\partial \theta_j}\right)
\E\left[
\E\left[
\frac{\rho'(d)}{d}
\bX^T\bV_P^{-1}(\by-\bX\bbeta_P)
\,\Bigg|\,
\bX
\right]
\right].
\end{split}
\end{equation}
According to Lemma~\ref{lem:Lemma 5.1},
the first term on the right hand side of~\eqref{eq:derivative Lambda-theta wrt beta}
has inner conditional expectation
\[
\begin{split}
&
\bX^T
\E_{\mathbf{0},\bI_k}
\left[
\frac{u'(\|\bz\|)}{\|\bz\|}
\bSigma^{-1/2}
\bz\bz^T\bSigma^{-1/2}
\bH_j
\bSigma^{-1/2}\bz
\right]\\
&\qquad=
\bX^T
\E_{\mathbf{0},\bI_k}
\left[
u'(\|\bz\|)\|\bz\|^2
\right]
\bSigma^{-1/2}
\E_{\mathbf{0},\bI_k}
\left[
\bu\bu^T
\bSigma^{-1/2}
\bH_j
\bSigma^{-1/2}
\bu
\right].
\end{split}
\]
Again, the second term on the right hand side is the expectation with respect to a spherical density
of an odd function of $\bu$, and is therefore equal to zero.
Similarly, the inner expectation of the second
term on the right hand side of~\eqref{eq:derivative Lambda-theta wrt beta} is equal to
\[
2
\bX^T
\bSigma^{-1}
\bH_j
\bSigma^{-1/2}
\E_{\mathbf{0},\bI_k}
\left[
u(\|\bz\|)
\bz
\right]
=
2
\bX^T
\bSigma^{-1}
\bH_j
\bSigma^{-1/2}
\E_{\mathbf{0},\bI_k}
\left[
u(\|\bz\|)\|\bz\|
\right]
\E_{\mathbf{0},\bI_k}
\left[
\bu
\right]
=\mathbf{0},
\]
and the inner expectation of the third
term on the right hand side of~\eqref{eq:derivative Lambda-theta wrt beta} is equal to
\[
\bX^T
\bSigma^{-1/2}
\E_{\mathbf{0},\bI_k}
\left[
\frac{\rho'(\|\bz\|)}{\|\bz\|}
\bz
\right]
=
\bX^T
\bSigma^{-1/2}
\E_{\mathbf{0},\bI_k}
\left[
\rho'(\|\bz\|)
\right]
\E_{\mathbf{0},\bI_k}
\left[
\bu
\right]
=\mathbf{0}.
\]
It follows that
\[
\frac{\partial\Lambda_{\btheta}(\bxi_P)}{\partial\bbeta}=\mathbf{0}.
\]
Finally,
to determine $\partial\Lambda_{\btheta,j}(\bxi_P)/\partial \theta_s$,
note that when $\bV$ is linear, we can write
\[
\Psi_{\btheta,j}(\bs,\bxi)
=
-
\vc\left(
\bV^{-1}
\bL_j
\bV^{-1}
\right)^T
\vc\left(\Psi_\bV(\bs,\bxi)\right)
\]
where $\Psi_\bV$ is defined in~\eqref{def:PsiV} and has the property that
\[
\int \Psi_\bV(\bs,\bxi_P)\,\dd P(\bs)=\mathbf{0},
\]
when $\by\mid\bX$ has an elliptically contoured density $f_{\bmu,\bSigma}$ with parameters $\bmu=\bX\bbeta_P$ and
$\bSigma=\bV_P$.
This means that
for each $j,s=1,\ldots,l$, we have
\[
\frac{\partial\Lambda_{\btheta,j}(\bxi_P)}{\partial \theta_s}
=
\int
\frac{\partial\Psi_{\btheta,j}(\bs,\bxi_P)}{\partial \theta_s}\,\dd P(\bs)
=
-\vc(\bSigma^{-1}\bL_j\bSigma^{-1})^T\vc\left(
\int
\frac{\partial\Psi_{\bV}(\bs,\bxi_P)}{\partial \theta_s}
\,\dd P(\bs)
\right)
\]
where $\Psi_{\bV}$ is defined in~\eqref{def:PsiV}.
As before, we find
\begin{equation}
\label{eq:derivative LambdaV-thetaj wrt thetas}
\begin{split}
\int
\frac{\partial\Psi_{\bV}(\bs,\bxi_P)}{\partial \theta_s}
\,\dd P(\bs)
=&
-\E
\left[
\E
\left[
\frac{ku'(d)}{2d}
\bz^T
\bSigma^{-1/2}
\bL_s
\bSigma^{-1/2}
\bz
\cdot
\bSigma^{1/2}\bz
\bz^T
\bSigma^{1/2}
\,\Bigg|\,
\bX
\right]
\right]\\
&\quad
+
\E
\left[
\E
\left[
\frac{v'(d)}{2d}
\bz^T
\bSigma^{-1/2}
\bL_s
\bSigma^{-1/2}\bz
\cdot
\bSigma
\,\Bigg|\,
\bX
\right]
\right]\\
&\qquad-
\E
\left[
\E
\left[
v(d)\bL_s
\,\Bigg|\,
\bX
\right]
\right].
\end{split}
\end{equation}
The first term on the right hand side of~\eqref{eq:derivative LambdaV-thetaj wrt thetas} is equal to
\[
\begin{split}
&
\E_{\mathbf{0},\bI_k}
\left[
\frac{ku'(\|\bz\|)}{2\|\bz\|}
\bz^T
\bSigma^{-1/2}\bL_s
\bSigma^{-1/2}
\bz
\bSigma^{1/2}\bz\bz^T\bSigma^{1/2}
\right]\\
&\qquad=
\E_{\mathbf{0},\bI_k}
\left[
\frac{ku'(\|\bz\|)\|\bz\|^3}{2}
\right]
\E_{\mathbf{0},\bI_k}
\left[
\bu^T\bSigma^{-1/2}\bL_s\bSigma^{-1/2}\bu
\bSigma^{1/2}\bu\bu^T\bSigma^{1/2}
\right].
\end{split}
\]
Furthermore, we can write
\[
\begin{split}
&
\vc(\bSigma^{-1}\bL_j\bSigma^{-1})^T
\vc\left(
\E_{\mathbf{0},\bI_k}
\left[
\bu^T\bSigma^{-1/2}\bL_s\bSigma^{-1/2}\bu
\bSigma^{1/2}\bu\bu^T\bSigma^{1/2}
\right]
\right)\\
&=
\vc(\bL_j)^T
\left(\bSigma^{-1}\otimes\bSigma^{-1}\right)
\E_{\mathbf{0},\bI_k}
\left[
\vc\left(
\bSigma^{1/2}\bu\bu^T\bSigma^{1/2}
\right)
\bu^T\bSigma^{-1/2}\bL_s\bSigma^{-1/2}\bu
\right]\\
&=
\vc(\bL_j)^T
\left(\bSigma^{-1}\otimes\bSigma^{-1}\right)
\left(\bSigma^{1/2}\otimes\bSigma^{1/2}\right)
\E_{\mathbf{0},\bI_k}
\left[
\vc\left(
\bu\bu^T
\right)
\vc\left(\bu\bu^T\right)^T
\right]
\vc\left(\bSigma^{-1/2}\bL_s\bSigma^{-1/2}\right)\\
&=
\vc\left(\bSigma^{-1/2}\bL_j\bSigma^{-1/2}\right)^T
\frac{1}{k(k+2)}
\Big(
\bI_{k^2}+\bK_{k,k}
+
\vc(\bI_k)\vc(\bI_k)^T
\Big)
\vc\left(\bSigma^{-1/2}\bL_s\bSigma^{-1/2}\right),
\end{split}
\]
using Lemma~\ref{lem:Lemma 5.1}.
Application of property~\eqref{eq:prop K} and the fact that $\vc(\bA)^T\vc(\bB)=\text{tr}(\bA\bB)$, yields
\[
\begin{split}
&
\vc(\bSigma^{-1}\bL_j\bSigma^{-1})^T
\vc\left(
\E_{\mathbf{0},\bI_k}
\left[
\bu^T\bSigma^{-1/2}\bL_s\bSigma^{-1/2}\bu
\bSigma^{1/2}\bu\bu^T\bSigma^{1/2}
\right]
\right)\\
&=
\frac{1}{k(k+2)}
\left(
2\text{tr}(\bSigma^{-1}\bL_j\bSigma^{-1}\bL_s)
+
\text{tr}(\bSigma^{-1}\bL_j)\text{tr}(\bSigma^{-1}\bL_s)
\right).
\end{split}
\]
It follows that the first term on the right hand side of~\eqref{eq:derivative LambdaV-thetaj wrt thetas}
leads to a first term in $\partial\Lambda_{\btheta,j}(\bxi_P)/\partial \theta_s$, which is equal to
\begin{equation}
\label{eq:first term}
\frac{\E_{\mathbf{0},\bI_k}
\left[
u'(\|\bz\|)\|\bz\|^3
\right]}{2(k+2)}
\Big(
2\text{tr}(\bSigma^{-1}\bL_j\bSigma^{-1}\bL_s)+\text{tr}(\bSigma^{-1}\bL_j)\text{tr}(\bSigma^{-1}\bL_s)
\Big).
\end{equation}
The second term on the right hand side of~\eqref{eq:derivative LambdaV-thetaj wrt thetas} is equal to
\[
\begin{split}
&
\E_{\mathbf{0},\bI_k}
\left[
\frac{v'(\|\bz\|)}{2\|\bz\|}
\bz^T\bSigma^{-1/2}\bL_s\bSigma^{-1/2}\bz \bSigma
\right]\\
&=
\E_{\mathbf{0},\bI_k}
\left[
\frac{v'(\|\bz\|)\|\bz\|}{2}
\right]
\E_{\mathbf{0},\bI_k}
\left[
\bu^T\bSigma^{-1/2}\bL_s\bSigma^{-1/2}\bu
\right]
\bSigma\\
&=
\E_{\mathbf{0},\bI_k}
\left[
\frac{v'(\|\bz\|)\|\bz\|}{2}
\right]
\vc\left(\bSigma^{-1/2}\bL_s\bSigma^{-1/2}\right)^T
\vc\left(\E_{\mathbf{0},\bI_k}
\left[
\bu\bu^T
\right]
\right)
\bSigma\\
&=
\E_{\mathbf{0},\bI_k}
\left[
\frac{v'(\|\bz\|)\|\bz\|}{2}
\right]
\vc\left(\bSigma^{-1/2}\bL_s\bSigma^{-1/2}\right)^T
\vc\left(
\frac1k\bI_k
\right)
\bSigma\\
&=
\frac{\E_{\mathbf{0},\bI_k}
\left[
v'(\|\bz\|)\|\bz\|
\right]}{2k}
\text{tr}\left(\bSigma^{-1}\bL_s\right)
\bSigma,\\
\end{split}
\]
using Lemma~\ref{lem:Lemma 5.1}.
This leads to a second term in $\partial\Lambda_{\btheta,j}(\bxi_P)/\partial \theta_s$, which is equal to
\begin{equation}
\label{eq:second term}
\begin{split}
&-
\frac{\E_{\mathbf{0},\bI_k}
\left[
v'(\|\bz\|)\|\bz\|
\right]}{2k}
\text{tr}\left(\bSigma^{-1}\bL_s\right)
\vc\left(\bSigma^{-1}\bL_j\bSigma^{-1}\right)^T
\vc(\bSigma)\\
&=
-
\frac{\E_{\mathbf{0},\bI_k}
\left[
v'(\|\bz\|)\|\bz\|
\right]}{2k}
\text{tr}\left(\bSigma^{-1}\bL_s\right)
\text{tr}\left(\bSigma^{-1}\bL_j\right).
\end{split}
\end{equation}
The third term on the right hand side of~\eqref{eq:derivative LambdaV-thetaj wrt thetas} leads to a third term
in $\partial\Lambda_{\btheta,j}(\bxi_P)/\partial \theta_s$, which is equal to
\begin{equation}
\label{eq:third term}
\E_{\mathbf{0},\bI_k}
\left[
v(\|\bz\|)
\right]
\vc\left(\bSigma^{-1}\bL_j\bSigma^{-1}\right)^T
\vc(\bL_s)
=
\E_{\mathbf{0},\bI_k}
\left[
v(\|\bz\|)
\right]
\text{tr}\left(\bSigma^{-1}\bL_j\bSigma^{-1}\bL_s\right).
\end{equation}
We conclude that $\partial\Lambda_{\btheta,j}(\bxi_P)/\partial \theta_s$ consists of three terms
given in~\eqref{eq:first term}, \eqref{eq:second term} and~\eqref{eq:third term}.
This means that $\partial\Lambda_{\btheta,j}(\bxi_P)/\partial \theta_s$
has a term $\text{tr}(\bSigma^{-1}\bL_j\bSigma^{-1}\bL_s)$
with coefficient
\[
\frac{\E_{\mathbf{0},\bI_k}
\left[
u'(\|\bz\|)\|\bz\|^3
\right]}{(k+2)}
+
\E_{\mathbf{0},\bI_k}
\left[
v(\|\bz\|)
\right]
=
\gamma_1,
\]
and a term $\text{tr}(\bSigma^{-1}\bL_s)\text{tr}(\bSigma^{-1}\bL_j)$ with coefficient
\[
\frac{\E_{\mathbf{0},\bI_k}
\left[
u'(\|\bz\|)\|\bz\|^3
\right]}{2(k+2)}
-
\frac{\E_{\mathbf{0},\bI_k}
\left[
v'(\|\bz\|)\|\bz\|
\right]}{2k}
=
-\gamma_2,
\]
where $\gamma_1$ and $\gamma_2$ are defined in~\eqref{def:gamma12},
and where we use that $u'(s)s^3=\rho''(s)s^2-\rho'(s)s$ and $v(s)=\rho'(s)s-\rho(s)+b_0$.
Finally, from the definition of $\bL$ in~\eqref{def:L} it follows that the $l\times l$ matrix with entries
\[
\begin{split}
&
\gamma_1\text{tr}(\bSigma^{-1}\bL_j\bSigma^{-1}\bL_s)-\gamma_2\text{tr}(\bSigma^{-1}\bL_s)\text{tr}(\bSigma^{-1}\bL_j)\\
&=
\gamma_1
\vc(\bL_j)^T
\left(\bSigma^{-1/2}\otimes\bSigma^{-1/2}\right)
\left(\bSigma^{-1/2}\otimes\bSigma^{-1/2}\right)
\vc(\bL_s)\\
&\qquad-
\gamma_2
\vc(\bL_j)^T
\left(\bSigma^{-1/2}\otimes\bSigma^{-1/2}\right)
\vc(\bI_k)
\vc(\bI_k)^T
\left(\bSigma^{-1/2}\otimes\bSigma^{-1/2}\right)
\vc(\bL_s)\\
&=
\gamma_1\mathbf{L}^T\left(\bSigma^{-1}\otimes\bSigma^{-1}\right)\mathbf{L}
-
\gamma_2
\vc(\bL_j)^T
\vc(\bSigma^{-1})
\vc(\bSigma^{-1})^T
\vc(\bL_s),
\end{split}
\]
is the matrix
\[
\gamma_1\mathbf{L}^T\left(\bSigma^{-1}\otimes\bSigma^{-1}\right)\mathbf{L}
-
\gamma_2\mathbf{L}^T
\vc(\bSigma^{-1})
\vc(\bSigma^{-1})^T
\mathbf{L}.
\]
This proves the lemma.
\end{proof}

\begin{lemma}
\label{lem:inverse}
Suppose that $\rho$ satisfies (R3)-(R4).
Let $\gamma_1$ and $\gamma_2$ defined in~\eqref{def:gamma12} and suppose that $\gamma_1>0$.
Then the inverse of $\partial\Lambda_{\btheta}(\bxi(P))/\partial\btheta^T$ exists and is given by
\[
\begin{split}
&
a
\Bigg(
\bL^T\left(\bSigma^{-1}\otimes\bSigma^{-1}\right)\bL
\Bigg)^{-1}\\
&+
b
\Bigg(
\bL^T\left(\bSigma^{-1}\otimes\bSigma^{-1}\right)\bL
\Bigg)^{-1}
\bL^T\vc(\bSigma^{-1})\vc(\bSigma^{-1})^T\bL
\Bigg(
\bL^T\left(\bSigma^{-1}\otimes\bSigma^{-1}\right)\bL
\Bigg)^{-1}
\end{split}
\]
where $a=1/\gamma_1$ and $b=\gamma_2/(\gamma_1(\gamma_1-k\gamma_2))$.
\end{lemma}
\begin{proof}
Together with Lemma~\ref{lem:block derivative}, first write
\[
\begin{split}
\frac{\partial\Lambda_{\btheta}(\bxi(P))}{\partial\btheta^T}
&=
\gamma_1
\bL^T\left(\bSigma^{-1}\otimes\bSigma^{-1}\right)\bL
-
\gamma_2\bL^T
\vc(\bSigma^{-1})
\vc(\bSigma^{-1})^T
\bL\\
&\qquad=
\gamma_1
\bE^T\bE
-
\gamma_2\bE^T
\vc(\bI_k)
\vc(\bI_k)^T
\bE,
\end{split}
\]
where
$\bE=\left(\bSigma^{-1/2}\otimes\bSigma^{-1/2}\right)\bL$.
Since $\bV(\btheta(P))=\bSigma$,
by definition of~$\bL$, it follows that for~$\btheta(P)=(\theta_1,\ldots,\theta_l)^T$,
\[
\bL\,\btheta(P)
=
\sum_{j=1}^l
\theta_j\vc(\bL_j)
=
\vc\left(\sum_{j=1}^l
\theta_j\bL_j
\right)
=
\vc(\bSigma).
\]
This means that
\begin{equation}
\label{eq:E from theta to V}
\bE\,\btheta(P)
=
\left(\bSigma^{-1/2}\otimes\bSigma^{-1/2}\right)\bL\,
\btheta(P)
=
\left(\bSigma^{-1/2}\otimes\bSigma^{-1/2}\right)
\vc(\bSigma)
=
\vc(\bI_k).
\end{equation}
Since $\bL$ has full rank, also $\bE$ has full rank.
This means that $(\bE^T\bE)^{-1}$ exists, and  satisfies
\begin{equation}
\label{eq:E from V to theta}
(\bE^T\bE)^{-1}\bE^T\vc(\bI_k)=\btheta(P).
\end{equation}
Now, write
\[
\begin{split}
&
\gamma_1
\bE^T\bE
-
\gamma_2\bE^T
\vc(\bI_k)
\vc(\bI_k)^T
\bE\\
&=
\Bigg(
\gamma_1
\bI_k
-
\gamma_2\bE^T
\vc(\bI_k)
\vc(\bI_k)^T
\bE(\bE^T\bE)^{-1}
\Bigg)
(\bE^T\bE).
\end{split}
\]
When we multiply the first matrix from the left with a matrix of the same type,
\begin{equation}
\label{eq:mult from left}
\begin{split}
&
\Bigg(
a
\bI_k
+
b
\bE^T
\vc(\bI_k)
\vc(\bI_k)^T
\bE(\bE^T\bE)^{-1}
\Bigg)\\
&\qquad\qquad\qquad\qquad
\Bigg(
\gamma_1
\bI_k
-
\gamma_2\bE^T
\vc(\bI_k)
\vc(\bI_k)^T
\bE(\bE^T\bE)^{-1}
\Bigg),
\end{split}
\end{equation}
then we find four terms.
A term $\bI_k$ with coefficient $a\gamma_1$, two terms
\[
\bE^T
\vc(\bI_k)
\vc(\bI_k)^T
\bE(\bE^T\bE)^{-1},
\]
with coefficient $-a\gamma_2+b\gamma_1$, and the term
\[
\bE^T
\vc(\bI_k)
\vc(\bI_k)^T
\bE(\bE^T\bE)^{-1}
\bE^T
\vc(\bI_k)
\vc(\bI_k)^T
\bE(\bE^T\bE)^{-1}
\]
with coefficient $-b\gamma_2$.
Consider the scalar valued inner product in the middle
\begin{equation}
\label{eq:inner product}
\vc(\bI_k)^T
\bE(\bE^T\bE)^{-1}
\bE^T
\vc(\bI_k)
=
\vc(\bI_k)^T
\bE\,\btheta(P)
=
\vc(\bI_k)^T\vc(\bI_k)=k,
\end{equation}
by application of~\eqref{eq:E from V to theta} and then~\eqref{eq:E from theta to V}.
It follows that the term with coefficient $-b\gamma_2$ reduces to
\[
k\bE^T
\vc(\bI_k)
\vc(\bI_k)^T
\bE(\bE^T\bE)^{-1}.
\]
Hence the matrix product in~\eqref{eq:mult from left} is equal to
\begin{equation}
\label{eq:result matrix prod}
a\gamma_1\bI_k
+
(-a\gamma_2+b\gamma_1-kb\gamma_2)\bE^T
\vc(\bI_k)
\vc(\bI_k)^T
\bE(\bE^T\bE)^{-1}
\end{equation}
When we multiply the same matrix from the right,
\[
\Bigg(
\gamma_1
\bI_k
-
\gamma_2\bE^T
\vc(\bI_k)
\vc(\bI_k)^T
\bE(\bE^T\bE)^{-1}
\Bigg)
\Bigg(
a
\bI_k
+
b
\bE^T
\vc(\bI_k)
\vc(\bI_k)^T
\bE(\bE^T\bE)^{-1}
\Bigg),
\]
we find the same result~\eqref{eq:result matrix prod}.
This matrix is equal to $\bI_k$ if and only if
$a\gamma_1=1$ and $-a\gamma_2+b\gamma_1-kb\gamma_2=0$,
or equivalently
\begin{equation}
\label{def:a en b}
\begin{split}
a&=1/\gamma_1\\
b&=\frac{a\gamma_2}{\gamma_1-k\gamma_2}
=
\frac{\gamma_2}{\gamma_1(\gamma_1-k\gamma_2)},
\end{split}
\end{equation}
where we use that $\gamma_1>0$ and
$\gamma_1-k\gamma_2=\mathbb{E}_{0,\bI_k}\left[\rho'(\|\bz\|)\right]/2>0$,
due to (R3)-(R4).
We conclude that the inverse of
$\gamma_1
\bI_k
-
\gamma_2\bE^T
\vc(\bI_k)
\vc(\bI_k)^T
\bE(\bE^T\bE)^{-1}$
exists and is equal to
\[
a\bI_k
+
b\bE^T
\vc(\bI_k)
\vc(\bI_k)^T
\bE(\bE^T\bE)^{-1}
\]
with $a$ and $b$ given in~\eqref{def:a en b}.
Hence, the inverse of the matrix
$\gamma_1
\bE^T\bE
-
\gamma_2\bE^T
\vc(\bI_k)
\vc(\bI_k)^T
\bE$
is equal to
\[
\begin{split}
&
\Big(\bE^T\bE\Big)^{-1}
\Bigg(
\gamma_1
\bI_k
-
\gamma_2\bE^T
\vc(\bI_k)
\vc(\bI_k)^T
\bE(\bE^T\bE)^{-1}
\Bigg)^{-1}\\
&=
\Big(\bE^T\bE\Big)^{-1}
\Bigg(
a
\bI_k
+b
\bE^T
\vc(\bI_k)
\vc(\bI_k)^T
\bE(\bE^T\bE)^{-1}
\Bigg)\\
&=
a
(\bE^T\bE)^{-1}
+b
(\bE^T\bE)^{-1}
\bE^T
\vc(\bI_k)
\vc(\bI_k)^T
\bE(\bE^T\bE)^{-1}.
\end{split}
\]
After inserting $\bE=\left(\bSigma^{-1/2}\otimes\bSigma^{-1/2}\right)\bL$,
this finishes the proof.
\end{proof}
\paragraph*{Proof of Corollary~\ref{cor:IF elliptical}}
\begin{proof}
Since $\rho$ is strictly increasing on~$[0,c_0]$,
the function $u(s)=\rho'(s)/s> 0$, for $0<s\leq c_0$ and zero for $s>c_0$.
This means that
\[
\alpha
=
\mathbb{E}_{\mathbf{0},\bI_k}
\left[
\left(1-\frac1{k}\right)
\frac{\rho'(\|\bz\|)}{\|\bz\|}
+
\frac1k
\rho''(\|\bz\|)
\right]
\geq
\frac1k
\mathbb{E}_{\mathbf{0},\bI_k}
\left[
\rho''(\|\bz\|)
\right]
>0.
\]
Furthermore, since $\bX$ has full rank, the inverse of $\E\left[\bX^T\bSigma^{-1}\bX\right]$ exists.
It follows that the matrix $\partial\Lambda_{\bbeta}(\bxi(P))/\partial\bbeta$ in~\eqref{def:derivative Lambda beta}
is non-singular.
According to Lemma~\ref{lem:inverse}, also
$\partial\Lambda_{\btheta}(\bxi(P))/\partial\btheta$ is non-singular.
Together, with Lemma~\ref{lem:diff Lambda} and Lemma~\ref{lem:block derivative}, we conclude that
$\partial\Lambda/\partial\bxi$ is continuously differentiable with a non-singular derivative at $\bxi(P)$,
so that Theorem~\ref{th:IF} applies.
Together with Lemma~\ref{lem:block derivative},
this implies that
\[
\begin{split}
\mathrm{IF}(\bs_0,\bbeta,P)
&=
-\left(\frac{\partial \Lambda_{\bbeta}(\bxi(P))}{\partial \bbeta}\right)^{-1}\Psi_{\bbeta}(\bs_0,\bxi(P))\\
&=
\frac{u(d_0)}{\alpha}
\left(
\E\left[\bX^T\bSigma^{-1}\bX\right]
\right)^{-1}
\bX_0^T\bSigma^{-1}(\by_0-\bX_0\bbeta)
\end{split}
\]
where $d_0^2=(\by_0-\bX_0\bbeta)^T\bSigma^{-1}(\by_0-\bX_0\bbeta)$.
From Theorem~\ref{th:IF}, together with Lemma~\ref{lem:block derivative}, it also follows that
\begin{equation}
\label{eq:decompose IF theta}
\begin{split}
&
\mathrm{IF}(\bs_0,\btheta,P)
=
-\left(\frac{\partial \Lambda_{\btheta}(\bxi(P))}{\partial \btheta}\right)^{-1}\Psi_{\btheta}(\bs_0,\bxi(P))\\
&=
\left(\frac{\partial \Lambda_{\btheta}(\bxi(P))}{\partial \btheta}\right)^{-1}
\bL^T
(\bSigma^{-1}\otimes\bSigma^{-1})\times\\
&\qquad\qquad\qquad
\times\vc\left(
ku(d_0)(\by_0-\bX_0\bbeta)(\by_0-\bX_0\bbeta)^T-v(d_0)\bSigma
\right)\\
&=
ku(d_0)
\left(\frac{\partial \Lambda_{\btheta}(\bxi(P))}{\partial \btheta}\right)^{-1}
\bL^T
(\bSigma^{-1}\otimes\bSigma^{-1})
\vc\left((\by_0-\bX_0\bbeta)(\by_0-\bX_0\bbeta)^T\right)\\
&\qquad-
v(d_0)
\left(\frac{\partial \Lambda_{\btheta}(\bxi(P))}{\partial \btheta}\right)^{-1}
\bL^T
(\bSigma^{-1}\otimes\bSigma^{-1})
\vc(\bSigma).
\end{split}
\end{equation}
Consider the first term on the right hand side of~\eqref{eq:decompose IF theta}.
We have that
\[
\begin{split}
&
(\bSigma^{-1/2}\otimes\bSigma^{-1/2})
\vc\left((\by_0-\bX_0\bbeta)(\by_0-\bX_0\bbeta)^T\right)\\
&\qquad=
\vc\left(\bSigma^{-1/2}(\by_0-\bX_0\bbeta)(\by_0-\bX_0\bbeta)^T\bSigma^{-1/2}\right)
\end{split}
\]
and from Lemma~\ref{lem:inverse},
\begin{equation}
\label{eq:expression in E}
\begin{split}
&
\left(\frac{\partial \Lambda_{\btheta}(\bxi(P))}{\partial \btheta}\right)^{-1}
\bL^T
(\bSigma^{-1/2}\otimes\bSigma^{-1/2})\\
&\quad=
a
(\bE^T\bE)^{-1}
\bE^T
+
b
(\bE^T\bE)^{-1}
\bE^T\vc(\bI_k)\vc(\bI_k)^T\bE
(\bE^T\bE)^{-1}
\bE^T,
\end{split}
\end{equation}
where $\bE=(\bSigma^{-1/2}\otimes\bSigma^{-1/2})\bL$.
This implies
\[
\begin{split}
&
\left(\frac{\partial \Lambda_{\btheta}(\bxi(P))}{\partial \btheta}\right)^{-1}
\bL^T
(\bSigma^{-1}\otimes\bSigma^{-1})
\vc\left((\by_0-\bX_0\bbeta)(\by_0-\bX_0\bbeta)^T\right)\\
&\quad=
a
(\bE^T\bE)^{-1}
\bE^T
\vc\left(\bSigma^{-1/2}(\by_0-\bX_0\bbeta)(\by_0-\bX_0\bbeta)^T\bSigma^{-1/2}\right)\\
&\qquad+
b
(\bE^T\bE)^{-1}
\bE^T\vc(\bI_k)\vc(\bI_k)^T\bE
(\bE^T\bE)^{-1}
\bE^T\times\\
&\qquad\qquad\qquad\times
\vc\left(\bSigma^{-1/2}(\by_0-\bX_0\bbeta)(\by_0-\bX_0\bbeta)^T\bSigma^{-1/2}\right).
\end{split}
\]
The first term on the right hand side is equal to
\[
\begin{split}
&
a\Big(\bL^T(\bSigma^{-1}\otimes\bSigma^{-1})\bL)\Big)^{-1}
\bL^T(\bSigma^{-1/2}\otimes\bSigma^{-1/2})\times\\
&\qquad\times
\vc\left(\bSigma^{-1/2}(\by_0-\bX_0\bbeta)(\by_0-\bX_0\bbeta)^T\bSigma^{-1/2}\right)
\end{split}
\]
and, with~\eqref{eq:E from theta to V} and~\eqref{eq:E from V to theta}, the second term on the right hand side is equal to
\[
\begin{split}
&
b\btheta(P)\btheta(P)^T\bE^T\vc\left(\bSigma^{-1/2}(\by_0-\bX_0\bbeta)(\by_0-\bX_0\bbeta)^T\bSigma^{-1/2}\right)\\
&=
b\btheta(P)\vc(\bI_k)^T\vc\left(\bSigma^{-1/2}(\by_0-\bX_0\bbeta)(\by_0-\bX_0\bbeta)^T\bSigma^{-1/2}\right)\\
&=
b\btheta(P)
\text{tr}\left(\bSigma^{-1/2}(\by_0-\bX_0\bbeta)(\by_0-\bX_0\bbeta)^T\bSigma^{-1/2}\right)\\
&=bd_0^2\btheta(P).
\end{split}
\]
It follows that the first term on the right hand side of~\eqref{eq:decompose IF theta} is equal to
\[
\begin{split}
&
aku(d_0)
\Big(\bL^T(\bSigma^{-1}\otimes\bSigma^{-1})\bL)\Big)^{-1}
\bL^T(\bSigma^{-1/2}\otimes\bSigma^{-1/2})\times\\
&\quad\times
\vc\left(\bSigma^{-1/2}(\by_0-\bX_0\bbeta)(\by_0-\bX_0\bbeta)^T\bSigma^{-1/2}\right)\\
&\qquad+
bku(d_0)d_0^2\btheta(P).
\end{split}
\]
Next consider the second term on the right hand side of~\eqref{eq:decompose IF theta}.
We have that
\[
(\bSigma^{-1/2}\otimes\bSigma^{-1/2})\vc(\bSigma)
=
\vc(\bI_k),
\]
and with~\eqref{eq:expression in E},
together with~\eqref{eq:E from theta to V} and~\eqref{eq:E from V to theta},
\[
\begin{split}
&\left(\frac{\partial \Lambda_{\btheta}(\bxi(P))}{\partial \btheta}\right)^{-1}
\bL^T
(\bSigma^{-1}\otimes\bSigma^{-1})
\vc(\bSigma)\\
&\quad=
a
(\bE^T\bE)^{-1}
\bE^T
\vc\left(\bI_k\right)
+
b
(\bE^T\bE)^{-1}
\bE^T\vc(\bI_k)
\vc(\bI_k)^T\bE
(\bE^T\bE)^{-1}
\bE^T
\vc\left(\bI_k\right)\\
&\quad=
a\btheta(P)
+
b\btheta(P)
\vc(\bI_k)^T\bE
(\bE^T\bE)^{-1}
\bE^T
\vc\left(\bI_k\right)\\
&\quad=
a\btheta(P)
+
b\btheta(P)
\vc(\bI_k)^T\bE
\btheta(P)\\
&\quad=
a\btheta(P)
+
b\btheta(P)
\vc(\bI_k)^T\vc(\bI_k)\\
&\quad=
(a+bk)\btheta(P).
\end{split}
\]
It follows that the second term on the right hand side of~\eqref{eq:decompose IF theta} is equal to
\[
-v(d_0)(a+bk)\btheta(P).
\]
Putting things together, we find that $\mathrm{IF}(\bs_0,\btheta,P)$ is equal to
\[
\begin{split}
&
aku(d_0)
\Big(\bL^T(\bSigma^{-1}\otimes\bSigma^{-1})\bL)\Big)^{-1}
\bL^T(\bSigma^{-1/2}\otimes\bSigma^{-1/2})\times\\
&\quad\times
\vc\left(\bSigma^{-1/2}(\by_0-\bX_0\bbeta)(\by_0-\bX_0\bbeta)^T\bSigma^{-1/2}\right)\\
&\qquad+
(bku(d_0)d_0^2-av(d_0)-bkv(d_0))\btheta(P).
\end{split}
\]
Since $v(d_0)=u(d_0)d_0^2-\rho(d_0)+b_0$, we have that
\[
\begin{split}
bku(d_0)d_0^2-av(d_0)-bkv(d_0)
&=
-\left(
\frac{v(d_0)}{\gamma_1}
-
\frac{k\gamma_2}{\gamma_1(\gamma_1-k\gamma_2)}
(\rho(d_0)-b_0)
\right)\\
&=
-\frac{u(d_0)d_0^2}{\gamma_1}
+
\left(
\frac{1}{\gamma_1}
+
\frac{k\gamma_2}{\gamma_1(\gamma_1-k\gamma_2)}
\right)
(\rho(d_0)-b_0)\\
&=
-\frac{u(d_0)d_0^2}{\gamma_1}
+
\frac{\rho(d_0)-b_0}{\gamma_1-k\gamma_2}\\
\end{split}
\]
We conclude that $\mathrm{IF}(\bs_0,\btheta,P)$ is given by
\[
\begin{split}
&
\frac{ku(d_0)}{\gamma_1}
\Big(\bL^T(\bSigma^{-1}\otimes\bSigma^{-1})\bL)\Big)^{-1}
\bL^T
\vc\left(\bSigma^{-1}(\by_0-\bX_0\bbeta)(\by_0-\bX_0\bbeta)^T\bSigma^{-1}\right)\\
&\qquad+
\left(
-\frac{u(d_0)d_0^2}{\gamma_1}
+
\frac{\rho(d_0)-b_0}{\gamma_1-k\gamma_2}
\right)
\btheta(P).
\end{split}
\]
This proves the corollary.
\end{proof}

\subsection{Proofs of Section~\ref{sec:asymptotic normality}}
As preparation we first establish the following lemma,
which is similar to Lemma~22 in~\cite{nolan&pollard1987}.
\begin{lemma}
\label{lem:nolan-pollard}
Let $\rho(\cdot)$ be a real-valued function of bounded variation on $\R^+$.
The class of all functions on $\R^p$ of the form
\[
\bs=(\by,\bX)\mapsto \rho(\|\bA(\by-\bX\bbeta)\|)
\]
with $\bA$ ranging over all $k\times k$ matrices and $\bbeta$ ranging over $\R^q$,
has polynomial discrimination.
\end{lemma}
\begin{proof}
Consider the class of functions
\[
g_{\bA,\bbeta}(\mathbf{s},t)
=
\|\bA(\by-\bX\bbeta)\|^2-\rho^{-1}(t)^2.
\]
According to Lemma 18 in~\cite{nolan&pollard1987},
it suffices to show that the functions $g_{\bA,\bbeta}(\cdot,\cdot)$ span a finite-dimensional vector space.
In order to do so, write
\[
\begin{split}
\|\bA(\by-\bX\bbeta)\|^2
&=
\sum_{i=1}^k
\sum_{j=1}^k
(\bA^T\bA)_{ij}(\by-\bX\bbeta)_i(\by-\bX\bbeta)_j\\
&=
\sum_{i=1}^k
\sum_{j=1}^k
\sum_{s=1}^k a_{is}a_{sj}
\left(y_i-\sum_{r=1}^q x_{ir}\beta_r\right)
\left(y_j-\sum_{w=1}^q x_{jw}\beta_w\right)
\end{split}
\]
This is a polynomial in $\mathbf{s}=(y_1,\ldots,y_k, x_{11},\ldots,x_{kq})$,
with coefficients in $\R$.
This means that the class of functions
$g_{\bA,\bbeta}(\mathbf{s},t)=
\|\bA(\by-\bX\bbeta)\|^2-\rho^{-1}(t)^2$ forms a finite dimensional vector space.
\end{proof}

A useful first step is the following lemma.
\begin{lemma}
\label{lem:subgraph class}
Let $u:\R\to\R$ be a function of bounded variation.
Let $\bs=(\by,\bX)=(s_1,\ldots,s_p)\in\R^p$, and define
\[
g(\bs,\bbeta,\bV)=u(\|\bV^{-1/2}(\by-\bX\bbeta)\|)
\quad
\text{for }
\bbeta\in\R^q, \bV\in\text{PDS}(k).
\]
Consider the classes of functions
\[
\begin{split}
\mathcal{F}
&=
\{g(\bs,\bbeta,\bV):\bbeta\in\R^q, \bV\in\text{PDS}(k)\},\\
\mathcal{F}_a
&=
\{g(\bs,\bbeta,\bV)s_a:g\in\mathcal{F}\},\\
\mathcal{F}_{ab}
&=
\{g(\bs,\bbeta,\bV)s_as_b:g\in\mathcal{F}\},
\end{split}
\]
for $a,b=1,\ldots,p$.
Denote by $\mathcal{G}$, $\mathcal{G}_a$, and $\mathcal{G}_{ab}$,
the corresponding classes of graphs of the functions in
$\mathcal{F}$, $\mathcal{F}_a$, and $\mathcal{F}_{ab}$, respectively.
Then $\mathcal{G}$, $\mathcal{G}_a$, and $\mathcal{G}_{ab}$, all have polynomial discrimination
for $a,b=1,\ldots,p$.
\end{lemma}
\begin{proof}
Because the function $u(\cdot)$ is of bounded variation, it follows from
Lemma~\ref{lem:nolan-pollard} that the class~$\mathcal{G}$ has polynomial discrimination.
To show the same for the class $\mathcal{G}_a$, suppose that~$\mathcal{G}_a$ is not of polynomial discrimination.
This means that for every integer $N$ there exists a set $V\subset\R^{p+1}$ of $N$ points,
such that all subsets of~$V$ can be written as $D_a\cap V$ for some $D_a\in\mathcal{G}_a$.
Let $V=\{\bv_1,\ldots,\bv_N\}$, where $\bv_i=(\bs_i,t_i)$, for $i=1,\ldots,N$.
Then for every $(\bs,t)\in V$ with $\bs=(s_1,\ldots,s_a,\ldots,s_p)$,
it must hold that $s_a\ne 0$, otherwise this point cannot be separated from the other points
by an element of the class $\mathcal{G}_a$.
Define the set $V_a=\{(\bs,t/s_a):(\bs,t)\in V\}$.
Note that for $D_a\in \mathcal{G}_a$, we have
\[
\begin{split}
(\bs,t)\in D_a
&\quad\Leftrightarrow\quad
0\leq t\leq g(\bs,\bbeta,\bV)s_a\text{ or }g(\bs,\bbeta,\bV)s_a\leq t\leq 0\\
&\quad\Leftrightarrow\quad
0\leq \frac{t}{s_a}\leq g(\bs,\bbeta,\bV)\text{ or }g(\bs,\bbeta,\bV)\leq \frac{t}{s_a}\leq 0\\
&\quad\Leftrightarrow\quad
(\bs,t/s_a)\in D,
\end{split}
\]
where $D\in \mathcal{G}$.
This implies that every subset of $V_a$ can be written as $D\cap V_a$, for some $D\in \mathcal{G}$.
However, this is in contradiction with the fact that $\mathcal{G}$ has polynomial discrimination.
We conclude that also~$\mathcal{G}_a$ has polynomial discrimination.
A similar argument yields that $\mathcal{G}_{ab}$ has polynomial discrimination.
\end{proof}

Lemma~\ref{lem:subgraph class} is comparable to Lemma~3 in~\cite{lopuhaa1997}, for similar classes of functions built from
functions~$g(\by,\bt,\bC)=u(\|\bC^{-1/2}(\by-\bt)\|)$, where $\bt\in\R^k$ and $\bC\in\text{PDS}(k)$.
With Lemma~\ref{lem:subgraph class} we can establish suitable bounds on the third term
of~\eqref{eq:decomposition estimator}.
This is provided by the following key lemma.
Once having established Lemma~\ref{lem:stoch equi}, asymptotic normality can be derived easily from~\eqref{eq:decomposition estimator}.

\begin{lemma}
\label{lem:stoch equi}
Let $\Psi=(\Psi_{\bbeta},\Psi_{\btheta})$ be defined in~\eqref{eq:Psi function}
and let $\bxi_n=\bxi(\mathbb{P}_n)$ and $\bxi_P=\bxi(P)$ be the solutions to
minimization problems~\eqref{def:Smin estimator structured} and~\eqref{def:Smin structured}.
Suppose that $\rho$ satisfies (R1)-(R4), such that $u(s)$ is of bounded variation,
and suppose that $\bV$ satisfies (V4).
Suppose that $\bxi_n\to\bxi_P$, in probability, and that
$\E\|\bs\|^2<\infty$.
Then
\begin{equation}
\label{eq:stoch equi}
\int
\left(
\Psi(\mathbf{s},\bxi_n)-\Psi(\mathbf{s},\bxi_P)
\right)
\,\dd (\mathbb{P}_n-P)(\mathbf{s})
=
o_P(1/\sqrt{n}).
\end{equation}
\end{lemma}
\begin{proof}
First write
$\Psi_{\btheta,j}(\bs,\bxi)
=
\Psi_{2,j}(\bs,\bxi)
-
\Psi_{3,j}(\bs,\bxi)$,
for $j=1,\ldots,l$, where
\begin{equation}
\label{eq:decompose psi-theta}
\begin{split}
\Psi_{2,j}(\bs,\bxi)
&=
u(d)
(\by-\bX\bbeta)^T\bV^{-1}
\bH_j
\bV^{-1}(\by-\bX\bbeta)\\
\Psi_{3,j}(\bs,\bxi)
&=
\mathrm{tr}\left(\bV^{-1}\frac{\partial \bV}{\partial \theta_j}\right)(\rho(d)-b_0),
\end{split}
\end{equation}
where $\bH_j$ and $d=d(\bs,\bxi)$ are defined in~\eqref{def:Hj} and~\eqref{def:Mahalanobis distance structured}.
It suffices to show that
\begin{eqnarray}
\int
\left(
\Psi_{\bbeta}(\bs,\bxi_n)-\Psi_{\bbeta}(\bs,\bxi_P)
\right)
\,\dd (\mathbb{P}_n-P)(\bs)
&=&\label{eq:stoch equi vector}
o_P(1/\sqrt{n}),
\\
\int
\left(
\Psi_{2,j}(\bs,\bxi_n)-\Psi_{2,j}(\bs,\bxi_P)
\right)
\,\dd (\mathbb{P}_n-P)(\bs)
&=&\label{eq:stoch equi matrix}
o_P(1/\sqrt{n}),
\\
\int
\left(
\Psi_{3,j}(\bs,\bxi_n)-\Psi_{3,j}(\bs,\bxi_P)
\right)
\,\dd (\mathbb{P}_n-P)(\bs)
&=&\label{eq:stoch equi scalar}
o_P(1/\sqrt{n}),
\end{eqnarray}
for $j=1,\ldots,l$.

To obtain~\eqref{eq:stoch equi scalar},
first write
$\bM_n=\bV(\btheta_n)^{-1}$ and $\bM_P=\bV(\btheta_P)^{-1}$,
so that $\bM_n\to\bM_P$, in probability,
according to condition~(V1).
Decompose as follows
\begin{equation}
\label{eq:decomp Psi3}
\begin{split}
&
\Psi_{3,j}(\bs,\bxi_n)-\Psi_{3,j}(\bs,\bxi_P)\\
&\qquad=
\mathrm{tr}\left(\bM_n^{-1}\frac{\partial \bV(\btheta_n)}{\partial \theta_j}\right)
(\rho(d(\bs,\bxi_n))-\rho(d(\bs,\bxi_P)))\\
&\qquad\qquad+
\left\{
\mathrm{tr}\left(\bM_n^{-1}\frac{\partial \bV(\btheta_n)}{\partial \theta_j}\right)
-
\mathrm{tr}\left(\bM_P^{-1}\frac{\partial \bV(\btheta_P)}{\partial \theta_j}\right)
\right\}
(\rho(d(\bs,\bxi_P))-b_0).
\end{split}
\end{equation}
After integration with respect to $\mathbb{P}_n-P$, the first term on the right hand side of~\eqref{eq:decomp Psi3}
becomes
\begin{equation}
\label{eq:first term scalar}
\mathrm{tr}\left(\bM_n^{-1}\frac{\partial \bV(\btheta_n)}{\partial \theta_j}\right)
\int
(\rho(d(\bs,\bxi_n))-\rho(d(\bs,\bxi_P)))
\,\dd(\mathbb{P}_n-P)(\bs),
\end{equation}
where with (V4),
\begin{equation}
\label{eq:conv trace}
\mathrm{tr}\left(\bM_n^{-1}\frac{\partial \bV(\btheta_n)}{\partial \theta_j}\right)
\to
\mathrm{tr}\left(\bM_P^{-1}\frac{\partial \bV(\btheta_P)}{\partial \theta_j}\right)
\end{equation}
in probability.
Furthermore, note that all functions
$\rho(d(\cdot,\bxi))$, for $\bxi\in\R^q\times\R^l$ are members of the
class
\[
\mathcal{F}=\left\{\rho(\|\bV^{-1/2}(\by-\bX\bbeta)\|):\bbeta\in\R^q, \bV\in\text{PDS}(k)\right\}.
\]
From (R1)-(R2) it follows that the function $\rho$ is of bounded variation (being the sum of two monotone functions).
According to Lemma~\ref{lem:subgraph class},
the class $\mathcal{G}$, consisting of subgraphs of functions in the class $\mathcal{F}$,
has polynomial discrimination.
Moreover, the class $\mathcal{F}$ has a constant envelope.
Then as a result of the empirical process theory developed in Pollard~\cite{pollard1984}
(e.g., see Theorem~1 in~\cite{lopuhaa1997}), it follows that for every $\delta>0$,
\[
\sup_{f_1,f_2\in[\delta]}
\sqrt{n}
\left|
\int
\left(f_1(\mathbf{s})-f_2(\mathbf{s})\right)\,\dd (\mathbb{P}_n-P)(\mathbf{s})
\right|
\to0
\]
in probability, where
\[
[\delta]=\left\{(f_1,f_2): f_1,f_2\in \mathcal{F}\text{ and } \int(f_1-f_2)^2\,\dd P\leq \delta^2\right\}.
\]
Because $\bxi_n\to\bxi_P$ in probability one, the pair of functions $f_1(\mathbf{s})=\rho(d(\mathbf{s},\bxi_n)$ and
$f_2(\mathbf{s})=\rho(d(\mathbf{s},\bxi_P)$ are in the set~$[\delta]$,
for sufficiently large $n$, with probability tending to one.
It follows that
\begin{equation}
\label{eq:order Psi3}
\begin{split}
&
\sqrt{n}
\left|
\int
\left(
\rho(d(\mathbf{s},\bxi_n)-\rho(d(\mathbf{s},\bxi_P)
\right)
\,\dd (\mathbb{P}_n-P)(\mathbf{\mathbf{s}})
\right|\\
&\quad\leq
\sup_{f_1,f_2\in[\delta]}
\sqrt{n}
\left|
\int
\left(f_1(\mathbf{s})-f_2(\mathbf{s})\right)\,\dd (\mathbb{P}_n-P)(\mathbf{s})
\right|
\to0,
\end{split}
\end{equation}
in probability.
Together with the fact that $\text{tr}(\bM_P^{-1}\partial\bV(\btheta_P)/\partial\theta_j)$ is bounded,
this proves that~\eqref{eq:first term scalar} is of the order $o_P(1/\sqrt{n})$.
For the second term on the right hand side of~\eqref{eq:decomp Psi3}, we have that according to the central limit theorem
\[
\int (\rho(d(\bs,\bxi_P))-b_0)\,\dd(\mathbb{P}_n-P)(\bs)=O_P(1/\sqrt{n}).
\]
Together with~\eqref{eq:conv trace} this implies that,
after integration with respect to $\mathbb{P}_n-P$,
the second term on the right hand side of~\eqref{eq:decomp Psi3} is of the order $o_P(1/\sqrt{n})$.
This proves~\eqref{eq:stoch equi scalar}.

To obtain~\eqref{eq:stoch equi vector},
decompose as follows
\begin{equation}
\label{eq:decomp Psi1}
\begin{split}
&
\Psi_{\bbeta}(\mathbf{s},\bxi_n)-\Psi_{\bbeta}(\mathbf{s},\bxi_P)\\
&\quad=
\bigg\{
u(d(\mathbf{s},\bxi_n))\bX^T\bM_n\by
-
u(d(\mathbf{s},\bxi_P))\bX^T\bM_P\by
\bigg\}\\
&\qquad
-
\bigg\{
u(d(\mathbf{s},\bxi_n))\bX^T\bM_n\bX\bbeta_n
-
u(d(\mathbf{s},\bxi_P))\bX^T\bM_P\bX\bbeta_P
\bigg\}.
\end{split}
\end{equation}
We will treat both terms on the right hand side separately.
For the first term on the right hand side of~\eqref{eq:decomp Psi1},
we write
\begin{equation}
\label{eq:decomp Psi1 term1}
\Big\{u(d(\mathbf{s},\bxi_n))-u(d(\mathbf{s},\bxi_P))\Big\}
\bX^T\bM_n\by
+
u(d(\mathbf{s},\bxi_P)
\bX^T(\bM_n-\bM_P)\by.
\end{equation}
First consider the first term in~\eqref{eq:decomp Psi1 term1}.
We consider each single element of the vector~$\bX^T\bM_n\by\in\R^q$ separately.
For $i=1,\ldots,q$ fixed, write
\[
\left(\bX^T\bM_n\by\right)_i
=
\sum_{j=1}^k x_{ji}\left(\bM_n\by\right)_j
=
\sum_{j=1}^k x_{ji}
\sum_{s=1}^k m_{n,js}y_s
=
\sum_{j=1}^k
\sum_{s=1}^k
x_{ji}y_s
m_{n,js}.
\]
Hence, after integration with respect to $\mathbb{P}_n-P$,
the $i$-th component of the first term in~\eqref{eq:decomp Psi1 term1}
can be written as a finite sum with summands
\[
\int
\Big(u(d(\mathbf{s},\bxi_n))-u(d(\mathbf{s},\bxi_P))\Big)
x_{ji}y_s
\,\dd (\mathbb{P}_n-P)(\mathbf{s})
m_{n,js},
\]
for $i=1,\ldots,q$ and $j,s\in\{1,\ldots,k\}$ fixed, where $m_{n,js}\to m_{P,js}$,
Because $u$ is of bounded variation and since
$\bs=(\by,\bX)=(y_1,\ldots,y_k,x_{11},\ldots,x_{kq})$,
all functions
$u(d(\bs,\bxi))x_{ji}y_s$, for $\bxi\in\R^q\times\R^l$, are members of the class
\[
\mathcal{F}_{ab}
=
\left\{
u(\|\bV^{-1/2}(\by-\bX\bbeta)\|)s_as_b:\bbeta\in\R^q, \bV\in\text{PDS}(k)
\right\}.
\]
According to Lemma~\ref{lem:subgraph class}, the corresponding class of subgraphs has polynomial discrimination,
so similar to~\eqref{eq:order Psi3} it follows that
for each $i=1,\ldots,q$ and $j,s\in\{1,\ldots,k\}$ fixed,
\[
\int
\Big(u(d(\mathbf{s},\bxi_n))-u(d(\mathbf{s},\bxi_P))\Big)
x_{ji}y_s
\,\dd (\mathbb{P}_n-P)(\mathbf{s})
=
o_P(n^{-1/2}),
\]
which means that
\begin{equation}
\label{eq:order Psi1 term1 term1}
\int
\Big\{u(d(\mathbf{s},\bxi_n))-u(d(\mathbf{s},\bxi_P))\Big\}
\bX^T\bM_P\by
\,\dd (\mathbb{P}_n-P)(\mathbf{s})
=
o_P(n^{-1/2}).
\end{equation}
Next, consider the second term on the right hand side of~\eqref{eq:decomp Psi1 term1}.
Similar to the first term,
after integration with respect to $\mathbb{P}_n-P$,
the $i$-th component of the second term on the right hand side of~\eqref{eq:decomp Psi1 term1}
can be written as a finite sum
of summands
\[
\int
u(d(\mathbf{s},\bxi_P))
x_{ji}y_s
\,\dd (\mathbb{P}_n-P)(\mathbf{s})
(m_{n,js}-m_{P,js}),
\]
for $i=1,\ldots,q$, and $j,s\in\{1,\ldots,k\}$ fixed.
According to the central limit theorem, the integral is of the order~$O_P(n^{-1/2})$,
and since $m_{n,js}\to m_{P,js}$, in probability,
 it follows that the product is of the order~$o_P(n^{-1/2})$.
We conclude that
\begin{equation}
\label{eq:order Psi1 term1 term2}
\int
u(d(\mathbf{s},\bxi_P))
\bX^T(\bM_n-\bM_P)\by
\,\dd (\mathbb{P}_n-P)(\mathbf{s})
=
o_P(n^{-1/2}).
\end{equation}
Putting together~\eqref{eq:order Psi1 term1 term1} and~\eqref{eq:order Psi1 term1 term2},
it follows for the first term on the right hand side of~\eqref{eq:decomp Psi1} that
\begin{equation}
\label{eq:order Psi1 term1}
\int
\Big\{
u(d(\mathbf{s},\bxi_n))\bX^T\bM_n\by
-
u(d(\mathbf{s},\bxi_P))\bX^T\bM_P\by
\Big\}
\,\dd (\mathbb{P}_n-P)(\mathbf{s})
=
o_P(n^{-1/2}).
\end{equation}
For the second term on the right hand side of~\eqref{eq:decomp Psi1},
we write
\begin{equation}
\label{eq:decomp Psi1 term2}
\begin{split}
&
u(d(\mathbf{s},\bxi_n))\bX^T\bM_n\bX\bbeta_n
-
u(d(\mathbf{s},\bxi_P))\bX^T\bM_P\bX\bbeta_P\\
&=
\Big(
u(d(\mathbf{s},\bxi_n))-u(d(\mathbf{s},\bxi_P))
\Big)
\bX^T\bM_n\bX\bbeta_n\\
&\qquad+
u(d(\mathbf{s},\bxi_P)
\Big(
\bX^T\bM_n\bX\bbeta_n-\bX^T\bM_P\bX\bbeta_P
\Big).
\end{split}
\end{equation}
Consider the first term on the right hand side of~\eqref{eq:decomp Psi1 term2}.
For $i=1,\ldots,q$ fixed, write
\[
\begin{split}
\left(
\bX^T\bM_n\bX\bbeta_n
\right)_i
&=
\sum_{j=1}^k
x_{ji}
(\bM_n\bX\bbeta_n)_j
=
\sum_{j=1}^k
x_{ji}
\sum_{s=1}^k
m_{n,js}
(\bX\bbeta_n)_s\\
&=
\sum_{j=1}^k
x_{ji}
\sum_{s=1}^k
m_{n,js}
\sum_{t=1}^k
x_{st}\bbeta_{n,t}
=
\sum_{j=1}^k\sum_{s=1}^k\sum_{t=1}^k
x_{ji}x_{st}m_{n,js}\bbeta_{n,t}.
\end{split}
\]
We see that, after integration with respect to $\mathbb{P}_n-P$, the $i$-th component of the first term
on the right hand side of~\eqref{eq:decomp Psi1 term2} can be written as a finite summation of summands
\[
\int
\Big(u(d(\mathbf{s},\bxi_n))-u(d(\mathbf{s},\bxi_P))\Big)
x_{ji}x_{st}
\,\dd (\mathbb{P}_n-P)(\mathbf{s})
m_{n,js}\bbeta_{n,t},
\]
for $i=1,\ldots,q$ and $s,t,j\in\{1,\ldots,k\}$, where $m_{n,js}\to m_{P,js}$ and
$\bbeta_{n,t}\to\bbeta_{P,t}$, in probability.
All functions
$u(d(\bs,\bxi))x_{ji}y_s$, for $\bxi\in\R^q\times\R^l$, are members of the class
\[
\mathcal{F}_{ab}
=
\left\{
u(\|\bV^{-1/2}(\by-\bX\bbeta)\|)s_as_b:\bbeta\in\R^q, \bV\in\text{PDS}(k)
\right\},
\]
and $u$ is of bounded variation.
According to Lemma~\ref{lem:subgraph class}, the corresponding class of subgraphs has polynomial discrimination,
so similar to~\eqref{eq:order Psi3} it follows that
\[
\int
\Big(u(d(\mathbf{s},\bxi_n))-u(d(\mathbf{s},\bxi_P))\Big)
x_{ji}x_{st}
\,\dd (\mathbb{P}_n-P)(\mathbf{s})
=
o_P(n^{-1/2}),
\]
which means that
\begin{equation}
\label{eq:order Psi1 term2 term1}
\int
\Big(
u(d(\mathbf{s},\bxi_n))-u(d(\mathbf{s},\bxi_P))
\Big)
\bX^T\bM_n\bX\bbeta_n
\,\dd (\mathbb{P}_n-P)(\mathbf{s})
=
o_P(n^{-1/2}).
\end{equation}
Next, consider the second term on the right hand side of~\eqref{eq:decomp Psi1 term2}.
We then have to deal with summands of the form
\[
\int
u(d(\mathbf{s},\bxi_P))
x_{ji}x_{st}
\,\dd (\mathbb{P}_n-P)(\mathbf{s})
(m_{n,js}\bbeta_{n,t}-m_{P,js}\bbeta_{P,t}),
\]
for $i=1,\ldots,q$, and $s,t,j\in\{1,\ldots,k\}$, where $m_{n,js}\to m_{P,js}$ and
$\bbeta_{n,t}\to\bbeta_{P,t}$, in probability.
According to the central limit theorem, the integral is of the order~$O_P(n^{-1/2})$.
Because, $m_{n,is}\bbeta_{n,t}\to m_{P,is}\bbeta_{P,t}$, in probability,
it follows that the product is of the order~$o_P(n^{-1/2})$.
We conclude,
\begin{equation}
\label{eq:order Psi1 term2}
\int
u(d(\mathbf{s},\bxi_P)
\Big(
\bX^T\bM_n\bX\bbeta_n-\bX^T\bM_P\bX\bbeta_P
\Big)
\,\dd (\mathbb{P}_n-P)(\mathbf{s})
=
o_P(n^{-1/2}).
\end{equation}
Putting together~\eqref{eq:order Psi1 term1} and~\eqref{eq:order Psi1 term2},
proves~\eqref{eq:stoch equi vector}

Finally, consider $\Psi_{2,j}$ in~\eqref{eq:decompose psi-theta}, with $\bH_j$ defined~\eqref{def:Hj}.
Write
\[
\begin{split}
\bM_n&=\bV(\btheta_n)^{-1}\bH_j(\btheta_n)\bV(\btheta_n)^{-1}\\
\bM_P&=\bV(\btheta_P)^{-1}\bH_j(\btheta_P)\bV(\btheta_P)^{-1},
\end{split}
\]
so that $\bM_n\to \bM_P$, in probability,
according to condition~(V4).
Decompose~$\Psi_{2,j}(\bs,\bxi_n)-\Psi_{2,j}(\bs,\bxi_P)$ as follows
\begin{equation}
\label{eq:decomp Psi2 term12}
\begin{split}
&
\Big\{
u(d(\bs,\bxi_n))-u(d(\bs,\bxi_P))
\Big\}
(\by-\bX\bbeta_n)^T\bM_n(\by-\bX\bbeta_n)\\
&\qquad+
u(d(\bs,\bxi_P))
\Big\{
(\by-\bX\bbeta_n)^T\bM_n(\by-\bX\bbeta_n)
-
(\by-\bX\bbeta_P)^T\bM_P(\by-\bX\bbeta_P)
\Big\}.
\end{split}
\end{equation}
The first term in~\eqref{eq:decomp Psi2 term12},
can be written as the trace of the matrix
\[
\Big\{
u(d(\bs,\bxi_n))-u(d(\bs,\bxi_P))
\Big\}
(\by-\bX\bbeta_n)(\by-\bX\bbeta_n)^T\bM_n,
\]
where $\bbeta_n\to\bbeta_P$ and $\bM_n\to \bM_P$, in probability.
As before, we consider each single entry of this $k\times k$ matrix.
The $(i,j)$-th element of $(\by-\bX\bbeta_n)(\by-\bX\bbeta_n)^T$ is equal to
\begin{equation}
\label{eq:decomp1}
\begin{split}
(\by-\bX\bbeta_n)_i(\by-\bX\bbeta_n)_j
=
y_iy_j
&-
\sum_{s=1}^k
y_jx_{is}\beta_{n,s}
-
\sum_{t=1}^k
y_ix_{jt}\beta_{n,t}\\
&+
\sum_{s=1}^k
\sum_{t=1}^k
x_{is}x_{jt}\beta_{n,s}\beta_{n,t}.
\end{split}
\end{equation}
We see that
the $(i,j)$-th entry of
\[
\Big\{
u(d(\bs,\bxi_n))-u(d(\bs,\bxi_P))
\Big\}
(\by-\bX\bbeta_n)(\by-\bX\bbeta_n)^T
\]
is a combination of four summations.
The last of these summations arising from~\eqref{eq:decomp1},
after integration with respect to $\mathbb{P}_n-P$, has summands
\[
\int
\Big(
u(d(\bs,\bxi_n))-u(d(\bs,\bxi_P))
\Big)
x_{is}x_{jt}
\,\dd (\mathbb{P}_n-P)(\bs)\,
\beta_{n,s}\beta_{n,t},
\]
where $\beta_{n,s}\beta_{n,t}\to\beta_{P,s}\beta_{P,t}$, in probability.
All functions
$u(d(\bs,\bxi))x_{ji}y_s$, for $\bxi\in\R^q\times\R^l$, are members of the class
\[
\mathcal{F}_{ab}
=
\left\{
u(\|\bV^{-1/2}(\by-\bX\bbeta)\|)s_as_b:\bbeta\in\R^q, \bV\in\text{PDS}(k)
\right\},
\]
where $u$ is of bounded variation
According to Lemma~\ref{lem:subgraph class}, the corresponding class of subgraphs has polynomial discrimination,
so similar to~\eqref{eq:order Psi3} it follows that
\[
\int
\Big\{u(d(\bs,\bxi_n))-u(d(\bs,\bxi_P))\Big\}
x_{is}x_{jt}
\,\dd (\mathbb{P}_n-P)(\bs)
=
o_P(n^{-1/2}).
\]
The other three summations that arise from the right hand side of~\eqref{eq:decomp1} can be handled in the same way.
It follows, that for each ~$i,j\in\{1,\ldots,k\}$ fixed,
\[
\int
\Big\{
u(d(\bs,\bxi_n))-u(d(\bs,\bxi_P))
\Big\}
(\by-\bX\bbeta_n)_i(\by-\bX\bbeta_n)_j
\,\dd
(\mathbb{P}_n-P)(\bs)
=
o_P(n^{-1/2}),
\]
which means that
\[
\int
\Big\{
u(d(\bs,\bxi_n))-u(d(\bs,\bxi_P))
\Big\}
(\by-\bX\bbeta_n)(\by-\bX\bbeta_n)^T
\bM_n
\,\dd
(\mathbb{P}_n-P)(\bs)
=
o_P(n^{-1/2}).
\]
After taking traces, we conclude that
\begin{equation}
\label{eq:order Psi2 term1}
\int
\Big\{
u(d(\bs,\bxi_n))-u(d(\bs,\bxi_P))
\Big\}
(\by-\bX\bbeta_n)^T\bM_n(\by-\bX\bbeta_n)
\,\dd
(\mathbb{P}_n-P)(\bs)
=
o_P(n^{-1/2}).
\end{equation}
Next, consider the second term on the right hand side of~\eqref{eq:decomp Psi2 term12}.
First, note that
\[
\begin{split}
(\by-\bX\bbeta)^T\bM(\by-\bX\bbeta)
&=
\sum_{i=1}^k
\sum_{j=1}^k
y_iy_jm_{ij}
+
\sum_{i=1}^k
\sum_{j=1}^k
y_i
\sum_{s=1}^k
x_{js}\beta_s
m_{ij}\\
&+
\sum_{i=1}^k
\sum_{j=1}^k
y_j
\sum_{t=1}^k
x_{it}\beta_t
m_{ij}
+
\sum_{i=1}^k
\sum_{j=1}^k
\sum_{s=1}^k
\sum_{t=1}^k
x_{is}x_{jt}\beta_s\beta_tm_{ij}.
\end{split}
\]
This means that
\begin{equation}
\label{eq:decomp2}
\begin{split}
&
(\by-\bX\bbeta_n)^T\bM_n(\by-\bX\bbeta_n)
-
(\by-\bX\bbeta_P)^T\bM_P(\by-\bX\bbeta_P)\\
&\quad=
\sum_{i=1}^k
\sum_{j=1}^k
y_iy_j
\left(
m_{n,ij}-m_{P,ij}
\right)
+
\sum_{i=1}^k
\sum_{j=1}^k
\sum_{s=1}^k
y_i
x_{js}
\left(
\beta_{n,s}m_{n,ij}
-
\beta_{P,s}m_{P,ij}
\right)\\
&\qquad+
\sum_{i=1}^k
\sum_{j=1}^k
\sum_{t=1}^k
y_j
x_{it}
\left(
\beta_{n,t}m_{n,ij}
-
\beta_{P,t}m_{P,ij}
\right)\\
&\qquad+
\sum_{i=1}^k
\sum_{j=1}^k
\sum_{s=1}^k
\sum_{t=1}^k
x_{is}x_{jt}
\left(
\beta_{n,s}\beta_{n,t}m_{n,ij}
-
\beta_{P,s}\beta_{P,t}m_{P,ij}
\right).
\end{split}
\end{equation}
We see that the $(i,j)$-th entry of
\[
u(d(\bs,\bxi_P))
\Big\{
(\by-\bX\bbeta_n)^T\bM_n(\by-\bX\bbeta_n)
-
(\by-\bX\bbeta_P)^T\bM_P(\by-\bX\bbeta_P)
\Big\}
\]
can be written as the combination of four summations.
The last of the summations arising from~\eqref{eq:decomp2},
after integration with respect to $\mathbb{P}_n-P$, has summands
\[
\int
u(d(\bs,\bxi_P))
x_{is}x_{jt}
\,\dd (\mathbb{P}_n-P)(\bs)
\left(
\beta_{n,s}\beta_{n,t}m_{n,ij}
-
\beta_{P,s}\beta_{P,t}m_{P,ij}
\right),
\]
where $\beta_{n,s}\beta_{n,t}m_{n,ij}\to\beta_{P,s}\beta_{P,t}m_{P,ij}$, in probability.
According to the central limit theorem, the integral is of the order
$O_P(n^{-1/2})$, whereas the second term tends to zero.
All functions
$u(d(\bs,\bxi))x_{is}x_{jt}$, for $\bxi\in\R^q\times\R^l$, are members of the class
\[
\mathcal{F}_{ab}
=
\left\{
u(\|\bV^{-1/2}(\by-\bX\bbeta)\|)s_as_b:\bbeta\in\R^q, \bV\in\text{PDS}(k)
\right\},
\]
where $u$ is of bounded variation
According to Lemma~\ref{lem:subgraph class}, the corresponding class of subgraphs has polynomial discrimination,
so similar to~\eqref{eq:order Psi3} it follows that
\[
\int
u(d(\bs,\bxi_P))
x_{is}x_{jt}
\,\dd (\mathbb{P}_n-P)(\bs)
\left(
\beta_{n,s}\beta_{n,t}m_{n,ij}
-
\beta_{P,s}\beta_{P,t}m_{P,ij}
\right)
=
o_P(n^{-1/2}).
\]
The other three summations that arise from the right hand side of~\eqref{eq:decomp2} can be handled in the same way, so that
\begin{equation}
\label{eq:order Psi2 term2}
\begin{split}
&
\int
u(d(\bs,\bxi_P))
\Big\{
(\by-\bX\bbeta_n)^T\bM_n(\by-\bX\bbeta_n)
-
(\by-\bX\bbeta_P)^T\bM_P(\by-\bX\bbeta_P)
\Big\}
\,\dd (\mathbb{P}_n-P)(\bs)\\
&\quad
=
o_P(n^{-1/2}).
\end{split}
\end{equation}
Putting together~\eqref{eq:order Psi2 term1} and~\eqref{eq:order Psi2 term2},
proves~\eqref{eq:stoch equi matrix} for each $j=1,\ldots,l$.
This finishes the proof of Lemma~\ref{lem:stoch equi}.
\end{proof}
\paragraph*{Proof of Corollary~\ref{cor:Asymp norm elliptical}}
\begin{proof}
As in the proof of Corollary~\ref{cor:IF elliptical}, it follows that
$\partial\Lambda/\partial\bxi$ is continuously differentiable with a non-singular derivative at $\bxi(P)$,
so that Theorem~\ref{th:asymp normal} applies.
According to Theorem~\ref{th:asymp normal} and Lemma~\ref{lem:block derivative},
it follows that $\sqrt{n}(\bbeta_n-\bbeta(P))$ is asymptotically normal with mean zero
and covariance matrix
\[
\frac1{\alpha^2}
\left(
\mathbb{E}\left[\mathbf{X}^T\bSigma^{-1}\mathbf{X}\right]
\right)^{-1}
\E\left[\Psi_{\bbeta}(\bs,\bxi_P)\Psi_{\bbeta}(\bs,\bxi_P)^T\right]
\left(
\mathbb{E}\left[\mathbf{X}^T\bSigma^{-1}\mathbf{X}\right]
\right)^{-1}
\]
where $\Psi_{\bbeta}$ is defined in~\eqref{eq:Psi linear}.
We find that
\[
\begin{split}
\E\left[\Psi_{\bbeta}(\bs,\bxi_P)\Psi_{\bbeta}(\bs,\bxi_P)^T\right]
&=
\E\left[
\bX^T
\E\left[
u(d)^2
\bSigma^{-1}
(\by-\bmu)(\by-\bmu)^T
\bSigma^{-1}
\Big|
\bX
\right]
\bX
\right],
\end{split}
\]
where $d^2=(\by-\bmu)^T\bSigma^{-1}(\by-\bmu)$ and
$u(s)=\rho'(s)/s$.
As before, with $\bz=\bSigma^{-1/2}(\by-\bmu)$ and $\bu=\bz/\|\bz\|$,
according to Lemma~\ref{lem:Lemma 5.1},
the inner conditional expectation can be written as
\[
\bSigma^{-1/2}
\E_{\mathbf{0},\bI_k}\left[u(\|\bz\|)^2\|\bz\|^2\right]
\E_{\mathbf{0},\bI_k}
\left[\bu\bu^T\right]
\bSigma^{-1/2}
=
\frac{\E_{\mathbf{0},\bI_k}\left[u(\|\bz\|)^2\|\bz\|^2\right]}{k}
\bSigma^{-1}.
\]
This implies that the asymptotic covariance of $\sqrt{n}(\bbeta_n-\bbeta(P))$ is given by
\[
\frac{\E_{\mathbf{0},\bI_k}\left[\rho'(\|\bz\|)^2\right]}{k\alpha^2}
\left(
\mathbb{E}\left[\mathbf{X}^T\bSigma^{-1}\mathbf{X}\right]
\right)^{-1}.
\]
Again, according to Theorem~\ref{th:asymp normal} and Lemma~\ref{lem:block derivative},
it follows that $\sqrt{n}(\btheta_n-\btheta(P))$ is asymptotically normal with mean zero
and covariance matrix
\[
\left(
\frac{\partial\Lambda_{\btheta}(\bxi(P))}{\partial \btheta}
\right)^{-1}
\E\left[
\vc\left(\Psi_{\btheta}(\bs,\bxi_P)\right)
\vc\left(\Psi_{\btheta}(\bs,\bxi_P)\right)^T
\right]
\left(
\frac{\partial\Lambda_{\btheta}(\bxi(P))}{\partial \btheta}
\right)^{-1}
\]
where $\Psi_{\btheta}$ is defined in~\eqref{eq:Psi linear}.
We have
\[
\begin{split}
&
\E\left[
\vc\left(\Psi_{\btheta}(\bs,\bxi_P)\right)
\vc\left(\Psi_{\btheta}(\bs,\bxi_P)\right)^T
\right]\\
&=
\bE^T
\E\left[
\vc\left(\bSigma^{-1/2}\Psi_{\bV}(\bs,\bxi_P)\bSigma^{-1/2}\right)
\vc\left(\bSigma^{-1/2}\Psi_{\bV}(\bs,\bxi_P)\bSigma^{-1/2}\right)^T
\right]
\bE
\end{split}
\]
where $\Psi_\bV$ is defined in~\eqref{def:PsiV}
and
$\bE=\left(\bSigma^{-1/2}\otimes\bSigma^{-1/2}\right)\bL$
and
\[
\begin{split}
&\E\left[
\vc\left(\bSigma^{-1/2}\Psi_{\btheta}(\bs,\bxi_P)\bSigma^{-1/2}\right)
\vc\left(\bSigma^{-1/2}\Psi_{\btheta}(\bs,\bxi_P)\bSigma^{-1/2}\right)^T
\right]\\
&=
k^2\E_{\mathbf{0},\bI_k}\left[u(\|\bz\|)^2\|\bz\|^4\right]
\E_{\mathbf{0},\bI_k}\left[\vc\left(\bu\bu^T\right)\vc\left(\bu\bu^T\right)^T\right]\\
&\quad-
k\E_{\mathbf{0},\bI_k}\left[u(\|\bz\|)v(\|\bz\|)\|\bz\|^2\right]
\E_{\mathbf{0},\bI_k}\left[\vc\left(\bu\bu^T\right)\vc\left(\bI_k\right)^T\right]\\
&\quad-
k\E_{\mathbf{0},\bI_k}\left[u(\|\bz\|)v(\|\bz\|)\|\bz\|^2\right]
\E_{\mathbf{0},\bI_k}\left[\vc\left(\bI_k\right)\vc\left(\bu\bu^T\right)^T\right]\\
&\quad+
\E_{\mathbf{0},\bI_k}\left[v(\|\bz\|)^2\right]
\E_{\mathbf{0},\bI_k}\left[\vc\left(\bI_k\right)\vc\left(\bI_k\right)^T\right].
\end{split}
\]
From Lemma~\ref{lem:Lemma 5.1}, the first term on the right hand side is equal to
\[
\frac{k\E_{\mathbf{0},\bI_k}\left[u(\|\bz\|)^2\|\bz\|^4\right]}{k+2}
\left(
\bI_{k^2}+\mathbf{K}_{k,k}+\vc(\bI_k)\vc(\bI_k)^T
\right).
\]
This leads to one term $\bI_{k^2}+\mathbf{K}_{k,k}$ with coefficient
\[
\frac{k\E_{\mathbf{0},\bI_k}\left[u(\|\bz\|)^2\|\bz\|^4\right]}{k+2}
\]
and using that, according to Lemma~\ref{lem:Lemma 5.1}, $\E_{\mathbf{0},\bI_k}\left[\bu\bu^T\right]=(1/k)\bI_k$,
we find a second term $\vc(\bI_k)\vc(\bI_k)^T$ with coefficient
\[
\frac{k\E_{\mathbf{0},\bI_k}\left[u(\|\bz\|)^2\|\bz\|^4\right]}{k+2}
-2\E_{\mathbf{0},\bI_k}\left[u(\|\bz\|)v(\|\bz\|)\|\bz\|^2\right]
+
\E_{\mathbf{0},\bI_k}\left[v(\|\bz\|)^2\right].
\]
Since $v(s)=u(s)s^2-\rho(s)+b_0$, we have that
\[
\frac{k}{k+2}u(s)^2s^4
-
2u(s)v(s)s^2
+
v(s)^2
=
-\frac{2}{k+2}u(s)^2s^4
+
(\rho(s)-b_0)^2.
\]
This means that
\[
\begin{split}
&\E\left[
\vc\left(\bSigma^{-1/2}\Psi_{\btheta}(\bs,\bxi_P)\bSigma^{-1/2}\right)
\vc\left(\bSigma^{-1/2}\Psi_{\btheta}(\bs,\bxi_P)\bSigma^{-1/2}\right)^T
\right]\\
&=
\delta_1
\left(\bI_{k^2}+\mathbf{K}_{k,k}\right)
+
\delta_2
\vc(\bI_k)\vc(\bI_k)^T
\end{split}
\]
where
\begin{equation}
\label{def:delta12}
\begin{split}
\delta_1
&=
\frac{k\E_{\mathbf{0},\bI_k}\left[u(\|\bz\|)^2\|\bz\|^4\right]}{k+2}\\
\delta_2
&=
-\frac{2}{k+2}
\E_{\mathbf{0},\bI_k}\left[u(\|\bz\|)^2\|\bz\|^4\right]
+
\E_{\mathbf{0},\bI_k}\left[
\left(\rho(\|\bz\|)-b_0\right)^2
\right]\\
&=
-\frac2k\delta_1
+
\E_{\mathbf{0},\bI_k}\left[
\left(\rho(\|\bz\|)-b_0\right)^2
\right]
\end{split}
\end{equation}
Since, $\bK_{k,k}(\bSigma^{-1/2}\otimes\bSigma^{-1/2})=(\bSigma^{-1/2}\otimes\bSigma^{-1/2})\bK_{k,k}$,
together with~\eqref{eq:prop K},
this implies that
\[
\begin{split}
\E\left[
\vc\left(\Psi_{\btheta}(\bs,\bxi_P)\right)
\vc\left(\Psi_{\btheta}(\bs,\bxi_P)\right)^T
\right]
&=
\bE^T
\left(
\delta_1
\left(\bI_{k^2}+\mathbf{K}_{k,k}\right)
+
\delta_2
\vc(\bI_k)\vc(\bI_k)^T
\right)
\bE\\
&=
2\delta_1\bE^T\bE
+
\delta_2
\bE^T\vc(\bI_k)\vc(\bI_k)^T\bE.
\end{split}
\]
Furthermore, according to Lemma~\ref{lem:inverse},
\[
\left(
\frac{\partial\Lambda_{\btheta}(\bxi(P))}{\partial \btheta}
\right)^{-1}
=
a
(\bE^T\bE)^{-1}
+b
(\bE^T\bE)^{-1}
\bE^T
\vc(\bI_k)
\vc(\bI_k)^T
\bE(\bE^T\bE)^{-1},
\]
with $a$ and $b$ defined in~\eqref{def:a en b}.
By application of~\eqref{eq:E from theta to V}, \eqref{eq:E from V to theta},
and~\eqref{eq:inner product}, we find that
\[
\begin{split}
&
\left(
\frac{\partial\Lambda_{\btheta}(\bxi(P))}{\partial \btheta}
\right)^{-1}
\E\left[
\vc\left(\Psi_{\btheta}(\bs,\bxi_P)\right)
\vc\left(\Psi_{\btheta}(\bs,\bxi_P)\right)^T
\right]\\
&=
2a\delta_1\bI_k
+
(a\delta_2+2b\delta_1+bk\delta_2)
(\bE^T\bE)^{-1}\bE^T\vc(\bI_k)
\vc(\bI_k)^T
\bE
\end{split}
\]
and
\[
\begin{split}
&
\left(
\frac{\partial\Lambda_{\btheta}(\bxi(P))}{\partial \btheta}
\right)^{-1}
\E\left[
\vc\left(\Psi_{\btheta}(\bs,\bxi_P)\right)
\vc\left(\Psi_{\btheta}(\bs,\bxi_P)\right)^T
\right]
\left(
\frac{\partial\Lambda_{\btheta}(\bxi(P))}{\partial \btheta}
\right)^{-1}
\\
&=
2\sigma_1(\bE^T\bE)^{-1}
+
\sigma_2(\bE^T\bE)^{-1}\bE^T\vc(\bI_k)
\vc(\bI_k)^T
\bE(\bE^T\bE)^{-1}\\
&=
2\sigma_1(\bE^T\bE)^{-1}
+
\sigma_2\btheta(P)\btheta(P)^T\\
\end{split}
\]
where
\[
\begin{split}
\sigma_1
&=
a^2\delta_1\\
\sigma_2
&=
2b(2a+kb)\delta_1+(a+kb)^2\delta_2.
\end{split}
\]
When we insert the expressions for $\delta_1$, $\delta_2$, $a$, and~$b$ given in~\eqref{def:delta12} and~\eqref{def:a en b},
then we find
\[
\begin{split}
\sigma_1
&=
\frac{k(k+2)\E_{\mathbf{0},\bI_k}\left[u(\|\bz\|)^2\|\bz\|^4\right]}{
\left(
\mathbb{E}_{0,\mathbf{I}_k}
\left[
\rho''(\|\bz\|)\|\bz\|^2+(k+1)\rho'(\|\bz\|)\|\bz\|
\right]
\right)^2}\\
\sigma_2
&=
-\frac2k\sigma_1
+
\frac{4\E_{\mathbf{0},\bI_k}[\left(\rho(\|\bz\|)-b_0\right)^2]}{
\left(\E_{\mathbf{0},\bI_k}\left[\rho'(\|\bz\|)\|\bz\|\right]\right)^2}
\end{split}
\]
By substituting $\bE=\left(\bSigma^{-1/2}\otimes\bSigma^{-1/2}\right)\bL$,
we find that the limiting covariance of $\sqrt{n}(\btheta_n-\btheta(P))$ is given by
\[
2\sigma_1
\Big(\bL^T\left(\bSigma^{-1}\otimes\bSigma^{-1}\right)\bL\Big)^{-1}
+
\sigma_2
\btheta(P)\btheta(P)^T
\]
This finishes the proof.
\end{proof}

\end{document}